\documentclass[10pt, reqno]{amsart}
\usepackage{graphicx, amssymb, amsmath, amsthm}
\numberwithin{equation}{section}

\let\Re=\undefined\DeclareMathOperator*{\Re}{Re}
\let\Im=\undefined\DeclareMathOperator*{\Im}{Im}
\let\LHS=\undefined\DeclareMathOperator*{\LHS}{LHS}

\newcommand{\R}{\mathbb{R}}
\newcommand{\C}{\mathbb{C}}
\newcommand{\wh}[1]{\widehat{#1}}
\newcommand{\wt}[1]{\widetilde{#1}}

\newcommand{\norm}[1]{\| #1\|}
\newcommand{\bnorm}[1]{\bigg\|#1\bigg\|}
\newcommand{\eps}{\varepsilon}

\def\smallint{\begingroup\textstyle \int\endgroup}

\newcommand{\xonorm}[3]{\|#1\|_{L_t^{#2}L_x^{#3}}}
\newcommand{\xnorm}[4]{\|#1\|_{L_t^{#2}L_x^{#3}(#4\times\R^4)}}
\newcommand{\bxnorm}[4]{\big\|#1\big\|_{L_t^{#2}L_x^{#3}(#4\times\R^4)}}
\newcommand{\xonorms}[2]{\|#1\|_{L_{t,x}^{#2}}}

\newcommand{\xnorms}[3]{\|#1\|_{L_{t,x}^{#2}(#3\times\R^4)}}

\newcommand{\nsc}{\vert\nabla\vert^{s_c}}

\let\al=\alpha

\newtheorem{theorem}{Theorem}[section]

\newtheorem{lemma}[theorem]{Lemma}

\newtheorem{corollary}[theorem]{Corollary}

\newtheorem{proposition}[theorem]{Proposition}

\theoremstyle{definition}
\newtheorem{definition}[theorem]{Definition}
\newtheorem{remark}[theorem]{Remark}

\theoremstyle{remark}

\makeatletter
\newcommand{\Extend}[5]{\ext@arrow0099{\arrowfill@#1#2#3}{#4}{#5}}
\makeatother

\begin{document}
\title[Defocusing energy-supercritical NLS]{The defocusing energy-supercritical NLS \\ in four space dimensions}

\author{Changxing Miao}
\address{Institute of Applied Physics and Computational Mathematics,
P. O. Box 8009,\ Beijing,\ China,\ 100088} \email{miao\_changxing@iapcm.ac.cn}

\author{Jason Murphy}
\address{Department of Mathematics, 
UCLA, Los Angeles, CA 90095-1555, USA} \email{murphy@math.ucla.edu}

\author{Jiqiang Zheng}
\address{The Graduate School of China Academy of Engineering Physics, P. O. Box 2101, Beijing, China, 100088}
\email{zhengjiqiang@gmail.com}

\begin{abstract}
We consider a class of defocusing energy-supercritical nonlinear Schr\"odinger equations in four space dimensions. Following a concentration-compactness approach, we show that for $1<s_c<3/2$, any solution that remains bounded in the critical Sobolev space $\dot{H}_x^{s_c}(\R^4)$ must be global and scatter. Key ingredients in the proof include a long-time Strichartz estimate and a frequency-localized interaction Morawetz inequality.  
\end{abstract}

 \maketitle

 \tableofcontents 
 
\section{Introduction}
\noindent We study the initial-value problem for defocusing nonlinear
 Schr\"odinger equations of the form
\begin{align} \label{equ1.1}
\begin{cases}    (i\partial_t+\Delta)u= F(u) \\
u(0,x)=u_0(x),
\end{cases}
\end{align}
where $u:\R_t\times\R_x^4\to \C$ and $F(u)=|u|^pu$ for some $p>0$.

The class of solutions to \eqref{equ1.1} is left invariant by the scaling
\begin{equation}\label{scale}
u(t,x)\mapsto \lambda^{2/p}u(\lambda^2t, \lambda x),\quad\lambda>0.
\end{equation}
This scaling defines a notion of \emph{criticality} for \eqref{equ1.1}. In particular, one can check that the only homogeneous $L_x^2$-based Sobolev space that is left invariant under \eqref{scale} is $\dot{H}_x^{s_c}(\R^4)$, where the \emph{critical regularity} $s_c$ is given by $s_c:=2-2/p$. If we take $u_0\in\dot{H}_x^s(\R^4)$, then for $s=s_c$ we call the problem \eqref{equ1.1} \emph{critical}. For $s>s_c$ we call the problem \emph{subcritical}, while for $s<s_c$ we call the problem \emph{supercritical}. 

We study the critical problem for \eqref{equ1.1} in the \emph{energy-supercritical} regime, that is, $s_c>1$. In four space dimensions, this corresponds to choosing $p>2$. We prove that for $1<s_c<3/2$ (i.e. $2<p<4$), any maximal-lifespan solution that remains uniformly bounded in $\dot{H}_x^{s_c}(\R^4)$ must be global and scatter. 

To begin, we need a few definitions. 

\begin{definition}[Solution]\label{def1.1}
 A function $u:I\times\R^4\to\mathbb{C}$ on a non-empty
time interval $I\ni 0$ is a \emph{solution} to
\eqref{equ1.1} if
it belongs to $C_t\dot{H}_x^{s_c}(K\times\R^4)\cap L_{t,x}^{3p}(K\times\R^4)$
for any compact interval $K\subset I$ and obeys the Duhamel formula
\begin{equation}\label{duhamel.1}
u(t)=e^{it\Delta}u_0-i\int_{0}^te^{i(t-s)\Delta}F(u(s))\,ds
\end{equation}
for each $t\in I$. We call $I$
the \emph{lifespan} of $u$. We say that $u$ is a \emph{maximal-lifespan solution}
if it cannot be extended to any strictly larger interval.
We call $u$ \emph{global} if $I=\R.$
\end{definition}

\begin{definition}[Scattering size and blowup]\label{def1.2}
For a solution $u:I\times\R^4\to\C$ to \eqref{equ1.1}, we define the \emph{scattering size} of $u$ on $I$ by
			\begin{equation}\nonumber
 	S_I(u):=\int_I\int_{\R^4}\vert u(t,x)\vert^{3p}\,dx\,dt.
				\end{equation}

If there exists $t_0\in I$ such that $S_{[t_0,\sup I)}(u)=\infty$, we say that $u$ \emph{blows up forward in time}. Similarly, if there exists $t_0\in I$ such that $S_{(\inf I,t_0]}(u)=\infty$, we say that $u$ \emph{blows up backward in time}. In particular, a solution may blow up in infinite time.	
		
On the other hand, standard arguments show that if $u$ is a global solution to \eqref{equ1.1} that obeys $S_{\R}(u)<\infty$, then $u$ \emph{scatters}, that is, there exist unique $u_{\pm}\in\dot{H}^{s_c}_x(\R^4)$ such that 
		$$\lim_{t\to\pm\infty}\norm{u(t)-e^{it\Delta}u_{\pm}}_{\dot{H}_x^{s_c}(\R^4)}=0.$$
\end{definition}

Our main result is the following

\begin{theorem}\label{theorem}
Let $1<s_c<3/2$. Suppose
$u:~I\times\R^4\to\C$ is a maximal-lifespan solution to
\eqref{equ1.1} such that
$u\in L_t^\infty\dot{H}_x^{s_c}(I\times\R^4)$. Then $u$ is global and scatters, with
	$$S_\R(u)\leq C(\norm{u}_{L_t^\infty\dot{H}_x^{s_c}})$$
for some function $C:[0,\infty)\to[0,\infty).$
\end{theorem}

Equivalently, Theorem \ref{theorem} states that failure to scatter must be accompanied by the divergence of the $\dot{H}_x^{s_c}$-norm. 

The motivation for Theorem \ref{theorem} originates in the study of the mass- and energy-critical nonlinear Schr\"odinger equations. Recall that for an arbitrary space dimension $d$, the equation
		$$(i\partial_t+\Delta)u=\pm\vert u\vert^pu$$
is called \emph{mass-critical} if $p=4/d$ and \emph{energy-critical} if $p=4/(d-2)$. In the mass-critical case, the rescaling \eqref{scale} leaves invariant the \emph{mass} of solutions, which is defined by 
		$$M[u(t)]=\int_{\R^d}\vert u(t,x)\vert^2\,dx$$
and is a conserved quantity for \eqref{equ1.1}. In the energy-critical case, the rescaling \eqref{scale} leaves invariant the \emph{energy} of solutions, which is defined by
		$$E[u(t)]=\int_{\R^d}\tfrac12\vert\nabla u(t,x)\vert^2\pm\tfrac{1}{p+2}\vert u(t,x)\vert^{p+2}\,dx$$
and is again a conserved quantity for \eqref{equ1.1}. The mass-critical NLS corresponds to $s_c=0$, while the energy-critical NLS corresponds to $s_c=1$. 

Due to the presence of conserved quantities at the critical regularity, the mass- and energy-critical equations have been the most widely studied instances of NLS. In this paper, we will apply some of the techniques that have been developed to study the mass- and energy-critical problems to the energy-supercritical regime. The assumption that $u$ stays bounded in the critical Sobolev space $\dot{H}_x^{s_c}$ in Theorem~\ref{theorem} will play the role of the `missing conservation law' at the critical regularity. 

In order to provide some context for our result and to introduce some of the techniques we will use, let us now briefly discuss some previous results for NLS at critical regularity.

For the defocusing energy-critical NLS, it is now known that arbitrary data in $\dot{H}_x^1$ lead to solutions that are global and scatter. This was proven first for radial initial data by Bourgain \cite{Bo99a}, Grillakis \cite{Grillakis}, and Tao \cite{TaoRadial}, and later for arbitrary data by Colliander--Keel--Staffilani--Takaoka--Tao \cite{CKSTT07}, Ryckman--Vi\c{s}an \cite{RV}, and Vi\c{s}an \cite{Visanphd, Visan2007}. (For results in the focusing case, see \cite{KM, KV20101}.) The chief difficulty in establishing these results stems from the fact that none of the known monotonicity formulas (i.e. Morawetz estimates) for NLS scale like the energy. Bourgain's `induction on energy' technique paved the way for how to proceed in such a scenario: by finding a bubble of concentration inside a solution, one can introduce a characteristic length scale into the problem, thus bringing the available Morawetz estimates back into play (despite their non-critical scaling).

We will follow the concentration-compactness approach to induction on energy, which entails the analysis of so-called minimal counterexamples. Minimal counterexamples were originally introduced in the context of the mass-critical NLS (see \cite{BeVa, Bourg, CarKer, Ker0, Ker, MerVeg}), although the first application of minimal counterexamples to establish a global well-posedness result was carried out in the focusing energy-critical setting by Kenig--Merle \cite{KM}.

For the defocusing mass-critical NLS, it has also been established that arbitrary data in $L_x^2$ lead to solutions that are global and scatter. This was proven through the use of minimal counterexamples, first for radial data in dimensions $d\geq 2$ (see \cite{KTV2009, KVZ2008, TVZ2007}) and later for arbitrary data in all dimensions by Dodson \cite{Dodson3, Dodson2, Dodson1}. (For results in the focusing case, see \cite{Dodson, KVZ2008, TVZ2007}.) 

Killip--Vi\c{s}an \cite{KV3D} and Vi\c{s}an \cite{Visan2011} have also revisited the defocusing energy-critical problem in dimensions $d\in\{3,4\}$ from the perspective of minimal counterexamples, incorporating techniques developed by Dodson \cite{Dodson3} in the mass-critical setting. Specifically, they prove a `long-time Strichartz estimate' for almost periodic solutions, which can then be used to preclude the existence of frequency-cascade solutions, as well as to establish a frequency-localized interaction Morawetz inequality (which may in turn be used to preclude the existence of soliton-like solutions). 

There has also been some work done on NLS at non-conserved critical regularity (that is, $s_c\notin\{0,1\}$). So far, it is not known how to treat the large-data case without some \emph{a priori} control of a critical norm. It is natural to conjecture that the analogue of Theorem~\ref{theorem} should hold for any $s_c> 0$ and in any dimension. The first result in this direction is due to Kenig--Merle \cite{KM2010}, who studied the $\dot{H}_x^{1/2}$-critical NLS in three dimensions. Using their concentration-compactness technique (as in \cite{KM}), along with the Lin--Strauss Morawetz inequality (which is scaling-critical in this case), they proved that any solution that stays bounded in $\dot{H}_x^{1/2}$ must be global and scatter. The same result was established in higher dimensions (by a similar approach) in \cite{Mu2}. Recently, the second author \cite{Mu} handled some other cases in the \emph{inter-critical} regime ($0<s_c<1$) in dimensions $d\in\{3,4,5\}$ by making use of the `long-time Strichartz estimate' approach described above.

Some work has also been done in the energy-supercritical regime $(s_c>1)$. In particular, Killip--Vi\c{s}an \cite{KV2010} handled the case of a cubic nonlinearity in dimensions $d\geq 5$, as well as some other cases for which $s_c>1$ in dimensions $d\geq 5$. Their restriction to high dimensions stems from their use of a `double Duhamel trick', which they use to prove that global almost periodic solutions belong to $H_x^1$. Once it is known that solutions belong to $H_x^1$, one can show that frequency-cascade solutions must have zero mass, while the interaction Morawetz inequality can be used directly to rule out soliton-like solutions (that is, no frequency localization is necessary). 

We pause here briefly to mention that similar problems have also been studied for the nonlinear wave equation. The interested reader may refer to \cite{Bulut1, Bulut2, Bulut3, DKM, KMNLW, KV, KVrad, Shen1, Shen2}.

In this paper, we treat the energy-supercritical regime in dimension $d=4$. We cannot employ the strategy of \cite{KV2010} described above, as the double Duhamel trick fails in dimensions $d<5.$ Our argument will instead be more in the spirit of \cite{Dodson3, KV3D, Mu, Visan2011}; that is, we will establish a long-time Strichartz estimate (Theorem~\ref{lse}) and a frequency-localized interaction Morawetz inequality (Theorem~\ref{uppbd}). We will also use arguments from \cite{KV2010, KV20101} to establish an `additional decay' result for almost periodic solutions (Proposition \ref{adecay1}), which we will use in the proof of the interaction Morawetz inequality. Our main results will apply to the range $1<s_c<3/2$. We will encounter this upper bound on $s_c$ both in the proof of the long-time Strichartz estimate (see Remark~\ref{restriction 1}) and in the proof of the interaction Morawetz inequality (see Section~\ref{flim section}).

Let us turn now to an outline of the arguments we will use to establish Theorem~\ref{theorem}.

\subsection{Outline of the proof of Theorem~\ref{theorem}} Before we can address the global-in-time theory for \eqref{equ1.1}, we need to have a good local-in-time theory in place. In particular, we have the following

\begin{theorem}[Local well-posedness]\label{lwp}
Let $1<s_c<3/2.$   Given $u_0\in \dot H_x^{s_c}(\R^4)$ and $t_0\in\R$, there exists a unique maximal-lifespan solution $u:I\times\R^4\to\C$ to \eqref{equ1.1} with $u(t_0)=u_0$. Moreover, this solution satisfies the following:
\begin{enumerate}
\item $($Local existence$)$ $I$ is an open neighborhood of $t_0$.
\item $($Blowup criterion$)$ If\ $\sup I$ is finite, then $u$ blows up
forward in time $($in the sense of Definition \ref{def1.2}$)$. Similarly, if
$\inf I$ is finite, then $u$ blows up backward in time.
\item $($Scattering$)$ If $\sup I=+\infty$ and $u$ does not blow up
forward in time, then $u$ scatters forward in time, that is, there exists unique $u_+\in\dot{H}_x^{s_c}(\R^4)$ such that
	\begin{equation}\label{1.2.1}
	\lim_{t\to\infty}\norm{u(t)-e^{it\Delta}u_+}_{\dot{H}_x^{s_c}(\R^4)}=0.
	\end{equation}
Conversely, given $u_+\in \dot H_x^{s_c}(\R^4)$, there
is a unique solution to \eqref{equ1.1} in a neighborhood of $t=\infty$
so that \eqref{1.2.1} holds. Analogous statements hold backward in time.
\item $($Small-data global existence and scattering$)$ If
$\|u_0\|_{\dot H_x^{s_c}}$ is sufficiently small, then $u$ is global and scatters, with $S_\R(u)\lesssim\norm{u_0}_{\dot{H}_x^{s_c}}^{3p}.$
\end{enumerate}
\end{theorem}
Theorem~\ref{lwp} follows from well-known arguments. In particular, one can first use arguments from \cite{CaW} to establish this theorem for data in the inhomogeneous Sobolev space $H_x^{s_c}$ (in fact, this is carried out in detail in \cite{KV2010}). To remove the $L_x^2$ assumption, one can use a stability result such as the following: 
\begin{theorem}[Stability]\label{pertu123}  Let $1<s_c<3/2$ and let $I$ be a compact time interval. Suppose $\tilde{u}:I\times\R^4\to\C$ is an approximate solution to \eqref{equ1.1} in the sense that 
\begin{align*}
(i\partial_t+\Delta)\tilde{u}=F(\tilde{u})+e
\end{align*}
for some function $e$. Assume that 
		\begin{align*}
		\norm{\tilde{u}}_{L_t^\infty \dot{H}_x^{s_c}(I\times\R^4)}&\leq E,
		\\ S_I(\tilde{u})&\leq L
		\end{align*}
for some $E>0$ and $L>0$.

Let $t_0\in I$ and $u_0\in\dot{H}_x^{s_c}(\R^4).$ Then there exists $\eps_0=\eps_0(E,L)$ such that if
		\begin{align*}
		\norm{u_0-\tilde{u}(t_0)}_{\dot{H}_x^{s_c}(\R^4)}&\leq\eps,
		\\ \norm{\nsc e}_{N^0(I)}&\leq\eps
		\end{align*}
for some $0<\eps<\eps_0,$ then there exists a solution $u:I\times\R^4\to\C$ to \eqref{equ1.1} with $u(t_0)=u_0$ satisfying 
		$$S_I(u-\tilde{u})\lesssim_{E,L} \eps.$$
For the $N^0(I)$ notation, see Definition \ref{def1} below.
\end{theorem}

Theorem~\ref{pertu123} also follows from well-known arguments, which are in fact similar in spirit to the arguments used to prove local well-posedness. It is worth noting that there are cases of NLS for which the stability theory can become quite delicate. In particular, this is the case for small power nonlinearities ($p<1$). One can refer to \cite[Section~3.4]{KVnote} for a further discussion and references. In our setting, however, we have $2<p<4$, and so the proof of Theorem~\ref{pertu123} is fairly straightforward. In particular, one can find most of the necessary ideas in \cite[Theorem~3.3]{KV2010}, wherein a cubic nonlinearity is considered in all dimensions $d\geq 2$. There is one small adjustment needed to deal with non-polynomial nonlinearities, but this technology exists as well. In particular, one can make use of \cite[Lemma~2.3]{KVrad} to estimate fractional derivatives of differences of non-polynomial nonlinearities. For an example of such an argument in the context of NLS, see for example \cite[Theorem~3.4]{Mu}. 

With the local theory in place, we are now in a position to sketch the proof of Theorem~\ref{theorem}. 

We argue by contradiction and suppose that Theorem~\ref{theorem} fails. Noting that Theorem~\ref{lwp} implies global existence and scattering for sufficiently small initial data, we can deduce the existence of a critical threshold size, below which the theorem holds but above which we can find solutions with arbitrarily large scattering size. Using a limiting argument, we can then deduce the existence of minimal counterexamples, that is, blowup solutions that live exactly at the critical threshold. As a consequence of their minimality, these solutions can be shown to possess compactness properties that are ultimately at odds with the dispersive nature of the equation. Hence, we can eventually preclude their existence altogether and conclude that Theorem~\ref{theorem} holds.

The key property of these minimal counterexamples is that of almost periodicity modulo the symmetries of the equation. Let us briefly discuss this property and some of its immediate consequences; for a more comprehensive treatment, one should refer to \cite{KVnote}.

\begin{definition}[Almost periodic solutions]\label{AP}
Let $s_c>0.$ A solution $u:I\times\R^4\to\C$ to \eqref{equ1.1} is called \emph{almost periodic }(\emph{modulo symmetries}) if 
			\begin{equation}
			\label{assume1.1}
			u\in L_t^\infty\dot{H}_x^{s_c}(I\times\R^4)
			\end{equation} 
and there exist functions $N:~I\to\R^+,$ $x:I\to\R^4$ and $C:\R^+\to\R^+$ such that
for all $t\in I$ and $\eta>0$,
\begin{equation}\label{apss}
\int_{|x-x(t)|>\frac{C(\eta)}{N(t)}}\big||\nabla|^{s_c}
u(t,x)\big|^2\,dx+\int_{|\xi|> C(\eta)N(t)}|\xi|^{2s_c}\vert\wh
u(t,\xi)\vert^2\,d\xi\leq\eta.
\end{equation}

We call $N(t)$ the \emph{frequency scale function}, $x(t)$ the \emph{spatial center function}, and
$C(\eta)$ the \emph{compactness modulus function}. 
\end{definition}

\begin{remark}\label{rem1.1}
The Arzel\`a--Ascoli theorem tells us that a family of functions $\mathcal{F}$ is precompact in $\dot{H}_x^{s_c}(\R^4)$ if and only if it is norm-bounded and there exists a compactness modulus function $C(\eta)$ such that
	$$\int_{\vert x\vert> C(\eta)}\big\vert\nsc f(x)\big\vert^2\,dx+\int_{\vert\xi\vert> C(\eta)}\vert\xi\vert^{2s_c}\vert \wh{f}(\xi)\vert^2\,d\xi\leq\eta$$
uniformly for $f\in\mathcal{F}$. Thus we see that a solution $u:I\times\R^4\to\C$ is almost periodic if and only if
	$$\{u(t):t\in I\}\subset\{\lambda^{2/p}f(\lambda(x+x_0)):\lambda\in(0,\infty),\ x_0\in\R^4,\ \text{and }f\in K\}$$
for some compact $K\subset \dot{H}_x^{s_c}(\R^4).$ We deduce the following: 

First, there exists a function $c:\R^+\to\R^+$ such that
\begin{equation}\label{xiaoc}
\int_{|x-x(t)|\leq\frac{c(\eta)}{N(t)}}\big||\nabla|^{s_c}
u(t,x)\big|^2\,dx+\int_{|\xi|\leq c(\eta)N(t)}|\xi|^{2s_c} |\wh{u}(t,\xi)|^2\,d\xi\leq\eta
\end{equation}
for all $t\in I$.

Second, using the Sobolev embedding $\dot{H}_x^{s_c}(\R^4)\hookrightarrow L_x^{2p}(\R^4)$, we can see that for a nonzero almost periodic solution $u:I\times\R^4\to\C$, there exists $C(u)>0$ such that
		$$\inf_{t\in I}\int_{\vert x-x(t)\vert\leq\frac{C(u)}{N(t)}}\vert u(t,x)\vert^{2p}\,dx\gtrsim_u 1.$$
\end{remark}

The modulation parameters of almost periodic solutions can be shown to obey the following local constancy property (see \cite[Lemma~5.18]{KVnote} for details).

\begin{lemma}[Local constancy]\label{local constancy} Let $u:I\times\R^4\to\C$ be a maximal-lifespan almost periodic solution to \eqref{equ1.1}. Then there exists $\delta=\delta(u)>0$ such that for all $t_0\in I$,
	$$[t_0-\delta N(t_0)^{-2},t_0+\delta N(t_0)^{-2}]\subset I.$$
Moreover,  $N(t)\sim_u N(t_0)$
for $\vert t-t_0\vert\leq\delta N(t_0)^{-2}.$ 
\end{lemma}

Using this local constancy property, we can divide the lifespan $I$ of an almost periodic solution $u$ into \emph{characteristic subintervals} $J_k$ on which we can take $N(t)$ to be constant and equal to some $N_k$, with $\vert J_k\vert\sim_u N_k^{-2}$. To do this, we need to modify the compactness modulus function by some time-independent multiplicative factor. 

Using the local constancy property, we can also deduce information about the behavior of the frequency scale function at the blowup time (see \cite[Corollary~5.19]{KVnote} for details).

\begin{corollary}[$N(t)$ at blowup]\label{constancy cor}
Let $u:I\times\R^4\to\C$ be a maximal-lifespan almost periodic solution to \eqref{equ1.1}. If $T$ is a finite endpoint of $I$, then $N(t)\gtrsim_u\vert T-t\vert^{-1/2}$. In particular, $\lim_{t\to T}N(t)=\infty.$ If $I$ is infinite or semi-infinite, then for any $t_0\in I$ we have $N(t)\gtrsim_u \langle t-t_0\rangle^{-1/2}.$ 
\end{corollary}

We also have the following result relating the frequency scale function of an almost periodic solution to its Strichartz norms. 
\begin{lemma}[Spacetime bounds]\label{fspssn} Let $u:I\times\R^4\to\C$ be an almost periodic solution to \eqref{equ1.1}. Then
$$\int_I N(t)^2\,dt\lesssim_u\bxnorm{\nsc u}{2}{4}{I}^2\lesssim_u 1+\int_I N(t)^2\,dt.$$ 
\end{lemma}

To prove Lemma~\ref{fspssn}, one can adapt the proof of \cite[Lemma~5.21]{KVnote}. The key is to notice that $\int_I N(t)^2\,dt$ approximately counts the number of characteristic subintervals inside $I$ and that $\xonorm{\nsc u}{2}{4}\sim_u 1$ on each such subinterval.

We are now in a position to state precisely the first main step in the proof of Theorem~\ref{theorem}.
\begin{theorem}[Reduction to almost periodic solutions]\label{reduction to ap} If Theorem~\ref{theorem} fails, then there exists a maximal-lifespan solution $u:I\times\R^4\to\C$ to \eqref{equ1.1} that is almost periodic, blows up forward and backward in time, and is minimal in the following sense:
		$$\norm{u}_{L_t^\infty\dot{H}_x^{s_c}(I\times\R^4)}\leq\norm{v}_{L_t^\infty\dot{H}_x^{s_c}(J\times\R^4)}$$
for all maximal-lifespan solutions $v:J\times\R^4\to\C$ that blow up in at least one time direction.
\end{theorem}

The reduction to almost periodic solutions is now widely regarded as a standard technique in the study of dispersive equations at critical regularity. Keraani \cite{Ker} was the first to prove the existence of minimal blowup solutions, while Kenig--Merle \cite{KM} were the first to use them to establish a global well-posedness result. Since then, the technique has proven to be extremely useful; see \cite{KM, KM2010, KMNLW, KTV2009, KV2010, KV20101, KV, KVrad, KVnote, KV3D, KVZ2008, Mu, Mu2, TVZ2007} for many more examples of this technique in action (and note that this is by no means an exhaustive list). For a good introduction to these methods, see \cite{KVnote}.

The proof of Theorem~\ref{reduction to ap} requires three main ingredients. First, one needs a linear profile decomposition for the linear propagator. The first such results were adapted to the mass- and energy-critical settings (see \cite{BeVa, CarKer, Ker0, MerVeg}), while the remaining cases were treated in \cite{Shao}. The second ingredient is a stability result for the nonlinear equation. We have already discussed such a result above (see Theorem~\ref{pertu123}). The third ingredient is a decoupling statement for nonlinear profiles. 

The last two ingredients are closely related, in the sense that the decoupling must hold in a space that is dictated by the stability theory. Most often, this means that the decoupling must hold in a space with $s_c$ derivatives. Keraani \cite{Ker0} showed how to prove such a decoupling statement in the context of the energy-critical NLS. The argument relies on pointwise estimates to bound the difference of nonlinearities and hence is also applicable in the mass-critical setting. For the case of non-conserved critical regularity, however, the nonlocal operator $\nsc$ enters the picture and prevents the direct use of this argument. In \cite{KV2010}, Killip--Vi\c{s}an present an argument for how to overcome this obstacle in the energy-supercritical setting. In particular, by making use of a square function of Strichartz that shares estimates with $\nsc$, one can arrive back in a position where the arguments of Keraani may be applied. (We note here that the argument in \cite{KV2010} does not work in the inter-critical setting. For some alternate approaches adapted to the inter-critical regime, see \cite{HR, KM2010, Mu, Mu2}.)

We therefore have all of the ingredients necessary for Theorem~\ref{reduction to ap}. Putting them together via the usual arguments, we can deduce that Theorem~\ref{reduction to ap} holds. 

With Theorem~\ref{reduction to ap} in place, we can now make some refinements to the class of solutions that we consider. First, a rescaling argument as in \cite{KTV2009, KV20101, TVZ2007} shows that we can restrict our attention to almost periodic solutions that do not escape to arbitrarily low frequencies on at least half of their maximal lifespan (say $[0,T_{max}))$. Second, following the lead of \cite{Dodson3}, we divide solutions into two classes based off of the interaction Morawetz inequality; these will correspond to `rapid frequency-cascade' solutions and `quasi-soliton' solutions. Third, as described above, we use Lemma~\ref{local constancy} to subdivide $[0,T_{max})$ into characteristic subintervals $J_k$ and set $N(t)\equiv N_k$ on each $J_k$, with $\vert J_k\vert \sim_u N_k^{-2}.$ Altogether, we have the following

\begin{theorem}[Two scenarios for blowup]\label{threenew} If Theorem~\ref{theorem} fails, then there exists an almost periodic solution $u:[0,T_{max})\times\R^4\to\C$ that blows up forward in time and satisfies
				$$N(t)\equiv N_k\geq 1$$
for $t\in J_k$, where $[0,T_{max})=\cup_k J_k$ and $\vert J_k\vert\sim_u N_k^{-2}$. Furthermore, 
$$\text{either}\quad
		\int_0^{T_{max}}N(t)^{3-4s_c}\,dt<\infty\quad
\text{or}\quad
		\int_0^{T_{max}}N(t)^{3-4s_c}\,dt=\infty.$$
	\end{theorem}

	Hence to prove Theorem~\ref{theorem} it suffices to preclude the existence of the solutions described in Theorem~\ref{threenew}. 

In Section \ref{frequency-cascade} we preclude the existence of \emph{rapid frequency-cascade solutions}, that is, solutions as in Theorem~\ref{threenew} for which $\int_0^{T_{max}}N(t)^{3-4s_c}\,dt<\infty.$ This requires two main ingredients. The first ingredient is a long-time Strichartz estimate, Theorem~\ref{lse}. As mentioned above, such estimates were first established by Dodson \cite{Dodson3} in the mass-critical setting, but have since appeared in the energy-critical and inter-critical settings \cite{KV3D, Mu, Visan2011}. We will establish a long-time Strichartz estimate for the first time in the energy-supercritical regime. It is in the proof of the long-time Strichartz estimate that we will first encounter the restriction $1<s_c<3/2$ (see Remark~\ref{restriction 1}). The second ingredient is the following reduced Duhamel formula, which is a consequence of almost periodicity (see \cite[Proposition~5.23]{KVnote}):

\begin{lemma} [Reduced Duhamel formula]\label{nowaste} Let $u: [0,T_{max})\times\R^4\to\C$ be a maximal-lifespan almost periodic solution to \eqref{equ1.1}. Then for all $t\in [0,T_{max})$ we have
\begin{equation}\label{nwd}
\begin{split}
u(t)=&\lim_{T\nearrow T_{max}}i\int_t^{T}e^{i(t-s)\Delta}F(u(s))\,ds\end{split}
\end{equation}
as a weak limit in $\dot{H}_x^{s_c}(\R^4)$.
\end{lemma}

Using the long-time Strichartz estimate and the reduced Duhamel formula, we can show that a rapid frequency-cascade solution must have finite mass. In fact, we can show that the solution must have \emph{zero} mass, contradicting the fact that the solution blows up.

In Section \ref{qs section}, we preclude the existence of \emph{quasi-soliton solutions}, that is, solutions as in Theorem~\ref{threenew} for which $\int_0^{T_{max}}N(t)^{3-4s_c}\,dt=\infty.$ The key ingredient in this case is a frequency-localized interaction Morawetz inequality (Theorem~\ref{uppbd}), which we will prove in Section \ref{flim section}. The idea of the proof of Theorem~\ref{uppbd} is to truncate to high frequencies in the standard interaction Morawetz inequality and to estimate the resulting error terms in an acceptable fashion. The main tool we will use to estimate these errors is the long-time Strichartz estimate (Theorem~\ref{lse}), although we will also need to make use of an additional decay result (Proposition~\ref{adecay1}) that we establish in Section \ref{decay section}, as well as a Strichartz-type estimate (Proposition~\ref{mxs}) that we discuss in Section \ref{sze}. 

One interesting feature of the interaction Morawetz inequality that we establish is the use of space localization in the proof. The only other setting in which we have seen both frequency and space localization used to prove an interaction Morawetz inequality is the energy-critical setting in three space dimensions \cite{CKSTT07, KV3D}. In fact, these sources provided a great deal of inspiration for the methods we employ. In particular, we will apply a standard interaction Morawetz identity (Proposition~\ref{morawetz identity}) with a weight that is chosen to equal the usual weight $a(x)=\vert x\vert$ in a large ball, but which is eventually constant. This additional spatial truncation is necessary to deal with certain error terms that result from the frequency localization. Of course, localizing in space introduces even more error terms to deal with, but in the end we can treat them all and arrive at a useful estimate. In fact, we will see that some terms that are quite troublesome in the energy-critical setting become relatively easy to deal with in the energy-supercritical regime. At various points in the argument, we will once again encounter the restriction $1<s_c<3/2$. For further discussion, see Section~\ref{flim section}.

To preclude the existence of quasi-soliton solutions and thereby complete the proof of Theorem~\ref{theorem}, we notice that we can use the frequency-localized interaction Morawetz inequality to get uniform control over $\int_I N(t)^{3-4s_c}\,dt$ for compact time intervals $I\subset[0,T_{max})$. We can thus derive a contradiction by taking $I$ sufficiently large inside $[0,T_{max})$. 

\subsection*{Acknowledgements} 
C. M. was supported by the NSF of China under grant No. 11171033 and 11231006. 
J. M. was supported in part by NSF grant DMS-1001531 (P.I. Rowan Killip). 
J. Z. was supported by the NSF of China under grant No. 11171033.

\section{Notation and useful lemmas}
\subsection{Some notation}
For nonnegative quantities $X$ and $Y$, we will write $X\lesssim Y$ to denote the estimate $X\leq C Y$ for some $C>0$. If $X\lesssim Y\lesssim X$, we will write $X\sim Y$. Dependence of implicit constants on the power $p$ or the dimension will be suppressed; dependence on additional parameters will be indicated by subscripts. For example, $X\lesssim_u Y$ indicates $X\leq CY$ for some $C=C(u)$. 

We will use the expression $\text{\O}(X)$ to denote a finite linear combination of terms that resemble $X$ up to Littlewood--Paley projections, complex conjugation, and/or maximal functions. We will use the expression $X\pm$ to denote $X\pm\eps$ for any $\eps>0$. 

For a spacetime slab $I\times\R^4$, we write $L_t^q L_x^r(I\times\R^4)$ for the Banach space of functions $u:I\times\R^4\to\C$ equipped with the norm
	$$\xnorm{u}{q}{r}{I}:=\bigg(\int_I \norm{u(t)}_{L_x^r(\R^4)}\,dt\bigg)^{1/q},$$
with the usual adjustments when $q$ or $r$ is infinity. When $q=r$, we abbreviate $L_t^qL_x^q=L_{t,x}^q$. We will also abbreviate $\norm{f}_{L_x^r(\R^4)}$ to $\norm{f}_{L_x^r}.$ For $1\leq r\leq\infty$, we use $r'$ to denote the dual exponent to $r$, i.e. the solution to $\tfrac{1}{r}+\tfrac{1}{r'}=1.$

We define the Fourier transform on $\mathbb{R}^4$ by
\begin{equation*}
\aligned \widehat{f}(\xi):= \tfrac{1}{4\pi^2}\int_{\mathbb{R}^4}e^{- ix\cdot \xi}f(x)\,dx.
\endaligned
\end{equation*}
We can then define the fractional differentiation operator $\vert\nabla\vert^s$ for $s\in \R$ via 
	$$\wh{\vert\nabla\vert^s f}(\xi):=\vert\xi\vert^s\wh{f}(\xi),$$
with the corresponding homogeneous Sobolev norm 
	$$\norm{f}_{\dot{H}_x^s(\R^4)}:=\norm{\vert\nabla\vert^s f}_{L_x^2(\R^4)}.$$

\subsection{Basic harmonic analysis}

We will make frequent use of the Littlewood--Paley projection operators. 
Specifically, we let $\varphi$ be a radial bump function supported
on the ball $|\xi|\leq 2$ and equal to 1 on the ball $|\xi|\leq 1$.
For $N\in 2^{\mathbb{Z}}$, we define the Littlewood--Paley projection operators by
\begin{align*}
&\widehat{P_{\leq N}f}(\xi) := \wh{f_{\leq N}}(\xi):= \varphi(\xi/N)\widehat{f}(\xi),
 \\ &\widehat{P_{> N}f}(\xi) := \wh{f_{>N}}(\xi) :=
(1-\varphi(\xi/N))\widehat{f}(\xi), 
\\ &\widehat{P_{N}f}(\xi) := \wh{f_N}(\xi):=
(\varphi(\xi/N)-\varphi({2\xi}/{N}))\widehat{f}(\xi).
\end{align*}
We may also define $$P_{M<\cdot\leq N}:=P_{\leq N}-P_{\leq M}=\sum_{M<N'\leq N}P_N'$$
for $M<N.$ All such summations should be understood to be over $N'\in 2^{\mathbb{Z}}.$ 

The Littlewood--Paley operators commute with derivative operators,
the free propagator, and the conjugation operation. These operators are
self-adjoint and bounded on every $L^p_x$ and $\dot{H}^s_x$ space
for $1\leq p\leq \infty$ and $s\geq 0$. They also obey the
following standard Bernstein estimates:
\begin{lemma}[Bernstein estimates] For $1\leq r\leq q\leq\infty$ and $s\geq 0$, 
\begin{align*}
\big\||\nabla|^{\pm s} P_{N} f \big\|_{L_x^r(\R^4)} & \thicksim  N^{\pm s}
\big\|
P_{N} f \big\|_{L_x^r(\R^4)},  \\
\big\||\nabla|^s P_{\leq N} f \big\|_{L_x^r(\R^4)} & \lesssim   N^{s}
\big\|
P_{\leq N} f \big\|_{L_x^r(\R^4)},  \\
 \big\| P_{> N} f \big\|_{L_x^r(\R^4)} & \lesssim  N^{-s} \big\|
|\nabla|^{s}P_{> N} f \big\|_{L_x^r(\R^4)}, \\
\big\| P_{\leq N} f \big\|_{L^q(\R^4)} & \lesssim
N^{\frac{4}{r}-\frac{4}{q}} \big\|
P_{\leq N} f \big\|_{L_x^r(\R^4)}.
\end{align*}

\end{lemma}

We will also need the following fractional calculus estimates from \cite{CW}. 

\begin{lemma}[Fractional product rule \cite{CW}]
Let $s\geq0$ and let $1<r,r_j,q_j<\infty$ satisfy
$\frac1r=\frac1{r_j}+\frac1{q_j}$ for $j=1,2$. Then
\begin{equation}\label{moser}
\big\||\nabla|^s(fg)\big\|_{L_x^r}\lesssim\|f\|_{{L_x^{r_1}}}\big\||\nabla|^sg
\big\|_{{L_x^{q_1}}}+\big\||\nabla|^sf\big\|_{{L_x^{r_2}}}\|g\|_{{L_x^{q_2}}}.
\end{equation}
\end{lemma}

\begin{lemma}[Fractional chain rule \cite{CW}]
Let $G\in
C^1(\mathbb{C}),~s\in(0,1],$ and $1<r,r_1,r_2<\infty$ satisfy $\frac1r=\frac1{r_1}+\frac1{r_2}.$
Then
\begin{equation}\label{fraclsfz}
\big\||\nabla|^s
G(u)\big\|_{L_x^r}\lesssim\|G'(u)\|_{L_x^{r_1}}\big\||\nabla|^su\big\|_{L_x^{r_2}}.
\end{equation}
\end{lemma}

Finally, we will also need a paraproduct estimate in the spirit \cite[Lemma~2.3]{Visan2011}. This estimate will serve as an important technical tool when we prove the long-time Strichartz estimate in Section \ref{lts section}, as well as when we prove the frequency-localized interaction Morawetz inequality in Section \ref{flim section}.
 
\begin{lemma}[Paraproduct estimate]\label{pde} Let $0<s<1$. If $1<r<r_1<\infty$ and $1<r_2<\infty$ satisfy
		$\tfrac{1}{r_1}+\tfrac{1}{r_2}=\tfrac{1}{r}+\tfrac{s}{4}<1,$
then 
		\begin{equation}\label{pde1}
		\norm{\vert\nabla\vert^{-s}(fg)}_{L_x^r(\R^4)}\lesssim
		\norm{\vert\nabla\vert^{-s}f}_{L_x^{r_1}(\R^4)}\norm{\vert\nabla\vert^{s}g}_{L_x^{r_2}(\R^4)}.
		\end{equation}
\end{lemma}
\begin{proof}
We will prove the equivalent estimate
	\begin{equation}\label{pdeeq}
	\norm{\vert\nabla\vert^{-s}(\vert\nabla\vert^sf\,\vert\nabla\vert^{-s}g)}_{L_x^r(\R^4)}
	\lesssim\norm{f}_{L_x^{r_1}(\R^4)}\norm{g}_{L_x^{r_2}(\R^4)}.
	\end{equation}

We decompose the left-hand side into low-high and high-low frequency interactions. In particular, we define the projections $\pi_{l,h}$ and $\pi_{h,l}$ by
	$$\pi_{l,h}(\phi,\psi):=\sum_{N\lesssim M}\phi_N\psi_M,\quad\pi_{h,l}(\phi,\psi):=\sum_{N\gg M}\phi_N\psi_M$$
for any pair of functions $\phi$ and $\psi$.

We first consider the low-high interactions. To begin, we apply Sobolev embedding to get
	\begin{equation}
	\norm{\vert\nabla\vert^{-s}\pi_{l,h}(\vert\nabla\vert^s f,\vert\nabla\vert^{-s}g)}_{L_x^{r}(\R^4)}
	\lesssim\norm{\pi_{l,h}(\vert\nabla\vert^{s}f,\vert\nabla\vert^{-s}g)}_{L_x^{\frac{4r}{4+rs}}(\R^4)}.
	\label{para sobo}
	\end{equation}
Note that the assumption $\tfrac{1}{r}+\tfrac{s}{4}<1$ guarantees $\tfrac{4r}{4+rs}>1$. We now notice that the multiplier of the operator
		$$T(f,g):=\pi_{l,h}(\vert\nabla\vert^{s}f,\vert\nabla\vert^{-s}g),$$
that is
		$$\sum_{N\lesssim M}\vert\xi_1\vert^s\wh{f_N}(\xi_1)\vert\xi_2\vert^{-s}\wh{g_M}(\xi_2),$$
is a symbol of order zero with $\xi=(\xi_1,\xi_2)$. Thus a theorem of Coifman--Meyer \cite{CM1975,CM1991} allows us to continue from \eqref{para sobo} to obtain 
	$$\norm{\vert\nabla\vert^{-s}\pi_{l,h}(\vert\nabla\vert^s f,\vert\nabla\vert^{-s}g)}_{L_x^{r}(\R^4)}
	\lesssim \norm{f}_{L_x^{r_1}(\R^4)}\norm{g}_{L_x^{r_2}(\R^4)}.
	$$

	We next consider the high-low interactions. Noting that the multiplier of the operator
		$$S(f,h):=\vert\nabla\vert^{-s}\pi_{h,l}(\vert\nabla\vert^s f,h),$$
that is,
		$$\sum_{N\gg M}\vert\xi_1+\xi_2\vert^{-s}\vert\xi_1\vert^s\wh{f_N}(\xi_1)\wh{h_M}(\xi_2),$$
is a symbol of order zero as well, we can use the result cited above (together with Sobolev embedding) to deduce
	\begin{align*}
	\norm{\vert\nabla\vert^{-s}\pi_{h,l}(\vert\nabla\vert^s f,\vert\nabla\vert^{-s} g)}_{L_x^r(\R^4)}
	&\lesssim\norm{f}_{L_x^{r_1}(\R^4)}\norm{\vert\nabla\vert^{-s}g}_{L_x^{\frac{rr_1}{r_1-r}}(\R^4)}
	\\ &\lesssim	\norm{f}_{L_x^{r_1}(\R^4)}\norm{g}_{L_x^{r_2}(\R^4)}.
	\end{align*}

Combining the low-high and high-low interactions yields \eqref{pdeeq}.
\end{proof}

 \subsection{Strichartz estimates}\label{sze}
 Let $e^{it\Delta}$ be the free Schr\"odinger propagator, given by
 	\begin{equation}\label{explicit formula}
	[e^{it\Delta}f](x)=\tfrac{-1}{16\pi^2 t^2}\int_{\R^4} e^{i\vert x-y\vert^2/4t}f(y)\,dy
	\end{equation}
for $t\neq 0$. From this explicit formula we can read off the dispersive estimate
	$$\norm{e^{it\Delta}f}_{L_x^\infty(\R^4)}\lesssim\vert t\vert^{-2}\norm{f}_{L_x^1(\R^4)}$$ 
for $t\neq 0$. Interpolating with $\norm{e^{it\Delta}f}_{L_x^2(\R^4)}\equiv \norm{f}_{L_x^2(\R^4)}$ then yields
\begin{equation}\label{dispers}
\big\|e^{it\Delta}f \big\|_{L^r_x(\R^4)} \leq C|t|^{-2(1-\frac2{r})}
\|f\|_{L^{r'}_x(\R^4)}
\end{equation}
for $t\neq 0$ and $2\leq r\leq\infty$, where $\tfrac{1}{r}+\tfrac{1}{r'}=1.$ This estimate implies the standard Strichartz estimates, which we will state below. First, we need the following

\begin{definition}[Admissible pairs]\label{def1} A pair of exponents $(q,r)$ is called \emph{Schr\"odinger admissible} if $2\leq q,r\leq\infty$ and $\tfrac{2}{q}+\tfrac{4}{r}=2.$ For a spacetime slab $I\times\R^4$, we define the Strichartz norm 
	$$\norm{u}_{S^0(I)}:=\sup\big\{\xnorm{u}{q}{r}{I}:(q,r)\text{ Schr\"odinger admissible}\big\}.$$
We denote $S^0(I)$ to be the closure of all
test functions under this norm and write $N^0(I)$ for the dual of
$S^0(I)$.
\end{definition}

We may now state the standard Strichartz estimates in the form that we will need them.

\begin{proposition}[Strichartz \cite{GV, KeT98, St}]\label{prop1}
Let $s\geq 0$ and suppose $u:I\times\R^4\to\C$ is a solution to $(i\partial_t+\Delta)u=F$. Then
	$$\norm{\vert\nabla\vert^s u}_{S^0(I)}\lesssim\norm{\vert\nabla\vert^s u(t_0)}_{L_x^2(\R^4)}+\norm{\vert\nabla\vert^s F}_{N^0(I)}$$
for any $t_0\in I$. 
\end{proposition}
As mentioned above, the key ingredient for the proof of Proposition~\ref{prop1} is the estimate \eqref{dispers}. For the endpoint $(q,r)=(2,4)$, see \cite{KeT98}.

Finally, we will need the following Strichartz-type estimate, which we will use in Section \ref{flim section} to control the mass of solutions over balls. One should compare this result to the $q=\infty$ case of \cite[Proposition~3.2]{KV3D}. In fact, all of the necessary ideas to prove the following proposition may be found in \cite[Section 3]{KV3D}; we only need to make small adjustments to adapt the argument to four dimensions.

\begin{proposition}\label{mxs} Let $v$ be a solution to $(i\partial_t+\Delta)v=F+G$ on a spacetime slab $[0,T]\times\R^4$. Let
				\begin{equation}\label{smudge}
				[S_\rho v](t,x):=\rho^{-2}\bigg(\int_{\R^4}\vert v(t,y)\vert^2e^{-\vert x-y\vert^2/\rho^2}\,dy\bigg)^{1/2}
				\end{equation}
for $\rho>0$. Then we have
				\begin{equation}\label{s strichartz}
				\rho\xonorm{S_{\rho}v}{2}{\infty}\lesssim
				\xonorm{v}{\infty}{2}+\xonorm{F}{2}{4/3}+\rho^{-1}\xonorm{G}{2}{1},
				\end{equation}
where all spacetime norms are taken over $[0,T]\times\R^4$. 
\end{proposition}
\begin{proof}
We follow exactly the arguments in \cite[Section~3]{KV3D}, making only slight adjustments in order to adapt the argument to four dimensions. 

We will begin by using a `double Duhamel' trick, which has its origins in \cite[Section~14]{CKSTT07}. In particular, we will use the Duhamel formula to write $v$ in two different ways, namely
				$$v(t)=a(t)+b(t)=c(t)+d(t),$$
where $$a(t)+b(t):=\bigg(e^{it\Delta}v(0)-i\int_0^t e^{i(t-s)\Delta}F(s)\,ds\bigg)-i\int_0^t e^{i(t-s)\Delta}G(s)\,ds$$
and	$$c(t)+d(t):=\bigg(e^{-i(T-t)\Delta}v(T)+i\int_t^Te^{-i(\tau-t)\Delta}F(\tau)\,d\tau\bigg)+i\int_t^T e^{-i(\tau-t)\Delta}G(\tau)\,d\tau.$$

We next note the following basic pointwise estimate, which is a consequence of Cauchy--Schwarz:
		$$\vert v\vert^2\lesssim \vert a\vert^2+\vert c\vert^2+\vert{b}\bar{d}\vert.$$
Hence
	\begin{align}
	\big\vert S_\rho v(t,x)\vert^2\lesssim&\ \rho^{-4}\int e^{-\vert x-z\vert^2/\rho^2}\{\vert a(t,z)\vert^2+\vert c(t,z)\vert^2\}\,dz
					\label{maxm1}
	\\ &+\rho^{-4}\int e^{-\vert x-z\vert^2/\rho^2}\vert{b}(t,z)\bar{d}(t,z)\vert\,dz.
					\label{maxm2}
	\end{align} 

We use Young's inequality and Strichartz to estimate the contribution of \eqref{maxm1} to \eqref{s strichartz}. We find
	\begin{align*}
	\rho^2&\int_0^T \sup_x\rho^{-4} \int e^{-\vert x-z\vert^2/\rho^2}\{\vert a(t,z)\vert^2+\vert c(t,z)\vert^2\}\,dz\,dt
	\\ &\lesssim\rho^{-2}\norm{e^{-\vert\cdot\vert^2/\rho^2}}_{L_x^2}\int_0^T \big\{\norm{a(t)}_{L_x^4}^2 +\norm{c(t)}_{L_x^4}^2\big\}\,dt
	\\ &\lesssim \xonorm{a}{2}{4}^2+\xonorm{c}{2}{4}^2
	\\ &\lesssim \big(\xonorm{v}{\infty}{2}+\xonorm{F}{2}{4/3}\big)^2.
	\end{align*}
Taking square roots, we see that the contribution of \eqref{maxm1} is acceptable.

We now turn to \eqref{maxm2}. Recalling the definitions of $b$ and $d$, performing a change of variables, and moving the propagator, we see that to estimate the contribution of \eqref{maxm2} to \eqref{s strichartz}, we need to show
	\begin{align}
	\rho^2\int_0^T\!&\sup_x\rho^{-4}\bigg\vert\int_0^{T-t}\int_0^t\int \overline{G}(t+\tau,z)e^{i\tau\Delta}e^{-\vert x-z\vert^2/\rho^2}e^{is\Delta}G(t-s,z)\,dz\,d\tau\,ds\bigg\vert\,dt\nonumber
	\\ &\lesssim\rho^{-2}\xonorm{G}{2}{1}^2.\label{kernel}
	\end{align}

Using the explicit formula for the free propagator \eqref{explicit formula}, we rewrite
	$$\rho^{-4}e^{i\tau\Delta}e^{-\vert x-z\vert^2/\rho^2}e^{is\Delta}G(t-s,z)=\int K_{\rho}(\tau,z;s,y;x)G(t-s,y)\,dy,$$
where 
	$$K_{\rho}(\tau,z;s,y;x)= (16\pi^2\rho^2\tau s)^{-2}\int\exp\{i\vert z-w\vert^2/4\tau-\vert x-w\vert^2/\rho^2+i\vert w-y\vert^2/4s\}\,dw.$$
Completing the square and evaluating the Gaussian integral yields
	$$\norm{K_\rho(\tau,z;s,y;x)}_{L_{x,y,z}^\infty}\sim [16s^2\tau^2+\rho^4(s+\tau)^2]^{-1}.$$
Using this estimate and changing variables via $\alpha=\rho^{-2}\tau$ and $\beta=\rho^{-2}s$, we find
	\begin{align*}
	&\LHS\eqref{kernel}
	\\ &\lesssim \rho^2\int_0^T\int_0^{T-t}\int_0^t[16s^2\tau^2+\rho^4(s+\tau)^2]^{-1}
		\norm{G(t-s)}_{L_x^1}\norm{G(t+\tau)}_{L_x^1}\,d\tau\,ds\,dt
	\\ &\lesssim \rho^{-2}\int_0^T\!\int_0^\infty\!\int_0^\infty\![\alpha^2\beta^2+(\alpha+\beta)^2]^{-1}
		\norm{G(t-\rho^2\beta)}_{L_x^1}\norm{G(t+\rho^2\alpha)}_{L_x^1}\,d\alpha\, d\beta\,dt.
	\end{align*}
We next claim that
	\begin{align}
	\int_0^\infty\int_0^\infty&[\alpha^2\beta^2+(\alpha+\beta)^2]^{-1}\norm{G(t-\rho^2\beta)}_{L_x^1}\norm{G(t+\rho^2\alpha)}_{L_x^1}\,d\alpha\, d\beta\nonumber
	\\& \lesssim [\mathcal{M}\big(\norm{G(\cdot)}_{L_x^1}\big)(t)]^2,\label{maximal fn}
	\end{align}
where $\mathcal{M}$ denotes the Hardy--Littlewood maximal function. Once \eqref{maximal fn} is established, \eqref{kernel} follows easily from the maximal function estimate. Thus, to complete the proof of Proposition~\ref{mxs}, it remains to prove \eqref{maximal fn}. 

The argument for \eqref{maximal fn} is now identical to the argument appearing in the proof of Lemma~3.4 in \cite{KV3D}. In particular, one uses the fundamental theorem of calculus to show that the function 
		$$[\alpha^2\beta^2+(\alpha+\beta)^2]^{-1}$$
can be bounded above by a convex combination of $L^1$-normalized characteristic functions of rectangles $[0,\ell]\times[0,w]$. For the details, see discussion following display (3.7) in \cite{KV3D} (with $q=\infty$).
	\end{proof}

\subsection{A Gronwall inequality}

We record here a technical result from \cite{KV20101}, which will be of use in Section \ref{decay section}.

\begin{lemma}[Acausal Gronwall inequality \cite{KV20101}]\label{gronwall1} 
Let $\gamma>0$, $0<\eta<\tfrac12(1-2^{-\gamma})$ and $\{b_k\}\in\ell^\infty(\mathbb{Z}^+)$. Suppose $\{x_k\}\in\ell^\infty(\mathbb{Z}^+)$ is a nonnegative sequence that satisfies
	$$x_k\leq b_k+\eta\sum_{\ell=0}^\infty 2^{-\gamma\vert k-\ell\vert}x_{\ell}$$ for all $k\geq 0.$ 
Then there exists $r=r(\eta)\in(2^{-\gamma},1)$ such that
	$$x_k\lesssim\sum_{\ell=0}^kr^{\vert k-\ell\vert}b_{\ell}$$
for all $k\geq 0$, with $r\to 2^{-\gamma}$ as $\eta\to 0$. 
\end{lemma}

%%%%%%%%%%%%%%%%%%%%%%%%%%%%%%%%%%%%%%%%%%%%%%%%%%%%%%%%%%%%%%%%%%%
\section{Additional decay}\label{decay section}
In this section, we will show that almost periodic solutions to \eqref{equ1.1} as in Theorem~\ref{threenew} enjoy `additional decay', in the sense that they belong to $L_t^\infty L_x^q$ for values of $q$ that are smaller than $2p$ (recall that $2p$ is the exponent obtained by applying Sobolev embedding to $\dot{H}_x^{s_c}$). In particular, we will show that the solutions we consider belong to $L_t^\infty L_x^{p+1}.$

The arguments we present in this section have their origin in \cite[Section 6]{KV20101}, wherein the authors show that almost periodic solutions to the focusing energy-critical NLS in dimensions $d\geq 5$ have finite mass (in fact, they belong to $L_t^\infty \dot{H}_x^{-\eps}$ for some $\eps>0$). The proof in \cite{KV20101} consists of two steps: first, establish additional decay in the $L_x^q$ sense mentioned above (also referred to as `breaking scaling'), and second, employ a `double Duhamel' trick to upgrade this information to negative regularity. A similar argument was also used in \cite{KV2010} to establish negative regularity in the defocusing energy-supercritical regime in dimensions $d\geq 5$. 

As we will see, the first step of the proof appearing in \cite{KV2010, KV20101} carries over directly to our setting. It consists of a bootstrap argument based off of the reduced Duhamel formula, wherein a certain `acausal' Gronwall inequality (Lemma~\ref{gronwall1}) is used to deal with the various frequency interactions. The double Duhamel trick, on the other hand, fails in dimensions $d<5$ due to the fact that the dispersive effect of $e^{it\Delta}$ is too weak in low dimensions (recall that the convolution kernel of $e^{it\Delta}$ decays like $\vert t\vert^{-d/2}$). Thus, while we can establish additional decay, we cannot use it to prove negative regularity. 

Fortunately, proving additional decay in the $L^q$-sense will suffice for our purposes. In particular, we will make use of the fact that $u\in L_t^\infty L_x^{p+1}$ in Section \ref{flim section} in the proof of the frequency-localized interaction Morawetz inequality. 

\begin{proposition}[Additional decay]\label{adecay1} Let $1<s_c<3/2$. Suppose $u:[0,T_{max})\times\R^4\to\C$ is an almost periodic solution to \eqref{equ1.1} such that $\inf_{t\in[0,T_{max})}N(t)\geq 1$. Then
		\begin{equation}u\in L_t^\infty L_x^{p+1}([0,T_{max})\times\R^4).
		\label{breaking scaling}
		\end{equation}
\end{proposition}

\begin{proof} In fact, we will establish 
					\begin{equation}\label{ad1}
				u\in L_t^\infty L_x^q([0,T_{max})\times\R^4)\quad\text{for}\quad \tfrac{4p(p+1)}{p^2+p+2}<q< 2p,
					\end{equation}
which implies \eqref{breaking scaling}. We first note that to prove \eqref{ad1}, it will suffice to show that for $\tfrac{4p(p+1)}{p^2+p+2}<q<{2p}$, there exist $N_1>0$ and $\alpha>0$ such that 
				\begin{equation}\label{quant}
				\|u_N\|_{L_t^\infty L_x^q}\lesssim_u N^{\al}\quad \text{for}\quad N\leq N_1.
				\end{equation}
Indeed, combining \eqref{quant} with Bernstein yields
\begin{align*}
\|u\|_{L_t^\infty L_x^q}&\lesssim\sum_{N\leq N_1}\|u_N\|_{L_t^\infty L_x^q}+\|u_{>N_1}\|_{L_t^\infty L_x^q}\\
&\lesssim_u\sum_{N\leq N_1}N^{\al}+\||\nabla|^{2-\frac{4}{q}}u_{>N_1}\|_{L_t^\infty L_x^2}\\
&\lesssim_u N_1^{\alpha}+N_1^{2-\frac{4}{q}-s_c}\big\||\nabla|^{s_c}u\big\|_{L_t^\infty L_x^2}\\
&\lesssim_u 1.
\end{align*}

We therefore turn to establishing \eqref{quant}. Let $\eta>0$ be a small parameter to be determined later. Using the fact that $\inf_{t\in[0,T_{max})} N(t)\geq 1$, we may find $N_0=N_0(\eta)>0$ such that
		\begin{equation}\label{jxsmall}
		\xonorm{\nsc u_{\leq N}}{\infty}{2}<\eta
		\end{equation}
for $N\leq N_0.$ 

For $q>4$, we now define the quantity
		\begin{equation}
		%\label{aqn1}
		\nonumber
		A_q(N):=N^{\frac{4}{q}-\frac{2}{p}}\xnorm{u_N}{\infty}{q}{[0,T_{max})}.
		\end{equation}
Using Bernstein and the fact that $u\in L_t^\infty\dot{H}_x^{s_c}$, we first note that $A_q(N)$ satisfies 
		\begin{equation}
		\label{A_q bounded}
		A_q(N)\lesssim_u 1.
		\end{equation}

We next claim that $A_q(N)$ satisfies the following recurrence relation:

\begin{lemma}[Recurrence relation for $A_q(N)$]\label{adrr} For $4<q<\frac{4p}{p-2}$ and $N\leq 10N_0,$ we have
\begin{align}
 A_q(N)&\nonumber
\lesssim_u(\tfrac{N}{N_0})^{2-\frac2p-\frac{4}q}+\eta^p\sum_{\frac{N}{10}<
M\leq N_0}(\tfrac{N}{M})^{2-\frac2p-\frac{4}q}A_q(M)
\\ &\quad\quad+\eta^p \sum_{M\leq\frac{N}{10}}(\tfrac{M}{N})^{-1+\frac{4}q+\frac2p}A_q(M).  \label{diedai125}
\end{align}
%Note $2-\frac2p-\frac4q>0$ and $-1+\frac4p+\frac2p>0$ for this range of $q.$
\end{lemma}			

We will prove Lemma~\ref{adrr} below. Let us first see how we can use it to complete the proof of Proposition~\ref{adecay1}. As in \cite{KV2010, KV20101}, the key to extracting a bound from this recurrence relation will be to make use of an `acausal' Gronwall inequality. Specifically, we will make use of Lemma~\ref{gronwall1}, which we have imported directly from \cite{KV20101}. 

We apply Lemma~\ref{gronwall1} with $N=2^{-k}\cdot 10N_0$, $x_k=A_q(N)$, and $\eta$ chosen sufficiently small. Note that $\{x_k\}\in \ell^\infty$ by \eqref{A_q bounded}. Using \eqref{diedai125} and Lemma~\ref{gronwall1}, we deduce
	$$A_q(N)\lesssim_u N^{(2-\frac2p-\frac{4}q)-},\quad\text{i.e.}\quad\xonorm{u_N}{\infty}{q}\lesssim_u N^{(2-\frac{8}{q})-} $$
for $4<q<\tfrac{4p}{p-2}$ and $N\leq 10N_0$. In particular, we find
	$$\xonorm{u_N}{\infty}{\frac{4p}{p-2}-}\lesssim_u N^{\frac{4}{p}-}$$
for $N\leq 10N_0$. Hence, by interpolation and Bernstein, we find that for $2<q<\tfrac{4p}{p-2}$ and $N\leq 10N_0$, 
	\begin{align*}
	\xonorm{u_N}{\infty}{q}&\lesssim\xonorm{u_N}{\infty}{\frac{4p}{p-2}-}^{\frac{2p(q-2)}{q(p+2)}-}
							\xonorm{u_{N}}{\infty}{2}^{\frac{4p+2q-pq}{q(p+2)}+}
				\\ &\lesssim_u N^{\frac{8(q-2)}{q(p+2)}-}N^{(-2+\frac{2}{p})(\frac{4p+2q-pq}{q(p+2)}+)}.
	\end{align*}
Thus
	$$\|u_N\|_{L_t^\infty L_x^{\frac{4p(p+1)}{p^2+p+2}+}}\lesssim_u N^{0+}$$
for $N\leq 10N_0.$ Interpolating this estimate with the fact that $u\in L_t^\infty L_x^{2p}$, we deduce that \eqref{quant} holds, which completes the proof of Proposition~\ref{adecay1}. \end{proof}

It remains to establish Lemma~\ref{adrr}. 

\begin{proof}[Proof of Lemma~\ref{adrr}]
Using time-translation symmetry, we see that it suffices to prove
	\begin{align}
	\nonumber N^{\frac{4}{q}-\frac{2}{p}}\norm{u_N(0)}_{L_x^q}&\nonumber
\lesssim_u(\tfrac{N}{N_0})^{2-\frac2p-\frac{4}q}+\eta^p\sum_{\frac{N}{10}<
M\leq N_0}(\tfrac{N}{M})^{2-\frac2p-\frac{4}q}A_q(M)
\\ &\quad\quad+\eta^p \sum_{M\leq\frac{N}{10}}(\tfrac{M}{N})^{-1+\frac{4}q+\frac2p}A_q(M)  \label{diedai1252}
	\end{align}
for $N\leq 10N_0$ and $4<q<\tfrac{4p}{p-2}$.

We will first use the reduced Duhamel formula to write $u_N(0)$ as an integral over $[0,T_{max})$, which we would then like to divide into a short-time piece $(0<t< N^{-2})$ and a long-time piece $(t>N^{-2})$. Of course, if $T_{max}<N^{-2}$, then there will be no long-time piece to estimate. In the following, we will only consider the case $T_{max}> N^{-2}$, so that we have two pieces to estimate. In the end, we will derive the same bounds for both pieces, and hence we lose no generality proceeding in this way.

Using the reduced Duhamel formula \eqref{nwd}, Bernstein, and the dispersive estimate \eqref{dispers}, we estimate
\begin{align*}
\|u_N(0)\|_{L_x^q}&\lesssim\Big\|\int_0^{T_{max}} e^{it\Delta}P_NF(u(t))\,dt\Big\|_{L_x^q}
\\ &\lesssim N^{4(\frac12-\frac1q)}\int_0^{N^{-2}}\| e^{it\Delta}P_NF(u(t))\|_{L_x^2}\,dt
	\\&\quad\quad	+\int_{N^{-2}}^{T_{max}} t^{-4(\frac12-\frac1{q})}\|P_NF(u(t))\|_{L_x^{q'}}\,dt
\\ &\lesssim N^{2-\frac{8}q}\|P_NF(u(t))\|_{L_t^\infty L_x^ {q'}}.
\end{align*}
Hence
\begin{equation}\label{dd11} N^{\frac{4}{q}-\frac{2}{p}}\norm{u_N(0)}_{L_x^q}\lesssim
N^{2-\frac2p-\frac{4}q}\|P_NF(u)\|_{L_t^\infty L_x^ {q'}}
\end{equation} 
for $q>4$. 

For $N\leq 10N_0$, we now use the fundamental theorem of calculus to decompose the nonlinearity $F(u)$ as follows:
\begin{align}
F(u)&
= \text{\O}(u_{>N_0}u^p) \label{ad high}
\\ &\quad+F(u_{\frac{N}{10}<\cdot\leq N_0}) \label{ad mid}
\\ \label{contru1} &\quad+u_{\leq\frac{N}{10}}\int_0^1F_z\big(u_{\frac{N}{10}<\cdot\leq N_0}+\theta u_{\leq\frac{N}{10}}\big)\,d\theta
\\ \label{contru2} &\quad+\overline{u_{\leq\frac{N}{10}}}\int_0^1F_{\bar{z}}\big(u_{\frac{N}{10}<\cdot\leq N_0}+\theta u_{\leq\frac{N}{10}}\big)\,d\theta.
\end{align}

We first consider the contribution of \eqref{ad high} to \eqref{dd11}. Using H\"older, Bernstein, and Sobolev embedding, we find
\begin{align}\nonumber
\big\|P_N\text{\O}(u_{>N_0} u^p)\big\|_{L_t^\infty L_x^ {q'}}
&\lesssim\|u_{> N_0}\|_{L_t^\infty L_x^{\frac{2q}{q-2}}}\|u\|_{L_t^\infty
L_x^{2p}}^p
\\ \nonumber
&\lesssim N_0^{-2+\frac2p+\frac{4}q}\|\nsc u\|_{L_t^\infty
L_x^2}^{p+1}\\\label{di1}
&\lesssim_u N_0^{-2+\frac2p+\frac{4}q}.
\end{align}
Comparing with \eqref{diedai1252}, we see that this term is acceptable.

We next estimate the contribution of \eqref{ad mid} to \eqref{dd11}.
Using H\"older, Bernstein, and \eqref{jxsmall}, we estimate
\begin{align}
	\nonumber
	\big\| F\big(&u_{\frac{N}{10}<\cdot\leq N_0}\big)\big\|_{L_t^\infty L_x^{q'}}
	\\ \nonumber
	&\lesssim\big\|u_{\frac{N}{10}<\cdot\leq N_0}\big\|_{L_t^\infty L_x^{2p}}^{p-1}
	\sum_{\frac{N}{10}< M_1\leq M_2\leq N_0}\|u_{M_1}\|_{L_t^\infty L_x^q}\|u_{M_2}\|_{L_t^\infty L_x^\frac{2pq}{q+qp-4p}}
	\\ \nonumber &\lesssim\eta^{p-1}\sum_{\frac{N}{10}< M_1\leq M_2\leq N_0}\|u_{M_1}\|_{L_t^\infty L_x^q}
	M_2^{-2+\frac{8}{q}}\|\nsc u_{\leq N_0}\|_{L_t^\infty L_x^2}
	\\ \label{di2} 
	&\lesssim\eta^p \sum_{\frac{N}{10}< M\leq N_0}M^{-(2-\frac2p-\frac{4}q)} A_q(M).
\end{align}
Comparing with \eqref{diedai1252}, we see that this term is acceptable as well.

Finally, we estimate the contribution of \eqref{contru1} and \eqref{contru2} to \eqref{dd11}. It will suffice to treat \eqref{contru1}, as the same arguments may be used to handle \eqref{contru2}. Using H\"older, Bernstein, and \eqref{jxsmall}, we estimate
\begin{align}
	\Big\|P_N&\Big(u_{\leq\frac{N}{10}}\int_0^1F_z\big(u_{\frac{N}{10}<\cdot\leq N_0}+\theta u_{\leq\frac{N}{10}}\big)\,d\theta\Big)\Big\|_{L_t^\infty L_x^{q'}}\nonumber
\\
&\lesssim\|u_{\leq\frac{N}{10}}\|_{L_t^\infty L_x^\frac{4q}{q-4}}\Big\|P_{>\frac{N}{10}}\Big(\int_0^1F_z\big(u_{\frac{N}{10}<\cdot\leq N_0}+\theta u_{\leq\frac{N}{10}}\big)\,d\theta\Big)\Big\|_{L_t^\infty L_x^{4/3}}\nonumber
\\
&\lesssim N^{-1}\|u_{\leq\frac{N}{10}}\|_{L_t^\infty L_x^\frac{4q}{q-4}}\big\|\nabla u_{\leq N_0}\big\|_{L_t^\infty L_x^\frac{4p}{p+2}}\|u_{\leq N_0}\|_{L_t^\infty L_x^{2p}}^{p-1}\nonumber
\\
&\lesssim\eta^p N^{-1} \sum_{M\leq\frac{N}{10}}{M}^{-1+\frac{4}q+\frac2p}A_q(M).
\label{refines}
\end{align}
Comparing with \eqref{diedai1252}, we find that this term is also acceptable.

Inserting the estimates \eqref{di1}, \eqref{di2}, and \eqref{refines} into \eqref{dd11}, we conclude that \eqref{diedai1252} holds. This completes the proof of Proposition~\ref{adecay1}. \end{proof}

%%%%%%%%%%%%%%%%%%

\section{Long-time Strichartz estimates}\label{lts section} 

In this section we establish a long-time Strichartz estimate for almost periodic solutions to \eqref{equ1.1} as in Theorem~\ref{threenew}. Such estimates were originally developed by Dodson \cite{Dodson3} in the context of the mass-critical NLS. They have since appeared in the energy-critical setting in dimensions three and four (see \cite{KV3D, Visan2011}), as well as in the inter-critical setting (that is, $0<s_c<1$) in dimensions $d\in\{3,4,5\}$ (see \cite{Mu}). In this paper, we prove a long-time Strichartz estimate for the first time in the energy-supercritical regime. These estimates will be used in Section \ref{frequency-cascade}, in which we rule out the existence of rapid frequency-cascade solutions, as well as in Section \ref{flim section}, in which we establish a frequency-localized interaction Morawetz inequality.

\begin{theorem}[Long-time Strichartz estimate]\label{lse}
Let $1<s_c<3/2$ and let $u:[0,T_{max})\times\R^4\to\C $ be an almost
periodic solution to \eqref{equ1.1} with $N(t)\equiv N_k\geq1$ on
each characteristic $J_k\subset[0,T_{max}).$ Then on
any compact time interval $I\subset[0,T_{max}),$ which is a union of
characteristic subintervals $J_k$, and for any $N>0,$ we
have
\begin{equation}\label{lse1}
\big\||\nabla|^{s_c}u_{\leq
N}\big\|_{L_t^2L_x^4(I\times\R^4)}\lesssim_u1+N^{2s_c-1/2}K^{1/2},
\end{equation}
where $K:=\int_IN(t)^{3-4s_c}\,dt.$ Moreover, for any $\eta>0$, there
exists $N_0=N_0(\eta)$ such that for all $N\leq N_0,$
\begin{equation}\label{lsesmall}
\big\||\nabla|^{s_c}u_{\leq
N}\big\|_{L_t^2L_x^4(I\times\R^4)}\lesssim_u\eta\big(1+N^{2s_c-1/2}K^{1/2}\big).
\end{equation}
Furthermore, the implicit constants in \eqref{lse1} and
\eqref{lsesmall} are independent of $I$.
\end{theorem}

To begin, we fix $I\subset[0,T_{max})$ to be a union of contiguous characteristic subintervals. Throughout the proof, all spacetime norms will be taken over $I\times\R^4$ unless explicitly stated otherwise. For $N>0$, we define the quantities
$$B(N):=\big\||\nabla|^{s_c}u_{\leq
N}\big\|_{L_t^2L_x^4(I\times\R^4)}\quad\text{and}
\quad B_k(N):=\big\||\nabla|^{s_c}u_{\leq
N}\big\|_{L_t^2L_x^4(J_k\times\R^4)}.$$ 

We will prove Theorem~\ref{lse} by induction.  For the base case, we have the following

\begin{lemma}\label{base}
The estimate \eqref{lse1} holds for $N\geq N_{max}:=\sup\limits_{J_k\subset
I}N_k.$
\end{lemma}
\begin{proof}
This is a simple consequence of Lemma~\ref{fspssn}. Indeed, we have
\begin{align*}
\big\||\nabla|^{s_c}u_{\leq
N}\big\|_{L_t^2L_x^4}^2&\lesssim_u1+\int_I N(t)^2\,dt
 \\ &\lesssim_u1+N_{max}^{4s_c-1}\int_IN(t)^{3-4s_c}\,dt
 \\ &\lesssim_u 1+N^{4s_c-1}K,
\end{align*}
which gives \eqref{lse1}.
\end{proof}

To complete the induction, we will establish a recurrence relation for $B(N)$. To this end, we first let $\eps_0>0$ and $\eps>0$ be small parameters to be determined later. We use Remark~\ref{rem1.1} to find $c=c(\eps)$ so that 
\begin{equation}\label{xcx}
\big\||\nabla|^{s_c}u_{\leq cN(t)}\big\|_{L_t^\infty L_x^2}<
\varepsilon.
\end{equation}
The recurrence relation we will use takes the following form.

\begin{lemma}[Recurrence relation for $B(N)$]\label{rfthm}
\begin{equation}\label{rf1}
\begin{split}
B(N)&\lesssim_u\inf_{t\in I}\big\||\nabla|^{s_c}u_{\leq
N}(t)\big\|_{L_x^2}+C(\varepsilon,\varepsilon_0)N^{2s_c-1/2}K^{1/2}\\
&\quad\quad+\varepsilon^{1/2}
B\big(N/\varepsilon_0\big)+\sum_{M>N/\varepsilon_0}(\tfrac{N}{M})^{\frac53s_c}B(M)
\end{split}
\end{equation}
uniformly in $N$ for some positive constant
$C(\varepsilon,\varepsilon_0)$. 

We also have the following refinement of \eqref{rf1}:
\begin{equation}\label{rf2}
\begin{split}
B(N)\lesssim_ug(N)\big(1+N^{2s_c-1/2}K^{1/2}\big)+\varepsilon^{1/2}
B\big(N/\varepsilon_0\big)+\sum_{M>N/\varepsilon_0}(\tfrac{N}{M})^{\frac53s_c}B(M),
\end{split}
\end{equation}
where
\begin{equation}\label{g(N)}
g(N):=\inf_{t\in I}\big\||\nabla|^{s_c}u_{\leq N}(t)\big\|_{
L_x^2}+C(\varepsilon,\varepsilon_0)\sup_{J_k\subset
I}\big\||\nabla|^{s_c}u_{\leq N/\varepsilon_0}\big\|_{L_t^\infty
L_x^2(J_k\times\R^4)}.
\end{equation}
\end{lemma}

Before we turn to the proof of Lemma~\ref{rfthm}, let us see that we can use it to complete the proof of Theorem~\ref{lse}.

\begin{proof}[Proof of Theorem~\ref{lse}] From Lemma~\ref{base}, we see that \eqref{lse1} holds for
$N\geq N_{max}$. That is, we have
\begin{equation}\label{xhdd}
B(N)\leq C(u)\big(1+N^{2s_c-1/2}K^{1/2}\big),
\end{equation}
for all $N\geq N_{max}$. Clearly this
inequality remains true if we replace $C(u)$ by any larger constant.

We now suppose that \eqref{xhdd} holds for frequencies above $N$ and use the recurrence formula \eqref{rf1} to show that \eqref{xhdd} holds at frequency $N/2$. Choosing $\eps_0<1/2$, we use \eqref{rf1} and \eqref{xhdd} to see
\begin{align*}
B\big(N/2\big)
	&\leq \wt{C}(u)\Big(1+C(\varepsilon,\varepsilon_0)\big(N/2\big)^{2s_c-1/2}K^{1/2}+\varepsilon^{1/2} B\big(N/2\varepsilon_0\big)
		\\ &\quad\quad\quad\quad+\!\!\sum_{M>N/2\varepsilon_0}(\tfrac{N}{M})^{\frac53s_c}B(M)\Big)
	\\ &\leq\wt{C}(u)\Big(1+C(\varepsilon,\varepsilon_0)\big(N/2\big)^{2s_c-1/2}K^{1/2} 
	\\ 	&\quad\quad\quad\quad+\varepsilon^{1/2} C(u)\big(1+(N/2\varepsilon_0)^{2s_c-1/2}K^{1/2}\big)
	\\	&\quad\quad\quad\quad+C(u)\sum_{M>N/2\varepsilon_0}(\tfrac{N}{2M})^{\frac53s_c}\big(1+M^{2s_c-1/2}K^{1/2}\big)\Big)
	\\ &\leq\wt{C}(u)\Big(1+C(\varepsilon,\varepsilon_0)\big(N/2\big)^{2s_c-1/2}K^{1/2} 
	\\	&\quad\quad\quad\quad+\varepsilon^{1/2} C(u)\big(1+(N/2\varepsilon_0)^{2s_c-1/2}K^{1/2}\big)
	\\	&\quad\quad\quad\quad+C(u)\varepsilon_0^{\frac53s_c}+C(u)\varepsilon_0^{\frac13(\frac32-s_c)}(N/2)^{2s_c-\frac12}K^\frac12\Big)
	\\&=\wt{C}(u)\Big(1+C(\varepsilon,\varepsilon_0)\big(N/2\big)^{2s_c-1/2}K^{1/2}\Big)
	\\	&\quad\!+C(u)\Big(\big(\varepsilon^{1/2}+\varepsilon_0^{\frac53s_c}\big)\wt{C}(u)
	\\ &\quad\quad\quad\quad\quad+\big(\varepsilon_0^{1/2-2s_c}\varepsilon^{1/2}
+\varepsilon_0^{\frac13(\frac32-s_c)}\big)\wt{C}(u)(N/2)^{2s_c-1/2}K^{1/2}\Big),
\end{align*}
where we use $s_c<3/2$ in the third inequality to guarantee the
convergence of the sum. 

If we now choose $\varepsilon_0$ possibly even smaller depending on $\wt{C}(u)$; $\varepsilon$ sufficiently small depending on $\wt{C}(u)$ and $\varepsilon_0$; and $C(u)$ possibly larger such that $C(u)\geq2\big(1+C(\varepsilon,\varepsilon_0)\big)\wt{C}(u),$ we find
\begin{align*}
B\big(N/2\big)
		&\leq\wt{C}(u)\Big(1+C(\varepsilon,\varepsilon_0)(N/2)^{2s_c-1/2}K^{1/2}\Big)+\tfrac12C(u)\Big(1+(N/2)^{2s_c-1/2}K^{1/2}\Big)
		\\ &\leq C(u)\Big(1+(N/2)^{2s_c-1/2}K^{1/2}\Big).
\end{align*}
Thus \eqref{xhdd} holds at frequency ${N}/2$, and hence we conclude \eqref{lse1} by induction.

Next, we will use the recurrence formula \eqref{rf2} to prove \eqref{lsesmall}. To begin, we note that for fixed $\eps,\eps_0>0$, we can use Remark~\ref{rem1.1} and the fact that $\inf_{t\in I}N(t)\geq1$ to see
		\begin{equation}\label{g(N) to zero}
		\lim_{N\to0}g(N)=0,
		\end{equation}
where $g(N)$ is as in \eqref{g(N)}. 

Now, arguing as above, we can use \eqref{lse1} and \eqref{rf2} to see
\begin{align*}
B(N)& \leq\wt{C}(u)\Big(g(N)+g(N)N^{2s_c-1/2}K^{1/2}+\varepsilon^{1/2}\big(1+(N/\varepsilon_0)^{2s_c-1/2}K^{1/2}\big)
	\\&\quad\quad\quad\quad+\varepsilon_0^{\frac53s_c}\big(1+(N/\varepsilon_0)^{2s_c-1/2}K^{1/2}\big)\Big)
	\\&\leq\wt{C}(u)\Big(g(N)+\varepsilon^{1/2}+\varepsilon_0^{\frac53s_c}
	\\&\quad\quad\quad\quad+\big(g(N)+\varepsilon^{1/2}\varepsilon_0^{1/2-2s_c}
	+\varepsilon_0^{\frac13(\frac32-s_c)}\big)N^{2s_c-1/2}K^{1/2}\Big).
\end{align*}

Now let $\eta>0$. We first choose $\varepsilon_0$ small enough depending on $\eta$ so that
$$\varepsilon_0^{\frac53s_c}+\varepsilon_0^{\frac13(\frac32-s_c)}<\eta.$$ 
Next, we choose $\varepsilon$ sufficiently small depending on $\eta$ and $\varepsilon_0$ so that
$$\varepsilon^{1/2}+\varepsilon^{1/2}\varepsilon_0^{1/2-2s_c}<\eta.$$ 
Finally, we choose $C(u)\geq 2\tilde{C}(u).$ Recalling \eqref{g(N) to zero}, we can now choose $N_0=N_0(\eta)$ such that $g(N)<\eta$ for $N\leq N_0$. Thus, continuing from above, we see that \eqref{lsesmall} holds for for $N\leq N_0$. This concludes the proof of Theorem~\ref{lse}.
\end{proof}

It remains to prove Lemma~\ref{rfthm}.

\begin{proof}[Proof of Lemma~\ref{rfthm}] We begin with an application of Strichartz to see
\begin{equation}\label{nte}
B(N)\lesssim \inf_{t\in I}\big\|{\nsc u_{\leq N}(t)}\big\|_{L_x^2}+\big\||\nabla|^{s_c}P_{\leq
N}F(u)\big\|_{L_t^2L_x^{4/3}}.
\end{equation}

We need only to estimate the nonlinear term. To do this, we will decompose the nonlinearity and estimate each piece individually. Recalling that $p>2$, we write

\begin{align}
F(u)	&=G(u)u_{>N/\varepsilon_0}			\label{fu1}
\\	&\quad+|u|^{p-2}\bar{u}\big(P_{\leq cN(t)}u_{\leq N/\varepsilon_0}\big)u_{\leq N/\varepsilon_0} \label{fu3}
\\	&\quad+|u|^{p-2}\bar{u}\big(P_{>cN(t)}u_{\leq N/\varepsilon_0}\big)u_{\leq N/\varepsilon_0},	\label{fu2}
\end{align}
where 
			$$G(u):=|u|^p+|u|^{p-2}\bar{u}u_{\leq N/\varepsilon_0}$$ and 
$c=c(\eps)$ is as in \eqref{xcx}.

We first consider the contribution of term \eqref{fu1}. Using H\"older, Bernstein, and
Lemma~\ref{pde}, we estimate
\begin{align}\nonumber
\big\||\nabla|^{s_c}P_{\leq N}&\big(G(u)u_{>N/\varepsilon_0}\big)\big\|_{L_t^2L_x^{4/3}}		
		\\ 	&\lesssim	\nonumber
N^{\frac53s_c}\big\||\nabla|^{-\frac23s_c}\big(G(u)u_{>N/\varepsilon_0}\big)\big\|_{L_t^2L_x^{4/3}}
		\\    \nonumber
			&\lesssim N^{\frac53s_c}\big\||\nabla|^{\frac23s_c}G(u)\big\|_{L_t^\infty L_x^\frac{6p}{5p-2}}\big\||\nabla|^{-\frac23s_c}u_{>N/\varepsilon_0}\big\|_{L_t^2L_x^4}
		\\	\label{xdd1}
			&\lesssim\big\||\nabla|^{\frac23s_c}G(u)\big\|_{L_t^\infty L_x^\frac{6p}{5p-2}}\sum_{M>N/\varepsilon_0}(\tfrac{N}{M})^{\frac53s_c}B(M).
\end{align}
On the other hand, by the fractional chain rule, fractional product rule, Sobolev embedding, and \eqref{assume1.1}, we can estimate
$$\big\||\nabla|^{\frac23s_c}|u|^p\big\|_{L_t^\infty
L_x^\frac{6p}{5p-2}}\lesssim \|u\|_{L_t^\infty
L_x^{2p}}^{p-1}\big\||\nabla|^{\frac23s_c}u\big\|_{L_t^\infty
L_x^\frac{6p}{2p+1}}\lesssim\big\||\nabla|^{s_c}u\big\|_{L_t^\infty
L_x^2}^p\lesssim_u1$$
 and
\begin{align*}
\big\||\nabla|^{\frac23s_c}&\big(|u|^{p-2}\bar{u}u_{\leq
N/\varepsilon_0}\big)\big\|_{L_t^\infty L_x^\frac{6p}{5p-2}}
	\\ &\lesssim\|u\|_{L_t^\infty L_x^{2p}}\big\||\nabla|^{\frac23s_c}\big(|u|^{p-2}\bar{u}\big)
		\big\|_{L_t^\infty L_x^\frac{6p}{5(p-1)}}+\|u\|_{L_t^\infty L_x^{2p}}^{p-1}
		\big\||\nabla|^{\frac23s_c}u\big\|_{L_t^\infty L_x^\frac{6p}{2p+1}}
	\\& \lesssim\big\||\nabla|^{s_c}u\big\|_{L_t^\infty L_x^2}\|u\|_{L_t^\infty L_x^{2p}}^{p-2}
		\big\||\nabla|^{\frac23s_c}u\big\|_{L_t^\infty L_x^\frac{6p}{2p+1}}+\big\||\nabla|^{s_c}u\big\|_{L_t^\infty L_x^2}^p
\\ &\lesssim\big\||\nabla|^{s_c}u\big\|_{L_t^\infty L_x^2}^p
\\ &\lesssim_u1.
\end{align*}
Plugging these to \eqref{xdd1}, we obtain
\begin{equation}\label{diyg}
\big\||\nabla|^{s_c}P_{\leq
N}\big(G(u)u_{>N/\varepsilon_0}\big)\big\|_{L_t^2L_x^{4/3}}\lesssim\sum_{M>N/\varepsilon_0}\big(\tfrac{N}{M}\big)^{\frac53s_c}B(M).
\end{equation}
\begin{remark}\label{restriction 1}
It is this term that gives rise to the restriction $s_c<3/2$. We can see this in the following way. Our approach to this estimating this term (namely, using the paraproduct estimate Lemma~\ref{pde}) yields a bound of the form $$\sum_{M>N/\eps_0}\big(\tfrac{N}{M}\big)^{s_c+s}B(M)$$ for some $s>0$. To carry out the induction, we need this sum to converge (as we saw above). As we are trying to prove $B(N)\lesssim_u 1+N^{2s_c-1/2}K^{1/2}$, this means we need to have $s>s_c-1/2$. On the other hand, one can check that the scaling constraints from the paraproduct estimate (Lemma~\ref{pde}) impose the condition $s<1$. Hence, our approach only allows us to estimate this term satisfactorily for $s_c<3/2$. In particular, we do so by choosing $s=\tfrac23 s_c$.

We note here that in three dimensions, the issue just discussed arises when $s_c=1$, that is, when the equation is energy-critical. The authors of \cite{KV3D}  devised a new strategy to establish a long-time Strichartz estimate in this setting. Specifically, they proved a `maximal' Strichartz estimate, allowing for control over the worst possible Littlewood--Paley piece at each point in time. See \cite{KV3D} for a further discussion.
\end{remark}

Next, we turn to estimating the contribution of the term \eqref{fu3} to \eqref{nte}. Using the fractional product rule, we estimate

\begin{align}\nonumber
\Big\||\nabla|^{s_c}&P_{\leq N}\Big[|u|^{p-2}\bar{u}\big(P_{\leq
cN_k}u_{\leq N/\varepsilon_0}\big)u_{\leq
N/\varepsilon_0}\Big]\Big\|_{L_t^2L_x^{4/3}}
	\\ \nonumber &\lesssim\big\||\nabla|^{s_c}\big(|u|^{p-2}\bar{u}\big)\big\|_{L_t^\infty L_x^\frac{p}{p-1}}
			\big\|P_{\leq cN_k}u_{\leq N/\varepsilon_0}\big\|_{L_t^4L_x^\frac{8p}{4-p}}\|u_{\leq N/\varepsilon_0}
				\big\|_{L_t^4L_x^\frac{8p}{4-p}}
	\\ \nonumber &\quad+\|u\|_{L_t^\infty L_x^{2p}}^{p-1}\big\||\nabla|^{s_c}P_{\leq cN_k}u_{\leq N/\varepsilon_0}\big\|_{L_t^4L_x^{8/3}}
				\|u_{\leq N/\varepsilon_0}\big\|_{L_t^4L_x^\frac{8p}{4-p}}
	\\ &\quad+\|u\|_{L_t^\infty L_x^{2p}}^{p-1}\big\|P_{\leq
cN_k}u_{\leq N/\varepsilon_0}\big\|_{L_t^\infty L_x^{2p}}\big\||\nabla|^{s_c}u_{\leq
N/\varepsilon_0}\big\|_{L_t^2L_x^4}. \label{xdd3}
\end{align}
Using the fractional chain rule, Sobolev embedding, and \eqref{assume1.1}, we find
$$
\big\||\nabla|^{s_c}\big(|u|^{p-2}\bar{u}\big)\big\|_{L_t^\infty
L_x^\frac{p}{p-1}}\lesssim\|u\|_{L_t^\infty
L_x^{2p}}^{p-2}\big\||\nabla|^{s_c}u\big\|_{L_t^\infty
L_x^2}\lesssim_u1.
$$
By Sobolev embedding, interpolation, and \eqref{xcx}, we also see
\begin{align*}
\big\|P_{\leq cN_k}u_{\leq N/\varepsilon_0}\big\|_{L_t^4L_x^\frac{8p}{4-p}}&\lesssim\big\||\nabla|^{s_c}P_{\leq cN_k}u_{\leq N/\varepsilon_0}\big\|_{L_t^4L_x^{8/3}}
\\ &\lesssim\big\||\nabla|^{s_c}P_{\leq cN_k}u_{\leq N/\varepsilon_0}\big\|_{L_t^\infty L_x^2}^{1/2}\big\||\nabla|^{s_c}P_{\leq cN_k}u_{\leq N/\varepsilon_0}\big\|_{L_t^2L_x^4}^{1/2}
\\ &\lesssim\varepsilon^{1/2} B\big(N/\varepsilon_0\big)^{1/2}.
\end{align*}
Similarly, we find
$$\|u_{\leq N/\varepsilon_0}\big\|_{L_t^4L_x^\frac{8p}{4-p}}\lesssim B\big(N/\varepsilon_0\big)^{1/2}.$$ 
Plugging these estimates into \eqref{xdd3} and applying a few more instances of Sobolev embedding, we obtain
\begin{equation}\label{xdd31}
\Big\||\nabla|^{s_c}P_{\leq N}\Big[|u|^{p-2}\bar{u}\big(P_{\leq
cN_k}u_{\leq N/\varepsilon_0}\big)u_{\leq
N/\varepsilon_0}\Big]\Big\|_{L_t^2L_x^{4/3}}\lesssim_u
\varepsilon^{1/2} B\big(N/\varepsilon_0\big).
\end{equation}

Finally, we estimate the contribution of term \eqref{fu2} to \eqref{nte}.  To begin, we restrict our attention to a single characteristic interval $J_k$. Note that it suffices to consider the case $cN_k\leq {N}/{\varepsilon_0}$. In this case, using H\"older, Bernstein, Sobolev embedding, interpolation, Lemma~\ref{fspssn}, and \eqref{assume1.1}, we find
\begin{align}\nonumber
\Big\||\nabla|^{s_c}&P_{\leq N}\Big[|u|^{p-2}\bar{u}\big(P_{>cN_k}u_{\leq N/\varepsilon_0}\big)u_{\leq N/\varepsilon_0}\Big]\Big\|_{L_t^2L_x^{4/3}(J_k\times\R^4)}
		\\ \nonumber &\lesssim N^{s_c}\Big\||u|^{p-2}\bar{u}\big(P_{>cN_k}u_{\leq N/\varepsilon_0}\big)u_{\leq N/\varepsilon_0}\Big\|_{L_t^2L_x^{4/3}(J_k\times\R^4)}
		\\ \nonumber &\lesssim N^{s_c}\|u\|_{L_t^\infty L_x^{2p}}^{p-1}\big\|P_{>cN_k}u_{\leq N/\varepsilon_0}
			\big\|_{L_t^4L_x^{8/3}(J_k\times\R^4)}\|u_{\leq N/\varepsilon_0}\|_{L_t^4L_x^\frac{8p}{4-p}(J_k\times\R^4)}
		\\ \nonumber
			&\lesssim N^{s_c}\big(cN_k\big)^{-s_c}\big\||\nabla|^{s_c}u_{\leq N/\eps_0}\big\|_{L_t^4L_x^{8/3}(J_k\times\R^4)}^2
		\\ \nonumber
		 	&\lesssim c^{-s_c}(\tfrac{N}{N_k})^{s_c}\big\||\nabla|^{s_c}u_{\leq N/\eps_0}\big\|_{L_t^\infty L_x^2(J_k\times\R^4)}\big\||\nabla|^{s_c}u\big\|_{L_t^2L_x^4(J_k\times\R^4)}
		\\ \nonumber &\lesssim_u C(\varepsilon,\eps_0)(\tfrac{N}{N_k})^{2s_c-1/2}.
\end{align}

Summing these estimates over $J_k\subset I$ now yields
\begin{equation}\label{xdd21}
\Big\||\nabla|^{s_c}P_{\leq
N}\Big[|u|^{p-2}\bar{u}\big(P_{>cN_k}u_{\leq
N/\varepsilon_0}\big)u_{\leq
N/\varepsilon_0}\Big]\Big\|_{L_t^2L_x^{4/3}(I\times\R^4)}\lesssim_u
C(\varepsilon,\varepsilon_0)N^{2s_c-1/2}K^{1/2}.
\end{equation}

We note that here that in the estimates above, we could have held on to the terms
		$$\xnorm{\nsc u_{\leq N/\eps_0}}{\infty}{2}{J_k}$$
on each $J_k$. In this case, summing the estimates over $J_k\subset I$ yields 

\begin{align}\nonumber
 \Big\||\nabla|^{s_c}&P_{\leq N}\Big[|u|^{p-2}\bar{u}\big(P_{>cN_k}u_{\leq N/\varepsilon_0}\big)u_{\leq N/\varepsilon_0}\Big]\Big\|_{L_t^2L_x^{4/3}}
 \\ &\lesssim_u
C(\varepsilon,\varepsilon_0)\sup_{J_k\subset
I}\big\||\nabla|^{s_c}u_{\leq N/\varepsilon_0}\big\|_{L_t^\infty L_x^2(J_k\times\R^4)}N^{2s_c-1/2}K^{1/2}.\label{xdd22}
\end{align}

Combining \eqref{nte} with \eqref{diyg}, \eqref{xdd31} and \eqref{xdd21}, we see that $B(N)$ satisfies the recurrence relation
\eqref{rf1}. If we use \eqref{xdd22} instead of \eqref{xdd21}, we can deduce the recurrence relation \eqref{rf2}. This completes the proof of Lemma~\ref{rfthm}.
\end{proof}

\section{The rapid frequency-cascade scenario}\label{frequency-cascade}

%%%%%%%%%%%%%%%%%%%%%%%%%%%%%%%%%%%%%%%%%%%%%%%%%%%%%%%%%%%%%%%%%%%

In this section, we preclude the existence of rapid frequency-cascade solutions, that is, almost periodic solutions as in Theorem~\ref{threenew} such that $\int_0^{T_{max}} N(t)^{3-4s_c}\,dt<\infty.$ The proof will rely primarily on the long-time Strichartz estimate proved in the previous section.

\begin{theorem}[No rapid frequency-cascades]\label{no frequency-cascades} Let $1<s_c<3/2.$ Then there are no almost periodic solutions $u:[0,T_{max})\times\R^4\to\C$ to \eqref{equ1.1} with $N(t)\equiv N_k\geq 1$ on each characteristic subinterval $J_k\subset[0,T_{max})$ that satisfy
	\begin{equation}\label{rf blowup}
	\xnorms{u}{3p}{[0,T_{max})}=\infty
	\end{equation}
and
	\begin{equation}\label{rf K}
	K:=\int_0^{T_{max}} N(t)^{3-4s_c}\,dt<\infty.
	\end{equation}
\end{theorem}

\begin{proof} We argue by contradiction. Suppose that $u$ were such a solution. Then, using \eqref{rf K} and Corollary~\ref{constancy cor}, we see 
\begin{equation}\label{N at blowup}
\lim_{t\to T_{max}}N(t)=\infty,
\end{equation}
whether $T_{max}$ is finite or infinite. Combining this with \eqref{xiaoc}, we see that 
\begin{equation}\label{gtx}
\lim_{t\to T_{max}}\big\||\nabla|^{s_c}u_{\leq
N}(t)\big\|_{L_x^2(\R^4)}=0
\end{equation}
for any $N>0$. 

We now let $I_n\subset [0,T_{max})$ be a nested sequence of compact time intervals, each of which is a contiguous union of characteristic subintervals. We claim that for any $N>0$, we have
\begin{equation}\label{gxx}
\big\||\nabla|^{s_c}u_{\leq N}\big\|_{L_t^2L_x^4(I_n\times\R^4)}\lesssim_u\inf_{t\in I_n}\big\||\nabla|^{s_c}u_{\leq N}(t)\big\|_{L_x^2}+N^{2s_c-1/2}.%K^{1/2}.
\end{equation}
Indeed, defining
		$$B_n(N):=\big\||\nabla|^{s_c}u_{\leq N}\big\|_{L_t^2L_x^4(I_n\times\R^4)},$$
we have by \eqref{rf1} and \eqref{rf K} the estimate
\begin{align*}
B_n(N) \lesssim_u\inf_{t\in I_n}\big\||\nabla|^{s_c}u_{\leq N}(t)\big\|_{L_x^2}+C(\varepsilon,\varepsilon_0)N^{2s_c-1/2}+\sum_{M>N/\varepsilon_0}(\tfrac{N}{M})^{\frac53s_c}B_n(M).
\end{align*}
Arguing as we did to obtain \eqref{lse1}, we derive \eqref{gxx}. 

Now, letting $n\to\infty$ in \eqref{gxx} and using \eqref{gtx}, we find
\begin{equation}\label{rfc11}
\big\||\nabla|^{s_c}u_{\leq
N}\big\|_{L_t^2L_x^4([0,T_{max})\times\R^4)}\lesssim_u
N^{2s_c-1/2}%K^{1/2}
\end{equation}
for all $N>0$. 

We now claim that \eqref{rfc11} implies
\begin{lemma}\label{rfclemma}
\begin{equation}\label{ngtcl}
\big\||\nabla|^{s_c}u_{\leq N}\big\|_{L_t^\infty
L_x^2([0,T_{max})\times\R^4)}\lesssim_u N^{2s_c-1/2}%K^{1/2}
\end{equation}
for all $N>0.$ 
\end{lemma}

We will prove \eqref{ngtcl} below; for now, let us assume that \eqref{ngtcl} holds and use it to derive a contradiction. Fixing $0<\alpha<s_c-1/2$, we can use Bernstein, \eqref{assume1.1}, and \eqref{ngtcl} to estimate
\begin{align*}
\|\vert\nabla\vert^{-\alpha} u\|_{L_t^\infty L_x^2}&\lesssim\|\vert\nabla\vert^{-\alpha} u_{\leq1}\|_{L_t^\infty L_x^2}
		+\|\vert\nabla\vert^{-\alpha}u_{>1}\|_{L_t^\infty L_x^2}
	\\& \lesssim\sum_{N\leq1}\big\||\nabla|^{-\al}u_N\big\|_{L_t^\infty L_x^2}+\|\vert\nabla\vert^{s_c} u_{>1}\|_{L_t^\infty L_x^2}
	\\& \lesssim_u\sum_{N\leq1}N^{-\alpha-s_c}\big\||\nabla|^{s_c}u_N\big\|_{L_t^\infty L_x^2}+1
	\\ &\lesssim_u \sum_{N\leq1}N^{s_c-1/2-\alpha}+1
	\\ &\lesssim_u 1.
\end{align*}
Hence $u\in L_t^\infty\dot{H}_x^{-\alpha}([0,T_{max})\times\R^4)$. For $\eta>0$, we can interpolate this bound with \eqref{xiaoc} to get
\begin{equation}\nonumber
\int_{|\xi|\leq c(\eta) N(t)}\big|\wh{u}(t,\xi)\big|^2d\xi\lesssim_u
\eta^\frac{\al}{s_c+\al}.
\end{equation}
Thus, by Plancherel we find
\begin{align*}
M[u_0]=M[u(t)]&=\int_{|\xi|\leq
c(\eta) N(t)}\big|\wh{u}(t,\xi)\big|^2d\xi+\int_{|\xi|>
c(\eta) N(t)}\big|\wh{u}(t,\xi)\big|^2d\xi\\
&\lesssim_u\eta^\frac{\al}{s_c+\al}+\big(c(\eta)N(t)\big)^{-2s_c}\big\||\nabla|^{s_c}u\big\|_{L_t^\infty L_x^2}^2\\
&\lesssim_u\eta^\frac{\al}{s_c+\al}+\big(c(\eta)N(t)\big)^{-2s_c}.
\end{align*}
Choosing $\eta$ small, sending $t\to T_{max}$, and recalling \eqref{N at blowup}, we can deduce that $M[u_0]=0$. Thus, we have $u\equiv 0$, which contradicts \eqref{rf blowup}. \end{proof}

It remains to prove Lemma~\ref{rfclemma}.
 
\begin{proof}[Proof of Lemma~\ref{rfclemma}]
We first use the reduced Duhamel formula and Strichartz to see
\begin{align}
\big\||\nabla|^{s_c}u_{\leq N}\big\|_{L_t^\infty
L_x^2([0,T_{max})\times\R^4)}\lesssim
\label{gfues}
\big\||\nabla|^{s_c}P_{\leq N}F(u)\big\|_{L_t^2L_x^{4/3}([0,T_{max})\times\R^4)}.
\end{align}
We now decompose the nonlinearity via
\begin{align}\label{fgu}
F(u)=|u|^{p-2}\bar{u}u_{\leq
N}^2+\big(|u|^{p-2}\bar{u}u_{>N}+2|u|^{p-2}\bar{u}u_{\leq
N}\big)u_{>N}
\end{align}
and estimate the contribution of each piece individually.

We begin by estimating the contribution of the first term in \eqref{fgu}
to \eqref{gfues}. Using H\"older, the fractional product rule, the fractional chain rule, Sobolev embedding, interpolation, \eqref{assume1.1}, and \eqref{rfc11}, we estimate
\begin{align*}
\big\||\nabla|^{s_c}&P_{\leq N}\big(|u|^{p-2}\bar{u}u_{\leq
N}^2\big)\big\|_{L_t^2 L_x^{4/3}}
\\
&\lesssim\big\||\nabla|^{s_c}\big(|u|^{p-2}\bar{u}\big)\big\|_{L_t^\infty L_x^\frac{p}{p-1}}\|u_{\leq N}\|_{L_t^4 L_x^\frac{8p}{4-p}}^2+\|u\|_{L_t^\infty L_x^{2p}}^{p-1}\big\||\nabla|^{s_c}(u_{\leq N}^2)\big\|_{L_t^2L_x^\frac{4p}{p+2}}
\\
&\lesssim\|u\|_{L_t^\infty L_x^{2p}}^{p-2}\big\||\nabla|^{s_c}u\big\|_{L_t^\infty L_x^2}\big\||\nabla|^{s_c}u_{\leq
N}\big\|_{L_t^4L_x^{8/3}}^2
	\\&\quad\quad+\|u\|_{L_t^\infty L_x^{2p}}^{p-1}\|u_{\leq N}\|_{L_t^\infty L_x^{2p}}\big\||\nabla|^{s_c}u_{\leq N}\big\|_{L_t^2L_x^4}
\\
&\lesssim\big\||\nabla|^{s_c}u\big\|_{L_t^\infty L_x^2}^{p}\big\||\nabla|^{s_c}u_{\leq N}\big\|_{L_t^2 L_x^4}
\\
&\lesssim_uN^{2s_c-1/2}.
\end{align*}

Next, we estimate the contribution of the second term in
\eqref{fgu} to \eqref{gfues}. Defining
		$$G(u):=\vert u\vert^{p-2}\bar{u}u_{>N}+2\vert u\vert^{p-2}\bar{u}u_{\leq
N},$$
we can use H\"older, Bernstein, Lemma~\ref{pde}, and \eqref{rfc11} to estimate
\begin{align*}
\big\||\nabla|^{s_c}&P_{\leq N}\big(G(u)u_{>N}\big)\big\|_{L_t^2
L_x^{4/3}} 
\\ &\lesssim
N^{\frac53s_c}\big\||\nabla|^{-\frac23s_c}\big(G(u)u_{>N}\big)\big\|_{L_t^2
L_x^{4/3}}
\\
&\lesssim N^{\frac53s_c}\big\||\nabla|^{\frac23s_c}G(u)\big\|_{L_t^\infty
L_x^\frac{6p}{5p-2}}\big\||\nabla|^{-\frac23s_c}u_{>N}\big\|_{L_t^2L_x^4}
\\
&\lesssim \big\||\nabla|^{\frac23s_c}G(u)\big\|_{L_t^\infty
L_x^\frac{6p}{5p-2}}\sum_{M>N}(\tfrac{N}{M})^{\frac53s_c}\big\||\nabla|^{s_c}u_M\big\|_{L_t^2L_x^4}
\\
&\lesssim_u N^{2s_c-1/2}\big\||\nabla|^{\frac23s_c}G(u)\big\|_{L_t^\infty
L_x^\frac{6p}{5p-2}}.
\end{align*}
Using the fractional chain rule, the fractional product rule, and Sobolev embedding, we can estimate the remaining term via
\begin{align*}
\big\||\nabla|^{\frac23s_c}G(u)\big\|_{L_t^\infty L_x^\frac{6p}{5p-2}}
\lesssim\big\||\nabla|^{s_c}u\big\|_{L_t^\infty L_x^{2p}}^{p}\lesssim_u 1.
\end{align*}

Plugging our estimates into \eqref{gfues}, we conclude that \eqref{ngtcl} holds, which completes the proof of Lemma~\ref{rfclemma}.
\end{proof}

%%%%%%%%%%%%%%%%%%%%%%%%%%%%%%%%%%%%%%%%%%%%%%%%%%%%%%%%%%%%%%%%%%%
\section{The frequency-localized interaction Morawetz inequality}\label{flim section}
In this section, we establish spacetime bounds for the high-frequency portions of almost periodic solutions to \eqref{equ1.1}. We will use these estimates in the next section to preclude the existence of quasi-soliton solutions in the sense of Theorem~\ref{threenew}.

\begin{theorem}[Frequency-localized interaction Morawetz inequality]\label{uppbd} Let $1<s_c<3/2$, and let $u:[0,T_{max})\times\R^4\to\C$ be an almost periodic solution to \eqref{equ1.1} such that $N(t)\equiv N_k\geq 1$ on each characteristic subinterval $J_k\subset[0,T_{max})$. Let $I\subset[0,T_{max})$ be a compact time interval, which is a contiguous union of characteristic subintervals $J_k$. Then for any $\eta>0$, there exists $N_0=N_0(\eta)$ such that for all $N\leq N_0$, we have
\begin{equation}\label{limemd}
\int_I\iint_{\R^4\times\R^4}\frac{\vert u_{> N}(t,y)\vert^2\vert u_{>
N}(t,x)\vert^2}{\vert x-y\vert^3}dx\,dy\,dt \lesssim_u\eta\big(N^{1-4s_c}+K\big),
\end{equation}
where $K:=\int_IN(t)^{3-4s_c}dt$. Furthermore, $N_0$ and the implicit constants above are independent of the interval $I$.
\end{theorem}

We begin with a general form of the interaction Morawetz identity, introduced originally in \cite{ckstt}. For a fixed weight $a:\R^4\to\R$ and a function $\varphi$ solving $$(i\partial_t+\Delta)\varphi=\mathcal{N},$$ we define the \emph{interaction Morawetz action} by
				\begin{equation}
				\label{morawetz action}
				M(t)=2\Im\iint_{\R^4\times\R^4}\vert\varphi(y)\vert^2a_k(x-y)(\varphi_k\bar{\varphi})(x)\,dx\,dy,
				\end{equation}
where subscripts denote spatial derivatives and repeated indices are summed. 

Defining the \emph{mass bracket} $$\{f,g\}_M:=\Im(f\bar{g})$$ and the \emph{momentum bracket} $$\{f,g\}_{P}:=\Re(f\nabla\bar{g}-g\nabla\bar{f}),$$ one can compute
\begin{proposition}[Interaction Morawetz identity]\label{morawetz identity}
			\begin{align}
			\partial_t&M(t)\nonumber
			\\=&-\iint\vert\varphi(y)\vert^2a_{jjkk}(x-y)\vert\varphi(x)\vert^2\,dx\,dy\label{mass-mass}
			\\&+\iint\! 4a_{jk}(x-y)\big[\vert\varphi(y)\vert^2\Re(\bar{\varphi}_k\varphi_j)(x)-\Im(\bar{\varphi}\varphi_j)(y)
						\!\Im(\bar{\varphi}\varphi_k)(x)\big]dx\,dy\label{scary}
			\\&+\iint\{\mathcal{N},\varphi\}_M(y)4a_k(x-y)\Im(\bar{\varphi}\varphi_k)(x)\,dx\,dy\label{mass bracket}
			\\&+\iint\vert\varphi(y)\vert^2\ 2\nabla a(x-y)\cdot\{\mathcal{N},\varphi\}_P(x)\,dx\,dy.\label{momentum bracket}
			\end{align}
\end{proposition}

It is well-known that if one chooses the weight $a(x)=\vert x\vert$ and the function $\varphi$ to be a solution to \eqref{equ1.1}, then one can use Proposition~\ref{morawetz identity} and the fundamental theorem of calculus to deduce the \emph{interaction Morawetz inequality}:
	$$\int_I\iint_{\R^4\times\R^4}\frac{\vert \varphi(t,y)\vert^2\vert \varphi(t,x)\vert^2}{\vert x-y\vert^3}dx\,dy\,dt\lesssim \xnorm{\varphi}{\infty}{2}{I}^3\xnorm{\nabla\varphi}{\infty}{2}{I}.$$

In our setting, however, we only assume control over the $\dot{H}_x^{s_c}$-norm of solutions to \eqref{equ1.1}. Hence, to make the right-hand side of this inequality finite, we need to truncate solutions to high frequencies. Such a truncation results in error terms that must then be handled in order to recover a useful estimate. In certain situations, it is possible to control all of the resulting error terms and derive a useful estimate using the weight $a(x)=\vert x\vert$ (see \cite{Mu, RV, Visanphd, Visan2007, Visan2011} for some examples). In other settings, however, it turns out that certain error terms cannot be controlled. This problem was first encountered in \cite{CKSTT07} in the context of the energy-critical NLS in three dimensions. The authors of \cite{CKSTT07} addressed this difficulty by truncating the weight $a$ in space. While this additional truncation results in even more error terms to estimate, it ultimately proved to be a winning strategy (at least when combined with an averaging argument). The authors of \cite{KV3D} (who revisited the result of \cite{CKSTT07}) were also able to succeed in this setting via spatial truncation; in fact, they were able to avoid the averaging argument present in \cite{CKSTT07} altogether by means of a careful design for the weight $a$.

In order to establish the frequency-localized interaction Morawetz for $1<s_c<3/2$, we will also need to truncate the weight $a$ in space. It is actually possible to move slightly beyond $s_c=1$ using the standard weight $a(x)=\vert x\vert$, but this spatial truncation is necessary to treat the full range $1<s_c<3/2$. Because we are in the energy-supercritical setting, we will not actually need to design the weight as carefully as the authors of \cite{KV3D} (see the discussions preceding Lemma~\ref{nssc} and Lemma~\ref{potential energy}); nonetheless, we will model our presentation largely after \cite[Section~6]{KV3D} and make use of many of the ideas contained therein. 

We choose our weight $a:\R^4\to\R$ to be a smooth, spherically-symmetric function, which we may regard either as function of $x$ or $r=\vert x\vert$. We choose $a$ so that it satisfies the following properties:
	\begin{equation}\label{properties of a}
	\left\{
	\begin{array}{ll}
		\bullet\quad a(x)=\vert x\vert\quad\text{for}\quad \vert x\vert\leq N^{-1},
		\vspace{.05in}
	\\ \vspace{.05in}	\bullet\quad a(x)\text{ is constant for }\vert x\vert>2N^{-1},
	\\ 	\bullet\quad \partial_ra\geq 0,\vspace{.05in}
	\\ \bullet \quad\vert{\partial_r^k \partial_ra}\vert\lesssim_k r^{-k}\quad\text{for }k\geq 0,
	\end{array}
	\right.
	\end{equation}
where $\partial_r$ denotes the radial derivative and $N$ is a parameter to be chosen shortly.

To prove Theorem~\ref{uppbd}, we will apply Proposition~\ref{morawetz identity} with $a(x)$ as above. We will take $\varphi=u_{>N}$ and $\mathcal{N}=P_{>N}(F(u))$, where $N$ is taken small enough that $u_{>N}$ captures `most' of the solution $u$. To make this idea more precise, we will need the following corollary of Theorem~\ref{lse} and Proposition~\ref{mxs}. 

%%%corollary%%%

\begin{corollary}[A priori bounds]\label{a priori bounds} Let $1<s_c<3/2$, and let $u:[0,T_{max})\times\R^4\to\C$ be an almost periodic solution to \eqref{equ1.1} with $N(t)\equiv N_k\geq 1$ on each characteristic subinterval $J_k\subset[0,T_{max})$. Let $I\subset[0,T_{max})$ be a compact time interval, which is a contiguous union of characteristic subintervals $J_k$. Setting
		$K:=\int_I N(t)^{3-4s_c}\,dt,$ 
we have the following:

First, for any frequency $N>0$, we have 
					\begin{equation}\label{hc}
					\xnorm{\vert\nabla\vert^{\eps}u_{>N}}{q}{\frac{2q}{q-1}}{I}\lesssim_u N^{\eps-s_c}(1+N^{4s_c-1}K)^{1/q}
					\end{equation}
for all $2\leq q\leq\infty$ and $\eps<s_c-\tfrac{4s_c-1}{q}$ (note that $\eps<0$ is permitted). 

Second, for any $\eta>0$, there exists $N_0=N_0(\eta)$ such that for all $N\leq N_0$, we have
					\begin{align}
				\xnorm{\nsc u_{\leq N}}{q}{\frac{2q}{q-1}}{I}&\lesssim_u \eta(1+N^{4s_c-1}K)^{1/q}\quad\text{for}\quad 2\leq q\leq\infty,\label{lc}
					\\ \xnorm{u_{>N}}{\infty}{2}{I}&\lesssim_u \eta N^{-s_c},\label{hc2}
					\\ \xnorm{\nabla u_{>N}}{\infty}{2}{I}&\lesssim_u \eta N^{1-s_c} \label{hc3}.
					\end{align}
					
Finally, for any $\rho\geq N^{-1}$, we have
			\begin{align}
			\int_I&\sup_{x\in\R^4}\int_{\vert x-y\vert\leq\rho}\vert u_{>N}(t,y)\vert^2\,dy\,dt\lesssim_u \rho^2N^{-2s_c}(1+N^{4s_c-1}K).\label{mass control}
			\end{align}

Moreover, $N_0$ and the implicit constants above do not depend on $I$. 					

\end{corollary}

\begin{proof} Throughout the proof, all spacetime norms will be taken over $I\times\R^4$.

We first prove \eqref{hc}. To begin, we fix $\alpha>s_c-1/2$ and note that by \eqref{lse1} and Bernstein, we have
	\begin{align}
	\xonorm{\vert\nabla\vert^{-\alpha}u_{>N}}{2}{4}
	&\lesssim \sum_{M> N}M^{-\alpha-s_c}\xonorm{\nsc u_M}{2}{4} \nonumber
	\\ &\lesssim_u\sum_{M> N}M^{-\alpha-s_c}(1+M^{2s_c-1/2}K^{1/2}) \nonumber
	\\ &\lesssim_u N^{-\alpha-s_c}(1+N^{4s_c-1}K)^{1/2}. \label{apriori1}
	\end{align}

Next we note that for $2\leq q\leq\infty$ and $\eps<s_c-\tfrac{4s_c-1}{q}$, we have $\tfrac{(q-2)s_c-q\eps}{2}>s_c-1/2$; thus, by interpolation, \eqref{assume1.1}, and \eqref{apriori1}, we get
\begin{align*}
				\xonorm{\vert\nabla\vert^{\eps}u_{>N}}{q}{\frac{2q}{q-1}}
				&\lesssim \xonorm{\vert\nabla\vert^{-(\frac{(q-2)s_c-q\eps}{2})}u_{>N}}{2}{4}^{2/q}
				\xonorm{\nsc u_{>N}}{\infty}{2}^{1-2/q}
			\\	&\lesssim_u N^{\eps-s_c}(1+N^{4s_c-1}K)^{1/q},
\end{align*}
which settles \eqref{hc}.

We next turn to \eqref{lc}. By \eqref{lsesmall}, we know that for any $\eta>0$, there exists
$N_0=N_0(\eta)$ such that for $N\leq N_0,$ we have
\begin{equation}\nonumber
\big\||\nabla|^{s_c}u_{\leq
N}\big\|_{L_t^2L_x^4}\lesssim_u\eta\big(1+N^{2s_c-1/2}K^{1/2}).
\end{equation}
Interpolating this estimate with \eqref{assume1.1} yields \eqref{lc}.

We now consider \eqref{hc2} and \eqref{hc3}. Using Remark~\ref{rem1.1} and the fact that \newline$\inf_{t\in I}N(t)\geq 1,$ we may find $c(\eta)>0$ such that 

\begin{equation}\nonumber
\big\||\nabla|^{s_c}u_{\leq c(\eta)}\big\|_{L_t^\infty
L_x^2}\leq \eta.
\end{equation}
Combining this estimate with Bernstein, we get
\begin{align*}
N^{s_c}\|u_{\geq N}\|_{L_t^\infty L_x^2}&\lesssim
		N^{s_c}\|u_{N\leq\cdot\leq c(\eta)}\|_{L_t^\infty L_x^2}+N^{s_c}\|u_{\geq c(\eta)}\|_{L_t^\infty L_x^2}
		\\& \lesssim\big\||\nabla|^{s_c}u_{\leq c(\eta)}\big\|_{L_t^\infty L_x^2}+\tfrac{N^{s_c}}{c(\eta)^{s_c}}\big\||\nabla|^{s_c}u\big\|_{L_t^\infty L_x^2}
		\\ & \lesssim\eta + N^{s_c}.
\end{align*}
Taking $N$ sufficiently small, we get \eqref{hc2}. A similar argument gives \eqref{hc3}.

For the proof of \eqref{mass control}, we will use Theorem~\ref{lse} together with Proposition~\ref{mxs}.  We begin by noting that
	\begin{align}
	\int_I\sup_{x\in\R^4}\int_{\vert x-y\vert\leq\rho}\vert u_{>N}(t,y)\vert^2\,dy\,dt\lesssim \rho^{4}\xonorm{S_{\rho}u_{>N}}{2}{\infty}^2,
	\label{mass control one}
	\end{align}
with $S_\rho$ defined as in Proposition~\ref{mxs}. We now write $(i\partial_t+\Delta)u_{>N}=F+G$ with 
		$$F=P_{>N}(\vert u\vert^{p-2}\bar{u}u_{\leq N}^2+2\vert u\vert^{p-2}\bar{u}{u}_{\leq N}u_{>N})$$ 
and 
		$$G=P_{>N}(\vert u\vert^{p-2}\bar{u}u_{>N}^2)$$ 
and apply Proposition~\ref{mxs} to get
	\begin{align*}
	\rho^2\xonorm{S_{\rho}u_{>N}}{2}{\infty}
	&\lesssim\rho(\xonorm{u_{>N}}{\infty}{2}+\xonorm{F}{2}{4/3}+\rho^{-1}\xonorm{G}{2}{1}).
	\end{align*}

We first note that Bernstein and \eqref{assume1.1} immediately give 
	$$\xonorm{u_{>N}}{\infty}{2}\lesssim_u N^{-s_c}.$$
We will show that
	\begin{align}
	\xonorm{F}{2}{4/3}\lesssim_u N^{-s_c}(1+N^{4s_c-1}K)^{1/2} \label{mxs F}
	\end{align}
and
	\begin{align}
	\xonorm{G}{2}{1}\lesssim_u N^{-1-s_c}(1+N^{4s_c-1}K)^{1/2}.\label{mxs G}
	\end{align}
Then continuing from \eqref{mass control one} and using $\rho\geq N^{-1}$, we will have \eqref{mass control}.

We now turn to \eqref{mxs F} and \eqref{mxs G}. We estimate $F$ in two pieces. First, we use H\"older, Bernstein, the fractional product rule, the fractional chain rule, Sobolev embedding, \eqref{assume1.1} and \eqref{lc} to see
	\begin{align*}
	\xonorm{&P_{>N}(\vert u\vert^{p-2}\bar{u}u_{\leq N}^2)}{2}{4/3}
	\\ &\lesssim N^{-s_c}\xonorm{\nsc(\vert u\vert^{p-2}\bar{u}u_{\leq N}^2)}{2}{4/3}
	\\ &\lesssim N^{-s_c}\bigg\{\xonorm{\nsc(\vert u\vert^{p-2}\bar{u})}{\infty}{\frac{p}{p-1}}
					\xonorm{u_{\leq N}}{4}{\frac{8p}{4-p}}^2
		\\ &\quad\quad\quad\quad\quad+\xonorm{u}{\infty}{2p}^{p-1}\xonorm{\nsc (u_{\leq N}^2)}{2}{\frac{4p}{p+2}}\bigg\}
	\\ &\lesssim N^{-s_c}\bigg\{\xonorm{u}{\infty}{2p}^{p-2}\xonorm{\nsc u}{\infty}{2}\xonorm{\nsc u_{\leq N}}{4}{8/3}^2
	\\ &\quad\quad\quad\quad\quad+\xonorm{\nsc u}{\infty}{2}^{p-1}\xonorm{u_{\leq N}}{4}{\frac{8p}{4-p}}\xonorm{\nsc u_{\leq N}}{4}{8/3}\bigg\}
	\\ &\lesssim_u N^{-s_c}\xonorm{\nsc u_{\leq N}}{4}{8/3}^2
	\\ &\lesssim_u N^{-s_c}(1+N^{4s_c-1}K)^{1/2}.
		\end{align*}
Second, we use H\"older, Sobolev embedding, \eqref{assume1.1}, \eqref{hc}, and \eqref{lc} to see
	\begin{align*}
	\xonorm{P_{>N}&(\vert u\vert^{p-2}\bar{u}u_{\leq N}u_{>N})}{2}{4/3}
	\\ &\lesssim
		\xonorm{u}{\infty}{2p}^{p-1}\xonorm{u_{\leq N}}{4}{\frac{8p}{4-p}}\xonorm{u_{>N}}{4}{8/3}
	\\ &\lesssim \xonorm{\nsc u}{\infty}{2}^{p-1}\xonorm{\nsc u_{\leq N}}{4}{8/3}\xonorm{u_{>N}}{4}{8/3}
	\\ &\lesssim_u N^{-s_c}(1+N^{4s_c-1}K)^{1/2}.
	\end{align*}
This completes the proof of \eqref{mxs F}.

To estimate $G$, we will use H\"older, Bernstein, Schur's test, \eqref{assume1.1}, and \eqref{lse1}. We find
	\begin{align*}
	\xonorm{&P_{>N}(\vert u\vert^{p-2}\bar{u}u_{>N}^2)}{2}{1}
	\\ &\lesssim\xonorm{u}{\infty}{2p}^{p-2}\xonorm{\bar u u_{>N}^2}{2}{\frac{2p}{p+2}}
	\\ &\lesssim\xonorm{\nsc u}{\infty}{2}^{p-2}
		\bnorm{ \sum_{{\tiny \begin{array}{c} M_1\geq M_2\geq M_3 \\ M_2>N\end{array}}}\norm{u_{M_1}(t)}_{L_x^2}\norm{u_{M_2}(t)}_{L_x^4}
					\norm{u_{M_3}(t)}_{L_x^{\frac{4p}{4-p}}} }_{L_t^2}
	\\ &\lesssim_u \bnorm{\!\!\sum_{{\tiny \begin{array}{c} M_1\geq M_2\geq M_3 \\ M_2>N\end{array}}}\!\!\!\!\! M_1^{-s_c}M_2^{-s_c}M_3^{s_c-1}
				\norm{\nsc u_{M_1}}_{L_x^2}\norm{\nsc u_{M_2}}_{L_x^4}\norm{\nsc u_{M_3}}_{L_x^2} }_{L_t^2}
	\\ &\lesssim_u \bnorm{\!\sum_{M_1>M_3}\!\!(\tfrac{M_3}{M_1})^{s_c-1}\norm{\nsc u_{M_1}}_{L_x^2}\norm{\nsc u_{M_3}}_{L_x^2}}_{L_t^\infty}
			\!\!	\sum_{M>N}M^{-s_c-1}\xonorm{\nsc u_{M}}{2}{4}
	\\ &\lesssim_u\xonorm{\nsc u}{\infty}{2}^2\sum_{M>N}M^{-s_c-1}(1+M^{2s_c-1/2}K^{1/2})
	\\ &\lesssim_u N^{-s_c-1}(1+N^{4s_c-1}K)^{1/2},
	\end{align*}
where the condition $s_c<3/2$ guarantees the convergence of the sum. This completes the proof of \eqref{mxs G}, which in turn completes the proof of Corollary~\ref{a priori bounds}. 
	\end{proof}

%%% massmass %%%%%%
We will now examine each of the terms appearing in Proposition~\ref{morawetz identity}. 

We first consider \eqref{mass-mass}. This is perhaps the most important term in Proposition~\ref{morawetz identity}, as it is responsible for the left-hand side of \eqref{uppbd}. Indeed, we have the following

\begin{lemma}[Mass-mass interactions]\label{mass-mass lemma} Let $\eta>0$ and take $N_0=N_0(\eta)$ as in Corollary~\ref{a priori bounds}. Then for $N<N_0$, 

\begin{align*}
&-\int_I\iint_{\R^4\times\R^4}\vert u_{>N}(t,y)\vert^2a_{jjkk}(x-y)\vert u_{>N}(t,x)\vert^2\,dx\,dy\,dt
\\ &\quad\quad\quad-3\int_I\iint_{\R^4\times\R^4}\frac{ \vert u_{>N}(t,y)\vert^2 u_{>N}(t,x)\vert^2}{\vert x-y\vert^3}\,dx\,dy\,dt
\\ &\lesssim_u{\eta}(N^{1-4s_c}+K).
\end{align*}

\end{lemma}

\begin{proof}
In four dimensions, one has $-\Delta\Delta(\vert x\vert)=3\vert x\vert^{-3}$. Thus, using \eqref{properties of a}, we find that to prove the lemma it will suffice to show
	\begin{align}
	\int_I\iint_{\vert x-y\vert>N^{-1}}\frac{\vert u_{>N}(t,y)\vert^2\vert u_{>N}(t,x)\vert^2}{\vert x-y\vert^3}\,dx\,dy\,dt\lesssim_u\eta(N^{1-4s_c}+K).
	\label{masmas}
	\end{align}

The proof of \eqref{masmas} will rely on \eqref{mass control} and \eqref{hc2}. We estimate
	\begin{align*}
	\int_I&\iint_{\vert x-y\vert>N^{-1}}\frac{\vert u_{>N}(t,y)\vert^2\vert u_{>N}(t,x)\vert^2}{\vert x-y\vert^3}\,dx\,dy\,dt
	\\ &\lesssim \sum_{j=0}^\infty \int_I\iint_{2^jN^{-1}<\vert x-y\vert\leq 2^{j+1}N^{-1}}\frac{\vert u_{>N}(t,y)\vert^2\vert u_{>N}(t,x)\vert^2}{\vert x-y\vert^3}\,dx\,dy\,dt
	\\ &\lesssim\sum_{j=0}^{\infty}2^{-3j}N^{3}\xonorm{u_{>N}}{\infty}{2}^2\int_I\sup_x\int_{\vert x-y\vert\leq 2^{j+1}N^{-1}}\vert u_{>N}(t,y)\vert^2\,dy\,dt
	\\ &\lesssim_u\eta^2N^{3-2s_c}\sum_{j=0}^\infty2^{-3j}(2^{j+1}N^{-1})^2 N^{-2s_c}(1+N^{4s_c-1}K)
	\\ &\lesssim_u \eta (N^{1-4s_c}+K),
	\end{align*}
as needed. \end{proof} 

We next turn to the contribution of \eqref{scary}. In the standard interaction Morawetz inequality (that is, when one takes $a(x)=\vert x\vert$), this term can be shown to be non-negative, and hence it can be thrown away (see \cite[Proposition~2.5]{ckstt}). To see this, one relies on the convexity of $\vert x\vert$, which one loses by truncating the weight $a(x)$. In particular, we can no longer ignore this term. 

To deal with this term requires the estimation of $\|\nabla u_{>N}\|_{L_t^\infty L_x^2}$. Hence, this is a potentially dangerous term in the energy-critical (or energy-subcritical) case. In our setting, however, we can still exhibit smallness in $\|\nabla u_{>N}\|_{L_t^\infty L_x^2}$ (cf. \eqref{hc3}), and so this term is not particularly dangerous. In fact, this is one of the reasons why our weight $a$ can be so much simpler than the weight $a$ used in \cite{KV3D} (see also the discussion preceding Lemma~\ref{potential energy}).

We estimate the contribution of \eqref{scary} as follows:

\begin{lemma}[Contribution of \eqref{scary}]\label{nssc} Let $\eta>0$ and take $N_0=N_0(\eta)$ as in Corollary~\ref{a priori bounds}. Defining
	\begin{align*}\Phi_{jk}(t;x,y)&:=\vert u_{>N}(t,y)\vert^2\Re(\partial_j \bar{u}_{>N}\partial_k u_{>N})(t,x)
	\\ &\quad\quad-\Im(\bar{u}_{>N}\partial_j u_{>N})(t,y)\Im(\bar{u}_{>N}\partial_k u_{>N})(t,x),
	\end{align*}
we have
	\begin{align}
	\int_I\iint_{\vert x-y\vert\leq N^{-1}} 4a_{jk}(x-y)\Phi_{jk}(t;x,y)\,dx\,dy\,dt&\geq 0, \label{positive}
	\\ \bigg\vert\int_I\iint_{\vert x-y\vert >N^{-1}} 4a_{jk}(x-y)\Phi_{jk}(t;x,y)\,dx\,dy\,dt\bigg\vert&\lesssim\eta(N^{1-4s_c}+K).\label{nscr}
	\end{align}	
\end{lemma}

\begin{proof} We follow the argument presented in \cite[Lemma~6.6]{KV3D}. 

We first prove \eqref{positive} by showing that in this term the integrand is positive.

Noting that $a_{jk}(x-y)=a_{jk}(y-x)$, we see that we may replace $\Phi$ by the Hermitian matrix $$M_{jk}(t;x,y):=\tfrac12\Phi_{jk}(t;x,y)+\tfrac12\Phi_{jk}(t;y,x).$$ Now, for fixed $t\in I$ and $x,y\in\R^4$, we claim that $M_{jk}(t;x,y)$ is a positive semi-definite quadratic form on $\R^4.$ Indeed, this follows from the following simple estimate: for $e\in\R^4$ and any function $\varphi$, we have
	\begin{align*}\vert e_ke_j\Im(\bar{\varphi}\varphi_j)(y)\Im(\bar{\varphi}\varphi_k)(x)\vert
	&\leq\vert\varphi(y)\vert\,\vert e\cdot\nabla\varphi(y)\vert\,\vert\varphi(x)\vert\,\vert e\cdot\nabla\varphi(x)\vert
	\\&\leq\tfrac12\vert\varphi(x)\vert^2\vert e\cdot\nabla\varphi(y)\vert^2+\tfrac12\vert\varphi(y)\vert^2\vert e\cdot\nabla\varphi(x)\vert^2.
	\end{align*}

Next, we note that as $a$ is a real symmetric matrix, it has real eigenvectors. Moreover, recalling that $a(x)$ is convex for $\vert x\vert\leq N^{-1}$ (it equals $\vert x\vert$ in this region), we see that $a_{jk}$ is positive semi-definite on the region of integration in \eqref{positive}. Thus we have that the integrand is positive in this region, as we needed to show.

We now turn to \eqref{nscr}. Using the fact that $a$ is a spherically symmetric function on $\R^4$, we get that the eigenvalues of $\nabla^2a$ are $a_{rr}$ and $\tfrac1r a_r$, where $\tfrac1r a_r$ has multiplicity $3$. Now, we have chosen $a$ so that $a_r\geq 0$ and $\vert a_{rr}\vert\lesssim r^{-1}$, and hence \eqref{nscr} reduces to showing
	\begin{align}
	\int_I\int_{N^{-1}< \vert x-y\vert\leq 2N^{-1}}\frac{\vert\nabla u_{>N}(t,x)\vert^2\vert u_{>N}(t,y)\vert^2}{\vert x-y\vert}\,dx\,dy\,dt\lesssim_u \eta(N^{1-4s_c}+K).
	\label{not so scary}
	\end{align}

In fact, using \eqref{hc3} and \eqref{mass control} we get
	\begin{align*}
	\LHS\eqref{not so scary}&
	\lesssim N\xonorm{\nabla u_{>N}}{\infty}{2}^2\int_I\sup_x\int_{\vert x-y\vert\leq 2N^{-1}}\vert u_{>N}(t,y)\vert^2\,dy\,dt
	\\ &\lesssim \eta^2 N^{3-2s_c}N^{-2-2s_c}(1+N^{4s_c-1}K)
	\\ &\lesssim\eta (N^{1-4s_c}+K).	
	\end{align*}
This concludes the proof of Lemma~\ref{nssc}.
\end{proof}

It remains to estimate the contribution of the mass and momentum bracket terms in Proposition~\ref{morawetz identity} to \eqref{uppbd}. The key to estimating these terms successfully will be to make use of the term appearing on the \emph{left}-hand side of \eqref{uppbd} (that is, we will ultimately prove \eqref{uppbd} via a bootstrap argument). In particular, we will make frequent use of the terms appearing in the following

\begin{lemma}[Bootstrap terms]\label{bootstrap}
\begin{align}
				\xnorms{\vert\nabla\vert^{-1/4}u_{>N}}{4}{I}^4\lesssim\xnorms{\vert\nabla\vert^{-1/2}\vert u_{>N}\vert^2}{2}{I}^2
				&\sim\LHS\eqref{limemd}				
\nonumber
\end{align}
\end{lemma}
\begin{proof}
This lemma follows from the general inequalities		
	\begin{equation}\label{bst1}
	\norm{\vert\nabla\vert^{-1/2}\vert f\vert^2}_{L_{x}^2(\R^4)}^2
	\sim\iint_{\R^4\times\R^4}\frac{\vert f(t,y)\vert^2\vert f(t,x)\vert^2}{\vert x-y\vert^3}\,dx\,dy
	\end{equation}
and
	\begin{equation}\label{bst2}
	\norm{\vert\nabla\vert^{-1/4}f}_{L_x^4(\R^4)}^4\lesssim\norm{\vert\nabla\vert^{-1/2}\vert f\vert^2}_{L_x^2(\R^4)}^2.
	\end{equation}
 
The estimate \eqref{bst1} follows from the fact that in four dimensions, convolution with $\vert x\vert^{-3}$ is essentially the fractional differentiation operator $\vert\nabla\vert^{-1}$.

The original proof of \eqref{bst2} may be found in \cite[Section 5.2]{Visan2007}. We sketch the argument briefly as follows: writing $S$ for the Littlewood--Paley square function, one can establish the following pointwise inequality (for positive Schwartz functions, say):
		$$\vert S(\vert\nabla\vert^{-1/4}f)(x)\vert^2\lesssim(\vert\nabla\vert^{-1/2}\vert f\vert^2)(x).$$
The estimate \eqref{bst2} then follows from an application of the square function estimate.
\end{proof}

Before turning to the mass and momentum bracket terms, we record the following useful pointwise inequalities:
		\begin{align}
		\label{pw1}
		\vert F(u)-F(u_{> N})-F(u_{\leq N})\vert&\lesssim \vert u_{> N}u_{\leq N}u^{p-1}\vert,
	\\	\big\vert \vert u\vert^{p+2}-\vert u_{>N}\vert^{p+2}-\vert u_{\leq N}\vert^{p+2}\big\vert&\lesssim \vert u_{>N}u_{\leq N}u^p\vert. \label{pw2}
		\end{align}

We now turn to the estimation of the mass bracket terms. 

\begin{lemma}[Mass bracket terms]\label{mbterms} Let $\eta>0$ and take $N_0=N_0(\eta)$ as in Corollary~\ref{a priori bounds}. Then for $N\leq N_0$, we have
\begin{align*}
&\int_I\iint_{\R^4\times\R^4}\{P_{>N}(F(u)),u_{>N}\}_M(t,y)\nabla a(x-y)\cdot\Im(\bar{u}_{>N}\nabla u_{>N})(t,x)\,dx\,dy\,dt
\\ &\lesssim_u \eta(N^{1-4s_c}+K)+\eta\xonorms{\vert\nabla\vert^{-1/4}u_{>N}}{4}^4
	+\eta\xonorms{\vert\nabla\vert^{-1/2}\vert u_{>N}\vert^2}{2}^2,
	\end{align*}
where all spacetime norms are taken over $I\times\R^4$. 

\end{lemma}
\begin{proof}[Proof of Lemma~\ref{mbterms}]
To begin, we note that $\{F(u_{>N}),u_{>N}\}_M=0$. Thus, we may write
	\begin{align*}
	\{P_{>N}(F(u)),u_{>N}\}_{M}=&\,\{P_{>N}(F(u))-F(u_{>N}),u_{>N}\}_M
						\\=&\,\{P_{>N}(F(u)-F(u_{>N})-F(u_{\leq N})),u_{>N}\}_M
						\\&	-\{P_{\leq N}(F(u_{>N})),u_{>N}\}_{M}
							+\{P_{>N}(F(u_{\leq N})),u_{>N}\}_M
					 	\\=:&\,I+II+III.
	\end{align*}
To estimate the contribution of $I$, we consider two cases. 

In the case $2<p<8/3$, we may choose $q$ such that $\tfrac{7p-8}{2(p-1)}<q<\tfrac{2p}{p-1}$. The lower bound on $q$ will be necessary to apply \eqref{hc} below with $\eps=0$, while the upper bound $q$ guarantees that the exponent $\tfrac{2pq}{2p+q-pq}$ (which will appear below) is finite. 

In this case, using H\"older, Sobolev embedding, \eqref{assume1.1}, \eqref{properties of a}, \eqref{hc}, \eqref{lc}, \eqref{hc3}, and \eqref{pw1}, we see that the contribution of $I$ is controlled by
	\begin{align*}
	\xonorm{&u_{>N}^2u_{\leq N}u^{p-1}}{\frac{q}{q-1}}{1}\cdot\sup_{y\in\R^4}\xonorm{u_{>N}\nabla u_{>N}\nabla a(x-y)}{q}{1}
	\\ &\lesssim\xonorm{u_{>N}}{q}{\frac{2q}{q-1}}^2\xonorm{u_{\leq N}}{\frac{q}{q-3}}{\frac{2pq}{2p+q-pq}}
			\xonorm{u}{\infty}{2p}^{p-1}
			\\ &\quad\quad\quad\times \xonorm{u_{>N}}{q}{\frac{2q}{q-1}}
			\xonorm{\nabla u_{>N}}{\infty}{2}
					\sup_{y\in\R^4}\norm{\nabla a(x-y)}_{L_x^{2q}}
	\\ &\lesssim_u N^{-2s_c}(1+N^{4s_c-1}K)^{2/q}N^{2/q}\xonorm{u_{\leq N}}{\frac{q}{q-3}}{\frac{2pq}{3p+q-pq}}
			\\ &\quad\quad\quad\times \eta N^{1-2/q-2s_c}(1+N^{4s_c-1}K)^{1/q}
	\\ &\lesssim_u \eta N^{1-4s_c}(1+N^{4s_c-1}K)^{3/q}\xonorm{\nsc u_{\leq N}}{\frac{q}{q-3}}{\frac{2q}{3}}
	\\ &\lesssim_u \eta(N^{1-4s_c}+K),
	\end{align*}
which is acceptable. 

In the case $8/3\leq p<4$ (i.e. $s_c\geq 5/4$), we need to work a bit harder. In light of \eqref{pw1}, we are faced with estimating a term of the form
	\begin{equation}\label{mass I}
	\iiint (u_{>N}^2u_{\leq N}u^{p-1})(t,y)\nabla a(x-y)\cdot ({u}_{>N}\nabla u_{>N})(t,x)\,dx\,dy\,dt.
	\end{equation}
We first use Plancherel, H\"older, Lemma~\ref{pde}, Bernstein, interpolation, \eqref{assume1.1}, and \eqref{properties of a} to estimate
	\begin{align*}
	\int_{\R^4}&\nabla a(x-y)\cdot(u_{>N}\nabla u_{>N})(t,x)\,dx
	\\ &=\int_{\R^4}\vert\nabla\vert^{-1/4}(\nabla a(x-y)u_{>N}(t,x))
					\cdot\vert\nabla\vert^{1/4}\nabla u_{>N}(t,x)\,dx
	\\ &\lesssim \norm{\vert\nabla\vert^{-1/4}\big(\nabla a(x-y) u_{>N}(t)\big)}_{L_x^2}\norm{\vert\nabla\vert^{5/4} u_{>N}(t)}_{L_x^2}
	\\ &\lesssim \norm{\vert\nabla\vert^{-1/4}u_{>N}(t)}_{L_x^4}\norm{\vert\nabla\vert^{1/4}\nabla a(x-y)}_{L_x^{16/5}}N^{5/4-s_c}\xonorm{\nsc u_{>N}}{\infty}{2}
	\\ &\lesssim_u N^{5/4-s_c}\norm{\vert\nabla\vert^{-1/4}u_{>N}(t)}_{L_x^4}
		\norm{\nabla a(x-y)}_{L_x^{36/11}}^{3/4}\norm{\Delta a(x-y)}_{L_x^3}^{1/4}
	\\ &\lesssim_u N^{1/4-s_c}\norm{\vert\nabla\vert^{-1/4}u_{>N}(t)}_{L_x^4}
	\end{align*}
uniformly for $y\in\R^4$. We can now use this estimate, together with H\"older, Bernstein, interpolation, Sobolev embedding, Young's inequality, \eqref{assume1.1}, \eqref{hc}, and \eqref{lc} to get
	\begin{align*}
	\eqref{mass I}&
	\lesssim_u N^{1/4-s_c}\xonorms{\vert\nabla\vert^{-1/4}u_{>N}}{4}\xonorm{u_{\leq N}}{4}{\infty}
					\xonorm{u_{>N}}{4}{\frac{4p}{p+1}}^2\xonorm{u}{\infty}{2p}^{p-1}
	\\ &\lesssim_u N^{7/4-2s_c}\xonorms{\vert\nabla\vert^{-1/4}u_{>N}}{4}
				\xonorm{u_{\leq N}}{4}{\frac{8p}{4-p}}
		\\ &\quad\quad\quad\times\bigg[\xonorms{\vert\nabla\vert^{-1/4}u_{>N}}{4}^{1-2/p}
					\xonorm{\vert\nabla\vert^{(p-2)/8}u_{>N}}{4}{8/3}^{2/p}\bigg]^2
	\\ &\lesssim_u N^{7/4-2s_c}\xonorms{\vert\nabla\vert^{-1/4}u_{>N}}{4}^{3-4/p}
				\xonorm{\nsc u_{\leq N}}{4}{8/3}\xonorm{\vert\nabla\vert^{(p-2)/8}u_{>N}}{4}{8/3}^{4/p}
	\\ &\lesssim_u \eta N^{7/4-2s_c}\big( N^{s_c-(p-2)/8}\big)^{-4/p}(1+N^{4s_c-1}K)^{1/4+1/p}\xonorms{\vert\nabla\vert^{-1/4} u_{>N}}{4}^{3-4/p}
	\\ &\lesssim_u \eta \big(N^{1-4s_c}(1+N^{4s_c-1}K)\big)^\frac{4+p}{4p}\xonorms{\vert\nabla\vert^{-1/4}u_{>N}}{4}^\frac{3p-4}{p}
	\\ &\lesssim_u \eta(N^{1-4s_c}+K)+\eta\xonorms{\vert\nabla\vert^{-1/4}u_{>N}}{4}^4,
	\end{align*}
which is acceptable.

To estimate the contribution of $II$, we write $u_{>N}=\nabla\cdot\nabla\Delta^{-1}u_{>N}$ and integrate by parts. We find that we are left to estimate
	\begin{align}
	\iiint\{\nabla P_{\leq N}(F(u_{>N})),\vert\nabla\vert^{-1}u_{>N}\}_M(t,y)\nabla a(x-y)\Im(\bar{u}_{>N}\nabla u_{>N})(t,x)\,dx\,dy\,dt \label{IIa}
	\\ +
	\iiint\{P_{\leq N}(F(u_{>N})),\vert\nabla\vert^{-1}u_{>N}\}_M(t,y)\Delta a(x-y)\Im(\bar{u}_{>N}\nabla u_{>N})(t,x)\,dx\,dy\,dt. \label{IIb}
	\end{align}
	
To estimate \eqref{IIa}, we use H\"older, Bernstein, Lemma~\ref{pde}, the fractional chain rule, Sobolev embedding, Young's inequality, \eqref{assume1.1}, \eqref{properties of a}, \eqref{hc}, \eqref{hc2}, and \eqref{hc3}. We find
	\begin{align*}
	\eqref{IIa}&\lesssim\xonorms{\vert\nabla\vert^{-1}u_{>N}\,\nabla P_{\leq N}(F(u_{>N}))}{1}
			\xonorm{u_{>N}}{\infty}{2}\xonorm{\nabla u_{>N}}{\infty}{2}
	\\ &\lesssim_u \eta^2 N^{1-2s_c}\xonorm{\vert\nabla\vert^{-1}u_{>N}}{2}{4}\xonorm{\nabla P_{\leq N}(F(u_{>N}))}{2}{4/3}
	\\ &\lesssim_u \eta^2N^{-3s_c}(1+N^{4s_c-1}K)^{1/2}N^{s_c+1/2}\xonorm{\vert\nabla\vert^{-1/2}(F(u_{>N}))}{2}{4/3}
	\\ &\lesssim_u \eta^2 (N^{1-4s_c}+K)^{1/2}\xonorms{\vert\nabla\vert^{-1/2}\vert u_{>N}\vert^2}{2} 
		\xonorm{\vert\nabla\vert^{1/2}(\vert u_{>N}\vert^{p-2}u_{>N})}{2}{\frac{2p}{2p-1}}
	\\ &\lesssim_u \eta^2 (N^{1-4s_c}+K)^{1/2}	\xonorms{\vert\nabla\vert^{-1/2}\vert u_{>N}\vert^2}{2}
		\xonorm{u}{\infty}{2p}^{p-2}\xonorm{\vert\nabla\vert^{1/2}u}{\infty}{\frac{8p}{p+4}}
	\\ &\lesssim_u \eta^2 (N^{1-4s_c}+K)^{1/2}\xonorms{\vert\nabla\vert^{-1/2}\vert u_{>N}\vert^2}{2}\xonorm{\nsc u}{\infty}{2}^{p-1}	
	\\ &\lesssim_u \eta^2(N^{1-4s_c}+K)+\eta^2\xonorms{\vert\nabla\vert^{-1/2}\vert u_{>N}\vert^2}{2}^2,
	\end{align*}
which is acceptable.

To estimate \eqref{IIb}, we will use H\"older, Hardy--Littlewood--Sobolev, Bernstein, Lemma~\ref{pde}, the fractional chain rule, Sobolev embedding, Young's inequality, \eqref{assume1.1}, \eqref{properties of a}, \eqref{hc}, and \eqref{hc2}. We find
	\begin{align*}
	\eqref{IIb}&\lesssim \xonorm{\vert\nabla\vert^{-1}u_{>N}}{2}{4} \xonorm{P_{\leq N}(F(u_{>N})}{4}{8/5}
					\xonorm{\vert x\vert^{-1}\ast u_{>N}\nabla u_{>N}}{4}{8} 
			\\ &\lesssim_u N^{-1-s_c}(1+N^{4s_c-1}K)^{1/2}N^{3/4}\xonorm{\vert\nabla\vert^{-1/4}(F(u_{>N}))}{4}{4/3}
					\\ &\quad\quad\quad\times	\xonorm{u_{>N}\nabla u_{>N}}{4}{8/7}
			\\ &\lesssim_u N^{-1/4-s_c}(1+N^{4s_c-1}K)^{1/2}\xonorms{\vert\nabla\vert^{-1/4}u_{>N}}{4}
						\xonorm{\vert\nabla\vert^{1/4}\vert u_{>N}\vert^p}{\infty}{16/9}
					\\ &\quad\quad\quad\times\xonorm{u_{>N}}{4}{8/3}\xonorm{\nabla u_{>N}}{\infty}{2}
			\\ &\lesssim_u\! \eta N^{3/4-3s_c}(1+N^{4s_c-1}K)^{3/4}\xonorms{\vert\nabla\vert^{-1/4}u_{>N}}{4}\!
						\xonorm{u}{\infty}{2p}^{p-1}\xonorm{\vert\nabla\vert^{1/4}u}{\infty}{\frac{16p}{p+8}}
			\\ &\lesssim_u \eta(N^{4s_c-1}+K)^{3/4}\xonorms{\vert\nabla\vert^{-1/4}u_{>N}}{4}\xonorm{\nsc u}{\infty}{2}^p
			\\ &\lesssim_u \eta(N^{4s_c-1}+K)+\eta\xonorms{\vert\nabla\vert^{-1/4}u_{>N}}{4},
	\end{align*}
which is acceptable.

We estimate the contribution of $III$ as follows. Using H\"older, Bernstein, the fractional chain rule, Sobolev embedding, \eqref{properties of a}, \eqref{lc}, \eqref{hc2}, and \eqref{hc3}, we see that the contribution of $III$ is controlled by
	\begin{align*}
	\xonorms{&P_{>N}(F(u_{\leq N}))u_{>N}}{1}\xonorm{u_{>N}}{\infty}{2}\xonorm{\nabla u_{>N}}{\infty}{2}
	\\ &\lesssim \xonorm{u_{>N}}{4}{8/3}\xonorm{P_{>N}(F(u_{\leq N})}{4/3}{8/5}\xonorm{u_{>N}}{\infty}{2}\xonorm{\nabla u_{>N}}{\infty}{2}
	\\ &\lesssim_u\eta^2N^{1-3s_c}(1+N^{4s_c-1}K)^{1/4}N^{-s_c}\xonorm{\nsc F(u_{\leq N})}{4/3}{8/5}
	\\ &\lesssim_u \eta^2N^{1-4s_c}(1+N^{4s_c-1}K)^{1/4}\xonorm{u_{\leq N}}{4p}{{8p}/{3}}^p\xonorm{\nsc u_{\leq N}}{2}{4}
	\\ &\lesssim_u\eta^3N^{1-4s_c}(1+N^{4s_c-1}K)^{3/4}\xonorm{\nsc u_{\leq N}}{4p}{\frac{8p}{4p-1}}^p
	\\ &\lesssim_u \eta^{p+3}N^{1-4s_c}(1+N^{4s_c-1}K)
	\\ &\lesssim_u\eta(N^{1-4s_c}+K),
	\end{align*}
which is acceptable.

\end{proof}

We now turn to the momentum bracket terms. The primary contribution of this term will be a quantity involving the potential energy of the solution (which we estimate separately in Lemma~\ref{potential energy} below); the remaining contribution will consist of error terms that can be controlled in an acceptable way by making use of Corollary~\ref{a priori bounds} and the terms appearing in Lemma~\ref{bootstrap}.  

\begin{lemma}[Momentum bracket terms]\label{Pb terms} Let $\eta>0$ and let $N_0=N_0(\eta)$ be as in Corollary~\ref{a priori bounds}. For $N\leq \min\{N_0,1\}$, we have
	\begin{align}
	\int_I&\iint_{\R^4\times\R^4}\vert u_{>N}(t,y)\vert^2\ 2\nabla a(x-y) \nonumber
	\cdot\{P_{>N}(F(u)),u_{>N}\}_{P}(t,x)\,dx
	\\-&\tfrac{2p}{p+2}\int_I\iint_{\R^4\times\R^4}\vert u_{>N}(t,y)\vert^2\Delta a(x-y)
		\vert u_{>N}(t,x)\vert^{p+2}\,dx\,dy\,dt\nonumber
	\\ &\lesssim_u \eta(N^{4s_c-1}+K)+\eta\xonorms{\vert\nabla\vert^{-1/4}u_{>N}}{4}^4+
	\eta\xonorms{\vert\nabla\vert^{-1/2}\vert u_{>N}\vert^2}{2}^2,	\label{Pb}
	\end{align}
where all spacetime norms taken over $I\times\R^4$.
\end{lemma}
\begin{proof}[Proof of Lemma~\ref{Pb terms}] To begin, we note that $\{F(u),u\}_{P}=-\tfrac{p}{p+2}\nabla\vert u\vert^{p+2}.$ Thus, we may write
	\begin{align}
	\{P_{>N}&(F(u)),u_{>N}\}_P& \nonumber
	\\ &=\{F(u),u\}_P-\{F(u_{\leq N}),u_{\leq N}\}_P-\{F(u)-F(u_{\leq N}),u_{\leq N}\}_P	\nonumber
	\\ &\quad-\{P_{\leq N}(F(u)),u_{>N}\}_P	\nonumber
	\\ &=-\tfrac{p}{p+2}\nabla(\vert u\vert^{p+2}-\vert u_{\leq N}\vert^{p+2})	\nonumber
		-\{F(u)-F(u_{\leq N}),u_{\leq N}\}_{P}
	\\ &\quad-\{P_{\leq N}(F(u)),u_{>N}\}_P	\nonumber
	\\ &=:I+II+III.	\label{momb decomp}
	\end{align}
	
We first estimate the contribution of term $I$ in \eqref{momb decomp}. Integrating by parts, we see that term $I$ contributes to the left-hand side of \eqref{Pb} 
	\begin{align}
	&\tfrac{2p}{p+2}\iiint\vert u_{>N}(t,y)\vert^2\Delta a(x-y)
		\vert u_{>N}(t,x)\vert^{p+2}\,dx\,dy\,dt\nonumber
	\\ &+\tfrac{2p}{p+2}\iiint\vert u_{>N}(t,y)\vert^2\Delta a(x-y)(\vert u\vert^{p+2}-\vert u_{> N}\vert^{p+2}-\vert u_{\leq N}\vert^{p+2})(t,x)\,dx\,dy\,dt. \label{momb I}
	\end{align}
The first term above appears on the left-hand side of \eqref{Pb}; thus, to control the contribution of $I$ we are left to bound the second term above. To estimate this term, we will need to make use of the additional decay established in Proposition~\ref{adecay1}. Specifically, we will use H\"older, interpolation, Sobolev embedding, Young's inequality, \eqref{assume1.1},  \eqref{breaking scaling}, \eqref{properties of a}, \eqref{hc}, \eqref{lc}, and \eqref{pw2}. Letting $\theta=\tfrac{p^2-p+2}{2p(p-1)}\in(0,1)$ and noting $N^{1-s_c}\geq1$, we find
	\begin{align*}
	&\iiint\vert u_{>N}(t,y)\vert^2\Delta a(x-y)(\vert u\vert^{p+2}-\vert u_{> N}\vert^{p+2}-\vert u_{\leq N}\vert^{p+2})(t,x)\,dx\,dy\,dt
	\\ &=\iiint \vert\nabla\vert^{-1/2}\vert u_{>N}(t,y)\vert^2\vert\nabla\vert^{1/2}\Delta a(x-y)
		\\ &\quad\quad\quad\times(\vert u\vert^{p+2}-\vert u_{>N}\vert^{p+2}-\vert u_{\leq N}\vert^{p+2})(t,x)\,dx\,dy\,dt
	\\ &\lesssim\sup_{x\in\R^4}\norm{(\vert\nabla\vert^{-1/2}\vert u_{>N}\vert^2)(\vert\nabla\vert^{1/2}\Delta a(x-y))}_{L_t^2L_y^1}
			\xonorm{\text{\O}(u_{\leq N}u_{>N}u^p)}{2}{1}
	\\ &\lesssim\xonorms{\vert\nabla\vert^{-1/2}\vert u_{>N}\vert^2}{2}\sup_{x\in\R^4}\norm{\vert\nabla\vert^{1/2}\Delta a(x-y)}_{L_y^2}
		\\ &\quad\quad\quad\times\xonorm{u_{\leq N}}{4}{\frac{8p}{4-p}}\xonorm{u_{>N}}{4}{8/3}\xonorm{u}{\infty}{2p}^{p\theta}\xonorm{u}{\infty}{p+1}^{p(1-\theta)}
	\\ &\lesssim_u  \xonorms{\vert\nabla\vert^{-1/2}\vert u_{>N}\vert^2}{2}\sup_{x\in\R^4}
	\norm{\Delta a(x-y)}_{L_y^{9/4}}^{1/2}\norm{\nabla\Delta a(x-y)}_{L_y^{9/5}}^{1/2}
	\\ &\quad\quad\quad\times \xonorm{\nsc u_{\leq N}}{4}{8/3} N^{-s_c}(1+N^{4s_c-1}K)^{1/4}
	\\ &\lesssim_u \eta \xonorms{\vert\nabla\vert^{-1/2}\vert u_{>N}\vert^2}{2}N^{-1/2-s_c}(1+N^{4s_c-1}K)^{1/2}
	\\ &\lesssim_u \eta\xonorms{\vert\nabla\vert^{-1/2}\vert u_{>N}\vert^2}{2}N^{1/2-2s_c}(1+N^{4s_c-1}K)^{1/2}
	\\ &\lesssim_u \eta\xonorms{\vert\nabla\vert^{-1/2}\vert u_{>N}\vert^2}{2}^2
		+\eta(N^{1-4s_c}+K),
	\end{align*}
which is acceptable.
 
Next, we estimate the contribution of term $II$ in \eqref{momb decomp}. We begin by writing
	\begin{align*} 
	\{F(u)-F(u_{\leq N}),u_{\leq N}\}_{P}
	&=\nabla\text{\O}\big( (F(u)-F(u_{\leq N}))u_{\leq N}\big)	
	\\ &\quad\quad+\text{\O}\big( (F(u)-F(u_{\leq N}))\nabla u_{\leq N}\big)
	\end{align*}
When the gradient falls on $\text{\O}\big ( (F(u)-F(u_{\leq N}))u_{\leq N}\big)$, we integrate by parts. We also write
			$$F(u)-F(u_{\leq N})=u_{>N}\text{\O}(u^p).$$		
In this way, we find that to estimate the contribution of term $II$ it will suffice to control the terms
	\begin{align}
	&\iiint \vert u_{>N}(t,y)\vert^2\big(u_{>N}\text{\O}(u^p)\big)(t,x)\nabla a(x-y)\cdot\nabla u_{\leq N}(t,x)\,dx\,dy\,dt \label{iiiiii}
	\\ &+\iiint \vert u_{>N}(t,y)\vert^{2}\big(u_{>N}\text{\O}(u^p)\big)(t,x)\Delta a(x-y)u_{\leq N}(t,x)\,dx\,dy\,dy.		\label{ivvvi}
	\end{align}
We will treat these two terms separately.

To estimate \eqref{iiiiii}, we first write $u_{>N}=\nabla\cdot\nabla\Delta^{-1}u_{>N}$ and integrate by parts. In this way, we see that to control \eqref{iiiiii}, it will suffice to control
	\begin{align}
	&\iiint \vert u_{>N}(t,y)\vert^2\vert\nabla\vert^{-1}u_{>N}(t,x)\nabla\text{\O}(u^p)(t,x)\nabla a(x-y)\nabla u_{\leq N}(t,x)\,dx\,dy\,dt
	\label{momb i}
	\\ &+\iiint\vert u_{>N}(t,y)\vert^2\vert\nabla\vert^{-1}u_{>N}(t,x)\text{\O}(u^p)(t,x)\Delta a(x-y)\nabla u_{\leq N}(t,x)\,dx\,dy\,dt
	\label{momb ii}
	\\ &+\iiint\vert u_{>N}(t,y)\vert^2\vert\nabla\vert^{-1}u_{>N}(t,x)\text{\O}(u^p)(t,x)\nabla a(x-y)\Delta u_{\leq N}(t,x)\,dx\,dy\,dt.
	\label{momb iii}
	\end{align}

We first turn to \eqref{momb i}. Using H\"older, Bernstein, Sobolev embedding, Young's inequality \eqref{assume1.1}, \eqref{properties of a}, \eqref{hc}, \eqref{lc}, and \eqref{hc2}, we find
	\begin{align*}
	&\eqref{momb i}
	\\ &\lesssim\sup_{x\in\R^4}\norm{\vert u_{>N}\vert^2\,\nabla a(x-y)}_{L_t^4L_y^1}\xonorm{\vert\nabla\vert^{-1}u_{>N}\nabla\text{\O}(u^p)\nabla u_{\leq N}}{4/3}{1}
	\\ &\lesssim\xonorm{u_{>N}}{4}{8/3}\xonorm{u_{>N}}{\infty}{2}\sup_{x\in\R^4}\norm{\nabla a(x-y)}_{L_y^8}
	 \\&\quad\quad\quad\times	\xonorms{\vert\nabla\vert^{-1}u_{>N}}{4}\xonorm{\nabla(u^p)}{\infty}{4/3}\xonorm{\nabla u_{\leq N}}{2}{\infty}
	 \\&\lesssim_u \eta N^{-2s_c}(1+N^{4s_c-1}K)^{1/4}N^{-1/2}
	  \\ &\quad\quad\quad\times 	N^{-3/4+2/p}\xonorms{\vert\nabla\vert^{-1/4}u_{>N}}{4}\xonorm{u}{\infty}{2p}^{p-1}\xonorm{\nabla u}{\infty}{\frac{4p}{p+2}}\xonorm{\nabla u_{\leq N}}{2}{2p}
	  \\ &\lesssim_u \eta N^{3/4-3s_c}(1+N^{4s_c-1}K)^{1/4}\xonorms{\vert\nabla\vert^{-1/4}u_{>N}}{4}\xonorm{\nsc u}{\infty}{2}^{p}
	  	\xonorm{\nsc u_{\leq N}}{2}{4}
	\\ &\lesssim_u\eta (N^{1-4s_c}+K)^{3/4}\xonorms{\vert\nabla\vert^{-1/4}u_{>N}}{4}
	\\ &\lesssim_u\eta(N^{1-4s_c}+K)+\eta\xonorms{\vert\nabla\vert^{-1/4}u_{>N}}{4}^4,
	\end{align*}
which is acceptable. 

We next turn to \eqref{momb ii}. Using H\"older, Hardy--Littlewood--Sobolev, Bernstein, Sobolev embedding, Young's inequality, \eqref{assume1.1}, \eqref{properties of a}, \eqref{hc}, \eqref{lc}, and \eqref{hc2}, we estimate
	\begin{align*} 
	\eqref{momb ii}
	 &\lesssim\xonorm{\vert x\vert^{-1}\ast\vert u_{>N}\vert^2}{4}{8}\xonorm{\vert\nabla\vert^{-1}u_{>N} \text{\O}(u^p)\nabla u_{\leq N}}{4/3}{8/7}
	\\ &\lesssim \xonorm{\vert u_{>N}\vert^2}{4}{8/7}\xonorms{\vert\nabla\vert^{-1}u_{>N}}{4}
		\xonorm{u}{\infty}{2p}^p\xonorm{\nabla u_{\leq N}}{2}{8}
	\\ &\lesssim \xonorm{u_{>N}}{\infty}{2}\xonorm{u_{>N}}{4}{8/3}
		N^{-3/4}\xonorms{\vert\nabla\vert^{-1/4}u_{>N}}{4}
		\\ &\quad\quad\quad\times \xonorm{\nsc u}{\infty}{2}^p N^{3/2-s_c}\xonorm{\nabla u_{\leq N}}{2}{2p}
	\\ &\lesssim_u \eta N^{3/4-3s_c}(1+N^{4s_c-1}K)^{1/4}\xonorms{\vert\nabla\vert^{-1/4}u_{>N}}{4} 
		\xonorm{\nsc u_{\leq N}}{2}{4}
	\\ &\lesssim_u \eta (N^{4s_c-1}+K)^{3/4}\xonorms{\vert\nabla\vert^{-1/4}u_{>N}}{4}
	\\ &\lesssim_u\eta(N^{4s_c-1}+K)+\eta\xonorms{\vert\nabla\vert^{-1/4}u_{>N}}{4}^4,
	\end{align*}
which is acceptable.

Finally, we turn to \eqref{momb iii}. Using H\"older, Bernstein, Sobolev embedding, \eqref{assume1.1}, \eqref{properties of a}, \eqref{hc}, \eqref{lc}, \eqref{hc2}
	\begin{align*}
	&\eqref{momb iii}
	\\ &\lesssim \sup_{x\in\R^4}\norm{\vert u_{>N}\vert^2\,\nabla a(x-y)}_{L_t^4L_y^1}
		\xonorm{\vert\nabla\vert^{-1}u_{>N}\,\text{\O}(u^p)\Delta u_{\leq N}}{4/3}{1}
	\\ &\lesssim\xonorm{u_{>N}}{\infty}{2}\xonorm{u_{>N}}{4}{8/3}\sup_{x\in\R^4}\norm{\nabla a(x-y)}_{L_y^8}
		\\ &\quad\quad\quad\times\xonorms{\vert\nabla\vert^{-1}u_{>N}}{4}\xonorm{u}{\infty}{2p}^p\xonorm{\Delta u_{\leq N}}{2}{4}
	\\ &\lesssim_u\eta N^{-1/2-2s_c}(1+N^{4s_c-1}K)^{1/4}
		\\&\quad\quad\quad\times N^{-3/4}\xonorms{\vert\nabla\vert^{-1/4}u_{>N}}{4}N^{2-s_c}\xonorm{\nsc u_{\leq N}}{2}{4}
	\\ &\lesssim_u\eta^2 (N^{1-4s_c}+K)^{3/4}\xonorms{\vert\nabla\vert^{-1/4}u_{>N}}{4}
	\\ &\lesssim_u\eta(N^{1-4s_c}+K)+\eta\xonorms{\vert\nabla\vert^{-1/4}u_{>N}}{4}^4,
	\end{align*}
which is acceptable. This completes our treatment of \eqref{iiiiii}.

We now turn to \eqref{ivvvi}. Once again, we write $u_{>N}=\nabla\cdot\nabla\Delta^{-1}u_{>N}$ and integrate by parts. In this way, we see that to control \eqref{ivvvi}, it will suffice to control
	\begin{align}
	&\iiint \vert u_{>N}(t,y)\vert^2\vert\nabla\vert^{-1}u_{>N}(t,x)\nabla\text{\O}(u^p)(t,x)\Delta a(x-y)u_{\leq N}(t,x)\,dx\,dy\,dt\label{momb iv}
	\\ &+\iiint \vert u_{>N}(t,y)\vert^2\vert\nabla\vert^{-1}u_{>N}(t,x)\text{\O}(u^p)(t,x)\Delta a(x-y)\nabla u_{\leq N}(t,x)\,dx\,dy\,dt\label{momb v}
	\\ &+\iiint \vert u_{>N}(t,y)\vert^2\vert\nabla\vert^{-1}u_{>N}(t,x)\text{\O}(u^p)(t,x)\nabla\Delta a(x-y) u_{\leq N}(t,x)\,dx\,dy\,dt.\label{momb vi}
	\end{align}

We first turn to \eqref{momb iv}. We choose a parameter $r$ such that $\tfrac{7p-8}{2p-2}<r<4$. The lower bound guarantees that we may apply \eqref{hc} with $\eps=0$, while the upper bound guarantees $\nabla a\in L_x^r(\R^4)$ (cf. \eqref{properties of a}). Note also that as $s_c<3/2$ (that is, $p<4$) and $r<4$, we have that the exponent $\tfrac{8pr}{8p+4r-3pr}$ (which appears below) is finite. Then, using H\"older, Bernstein, Sobolev embedding, Young's inequality, \eqref{assume1.1}, \eqref{properties of a}, \eqref{hc}, and \eqref{lc}, we estimate
	\begin{align*}
	&\eqref{momb iv}
	\\ &\lesssim\sup_{x\in\R^4}\norm{\vert u_{>N}\vert^2\,\Delta a(x-y)}_{L_t^{r/2}L_y^{1}}\xonorm{\vert\nabla\vert^{-1} u_{>N}\nabla\text{\O}(u^p) u_{\leq N}}{\frac{r}{r-2}}{1}
	\\ &\lesssim\sup_{x\in\R^4}\norm{\Delta a(x-y)}_{L_y^r}\xonorm{u_{>N}}{r}{\frac{2r}{r-1}}^2\xonorms{\vert\nabla\vert^{-1}u_{>N}}{4}
	\\ &\quad\quad\quad\times\xonorm{\nabla (u^p)}{\infty}{4/3}\xonorm{u_{\leq N}}{\frac{4r}{3r-8}}{\infty}
	\\ &\lesssim_u N^{1-4/r-2s_c}(1+N^{4s_c-1}K)^{2/r}N^{-3/4}\xonorms{\vert\nabla\vert^{-1/4}u_{>N}}{4}
	\\ &\quad\quad\quad\times\xonorm{\nsc u}{\infty}{2}^p N^{1/2-s_c+4/r}\xonorm{u_{\leq N}}{\frac{4r}{3r-8}}{\frac{8pr}{8p+4r-3pr}}
	\\ &\lesssim_u N^{3/4-3s_c}(1+N^{4s_c-1}K)^{2/r}\xonorms{\vert\nabla\vert^{-1/4}u_{>N}}{4}
	\\ &\quad\quad\quad\times\xonorm{\nsc u_{\leq N}}{\frac{4r}{3r-8}}{\frac{8r}{r+8}}
	\\ &\lesssim_u \eta (N^{4s_c-1}+K)^{3/4}+\eta\xonorms{\vert\nabla\vert^{-1/4}u_{>N}}{4}^4,
	\end{align*}
which is acceptable. 

Next, we notice that we have already seen how to handle \eqref{momb v}, as it is identical to \eqref{momb ii}.

We turn to \eqref{momb vi}. We use H\"older, the Mikhlin multiplier theorem, Hardy--Littlewood--Sobolev, Bernstein, Sobolev embedding, Young's inequality, \eqref{assume1.1}, \eqref{properties of a}, \eqref{hc}, and \eqref{lc} to estimate

	\begin{align*}
	&\eqref{momb vi}
	\\ &\lesssim\xonorm{\vert x\vert^{-2}\ast\vert u_{>N}\vert^2}{2}{4}\xonorm{\vert\nabla\vert^{-1} u_{>N}\text{\O}(u^p)u_{\leq N}}{2}{4/3}
	\\ &\lesssim\xonorm{u_{>N}}{4}{8/3}^2 \xonorms{\vert\nabla\vert^{-1}u_{>N}}{4}
	\xonorm{u}{\infty}{2p}^p\xonorm{u_{\leq N}}{4}{\infty}
	\\ &\lesssim_u N^{-2s_c}(1+N^{4s_c-1}K)^{1/2}N^{-3/4}
		\xonorms{\vert\nabla\vert^{-1/4}u_{>N}}{4}
	N^{3/2-s_c}	
	\xonorm{u_{\leq N}}{4}{\frac{8p}{4-p}}
	\\ &\lesssim_u N^{3/4-3s_c}(1+N^{4s_c-1}K)^{1/2}\xonorms{\vert\nabla\vert^{-1/4}u_{>N}}{4}\xonorm{\nsc u_{\leq N}}{4}{8/3}
	\\ &\lesssim_u \eta(N^{4s_c-1}+K)^{3/4}\xonorms{\vert\nabla\vert^{-1/4}u_{>N}}{4}
	\\ &\lesssim_u\eta(N^{4s_c-1}+K)+\eta\xonorms{\vert\nabla\vert^{-1/4}u_{>N}}{4}^4,
	\end{align*}
which is acceptable. We have now finished estimating \eqref{ivvvi}, which in turn completes our estimation of the contribution of term $II$ in \eqref{momb decomp}.

Finally, we turn to estimating the contribution of term $III$ in \eqref{momb decomp}. We begin by writing
	\begin{align*}
	\{P_{\leq N}(F(u)),u_{>N}\}_P&=\nabla\text{\O}\big( P_{\leq N}(F(u))u_{>N}\big)
	\\ &\quad+\text{\O}\big(P_{\leq N}(F(u))\nabla u_{>N}\big).
	\end{align*}
Integrating by parts in both terms, we find that it will suffice to control the terms
	\begin{align}
	&\iiint\vert u_{>N}(t,y)\vert^2\text{\O}(P_{\leq N}(F(u))u_{>N})(t,x)\Delta a(x-y)\,dx\,dy\,dt \label{momb a}
	\\ &+\iiint\vert u_{>N}(t,y)\vert^2u_{>N}(t,x)\nabla P_{\leq N}(F(u))(t,x)\nabla a(x-y)\,dx\,dy\,dt.\label{momb b}
	\end{align}

To treat \eqref{momb a}, we further decompose
	\begin{align}
	&\eqref{momb a}\nonumber
	\\&=\iiint\vert u_{>N}(t,y)\vert^2 \text{\O}\big(P_{\leq N}(F(u_{\leq N})u_{>N}\big)(t,x)\Delta a(x-y)\,dx\,dy\,dt \label{momb a1}
	\\&+\iiint\vert u_{>N}(t,y)\vert^2\text{\O}\big(P_{\leq N}(F(u_{>N})u_{>N}\big)(t,x)\Delta a(x-y)\,dx\,dy\,dt\label{momb a2}
	\\&+\iiint\vert u_{>N}(t,y)\vert^2\text{\O}\big(P_{\leq N}(F(u)-F(u_{>N})-F(u_{\leq N}))u_{>N}\big)(t,x)\nonumber
	\\ &\quad\quad\quad\quad\times\Delta a(x-y)\,dx\,dy\,dt\label{momb a3}.
	\end{align}

We have already seen how to handle \eqref{momb a1} and \eqref{momb a3} when dealing with \eqref{momb I}. Indeed, we only need to note that 
	$$P_{\leq N}(F(u_{\leq N})u_{>N}\lesssim\text{\O}(u_{\leq N}u_{>N}u^p)$$
and
	$$P_{\leq N}(F(u)-F(u_{>N})-F(u_{\leq N}))u_{>N}
		\lesssim\text{\O}(u_{\leq N}u_{>N}u^p).$$ 

To treat \eqref{momb a2}, we write $u_{>N}=\nabla\cdot\nabla\Delta^{-1}u_{>N}$ and integrate by parts. In this way, we see that we need to estimate the terms
	\begin{align}
	&\iiint \vert u_{>N}(t,y)\vert^2\vert\nabla\vert^{-1}u_{>N}(t,x)\nabla P_{\leq N}(F(u_{>N}))(t,x)\Delta a(x-y)\,dx\,dy\,dt\label{a2i}
	\\ &\!\!+\!\!\iiint \vert u_{>N}(t,y)\vert^2\vert\nabla\vert^{-1}u_{>N}(t,x)P_{\leq N}(F(u_{>N}))(t,x)\nabla\Delta a(x-y)\,dx\,dy\,dt\label{a2ii}.
	\end{align}
	
To estimate \eqref{a2i}, we will use H\"older, Hardy--Littlewood--Sobolev, Bernstein, Lemma~\ref{pde}, the fractional chain rule, Sobolev embedding, \eqref{assume1.1}, \eqref{properties of a}, \eqref{hc}, and \eqref{hc2}. We find
	\begin{align*}
	&\eqref{a2i}
	\\&\lesssim\xonorm{\vert x\vert^{-1}\ast\vert u_{>N}\vert^2}{10/3}{10}\xonorm{\vert\nabla\vert^{-1}u_{>N}\nabla P_{\leq N}(F(u_{>N}))}{10/7}{10/9}
	\\&\lesssim\xonorm{\vert u_{>N}\vert^2}{10/3}{20/17}\xonorms{\vert\nabla\vert^{-1}u_{>N}}{4}\xonorm{\nabla P_{\leq N}(F(u_{>N}))}{20/9}{20/13}
	\\ &\lesssim\xonorm{u_{>N}}{\infty}{2}\xonorm{u_{>N}}{10/3}{20/7}
		N^{-3/4}\xonorms{\vert\nabla\vert^{-1/4}u_{>N}}{4}
	\\ &\quad\quad\quad\times N^{1+s_c}\xonorm{\vert\nabla\vert^{-1/4}F(u_{>N})}{20/9}{\frac{80p}{87p-40}}
	\\ &\lesssim_u \eta N^{1/4-s_c}(1+N^{4s_c-1}K)^{3/10}
		\xonorms{\vert\nabla\vert^{-1/4}u_{>N}}{4}^2
			\xonorm{\vert\nabla\vert^{1/4}\vert u_{>N}\vert^p}{5}{\frac{10p}{9p-5}}
	\\ &\lesssim_u \eta N^{1/4-s_c}(1+N^{4s_c-1}K)^{3/10}
		\xonorms{\vert\nabla\vert^{-1/4}u_{>N}}{4}^2
	\\ &\quad\quad\quad\times\xonorm{u}{\infty}{2p}^{p-1}\xonorm{\vert\nabla\vert^{1/4}u_{>N}}{5}{5/2}
	\\ &\lesssim_u \eta N^{1/2-2s_c}(1+N^{4s_c-1}K)^{1/2}\xonorms{\vert\nabla\vert^{-1/4}u_{>N}}{4}^2
	\\ &\lesssim_u \eta(N^{1-4s_c}+K)+\eta\xonorms{\vert\nabla\vert^{-1/4}u_{>N}}{4}^4,
	\end{align*}
which is acceptable. 

To estimate \eqref{a2ii}, we will estimate similarly. This time, we will use H\"older, the Mikhlin multiplier theorem, Hardy--Littlewood--Sobolev, Bernstein, Lemma~\ref{pde}, Sobolev embedding, the fractional chain rule, \eqref{assume1.1}, \eqref{properties of a}, \eqref{hc}, and \eqref{hc2}. We find
	\begin{align*}
	&\eqref{a2ii}
	\\ &\lesssim\xonorm{\vert x\vert^{-2}\ast\vert u_{>N}\vert^2}{10/3}{20/7}
		\xonorm{\vert\nabla\vert^{-1}u_{>N} P_{\leq N}(F(u_{>N}))}{10/7}{20/13}
	\\ &\lesssim\xonorm{\vert u_{>N}\vert^2}{10/3}{20/17}\xonorms{\vert\nabla\vert^{-1}u_{>N}}{4}\xonorm{P_{\leq N}(F(u_{>N}))}{20/9}{5/2}
	\\ &\lesssim	\xonorm{u_{>N}}{\infty}{2}\xonorm{u_{>N}}{10/3}{20/7}
		N^{-3/4}\xonorms{\vert\nabla\vert^{-1/4}u_{>N}}{4}
	\\&\quad\quad\quad\times	N^{1+s_c}\xonorm{\vert\nabla\vert^{-1/4}(F(u_{>N}))}{20/9}{\frac{80p}{87p-40}}
	\\&\lesssim_u \eta(N^{4s_c-1}+K)+\eta\xonorms{\vert\nabla\vert^{-1/4}u_{>N}}{4}^4
	\end{align*}
exactly as above. Thus \eqref{a2ii} is acceptable is as well. This completes the treatment of \eqref{momb a2}. We have now finished the estimation of \eqref{momb a}. 

We now turn to \eqref{momb b}. Again, we begin by writing $u_{>N}=\nabla\cdot\nabla\Delta^{-1}u_{>N}$ and integrating by parts. In this way, we see that to estimate \eqref{momb b} it will suffice to control the terms
	\begin{align}
	&\iiint \vert u_{>N}(t,y)\vert^2 \vert\nabla\vert^{-1}u_{>N}(t,x)\nabla P_{\leq N}(F(u))(t,x)\Delta a(x-y)\,dx\,dy\,dt \label{momb b one}
	\\ &+\iiint \vert u_{>N}(t,y)\vert^2\vert\nabla\vert^{-1}u_{>N}(t,x)\Delta P_{\leq N}(F(u))(t,x)\nabla a(x-y)\,dx\,dy\,dt\label{momb b two}.
	\end{align}
	
We first turn to \eqref{momb b one}. We decompose the nonlinearity to write
	\begin{align}
	\nonumber &\eqref{momb b one}
	\\ &=\!\!\iiint\vert u_{>N}(t,y)\vert^2 \vert\nabla\vert^{-1}u_{>N}(t,x)\nabla P_{\leq N}(F(u_{\leq N}))(t,x)\Delta a(x-y)\,dx\,dy\,dt \label{b2}
	\\ &+\!\!\iiint\vert u_{>N}(t,y)\vert^2 \vert\nabla\vert^{-1}u_{>N}(t,x)\nabla P_{\leq N}(F(u_{>N}))(t,x)\Delta a(x-y)\,dx\,dy\,dt \label{b1}
	\\ &+\!\!\iiint\vert u_{>N}(t,y)\vert^2 \vert\nabla\vert^{-1}u_{>N}(t,x)\nabla P_{\leq N}(F(u)-F(u_{\leq N})-F(u_{>N}))(t,x)\nonumber
	\\ &\quad\quad\quad\times\Delta a(x-y)\,dx\,dy\,dt. \label{b3}
	\end{align}
	
To estimate \eqref{b2}, we will use H\"older, Hardy--Littlewood--Sobolev, Bernstein, Sobolev embedding, Young's inequality, \eqref{hc}, \eqref{assume1.1}, \eqref{lc}, and\eqref{hc2}. We find
	\begin{align*}
	&\eqref{b2}
	\\&\lesssim\xonorm{\vert x\vert^{-1}\ast \vert u_{>N}\vert^2}{4}{8}\xonorm{\vert\nabla\vert^{-1}u_{>N}\nabla P_{\leq N}(F(u_{\leq N}))}{4/3}{8/7}
	\\ &\lesssim\xonorm{u_{>N}}{\infty}{2}\xonorm{u_{>N}}{4}{8/3}\xonorms{\vert\nabla\vert^{-1}u_{>N}}{4}\xonorm{\nabla P_{\leq N}(F(u_{\leq N}))}{2}{8/5}
	\\ &\lesssim_u \eta N^{-2s_c}(1+N^{4s_c-1}K)^{1/4} N^{-3/4}
		\xonorms{\vert\nabla\vert^{-1/4}u_{>N}}{4} 
	 N^{3/2-s_c}\xonorm{\nabla F(u_{\leq N})}{2}{\frac{2p}{p+1}}
	\\ &\lesssim_u \eta N^{3/4-3s_c}(1+N^{4s_c-1}K)^{1/4}\xonorms{\vert\nabla\vert^{-1/4}u_{>N}}{4}\xonorm{u}{\infty}{2p}^p\xonorm{\nabla u_{\leq N}}{2}{2p}
	\\ &\lesssim_u \eta N^{3/4-s_c}(1+N^{4s_c-1}K)^{1/4}\xonorms{\vert\nabla\vert^{-1/4}u_{>N}}{4}\xonorm{\nsc u_{\leq N}}{2}{4}
	\\ &\lesssim_u \eta^2 (N^{1-4s_c}+K)^{3/4}\xonorms{\vert\nabla\vert^{-1/4}u_{>N}}{4}
	\\ &\lesssim_u \eta(N^{1-4s_c}+K)+\eta\xonorms{\vert\nabla\vert^{-1/4}u_{>N}}{4}^4,
	\end{align*}
which is acceptable.

We next note that have already encountered \eqref{b1} before in the form of \eqref{a2i}; hence we will move on to \eqref{b3}. 

To estimate \eqref{b3}, we will use H\"older, Hardy--Littlewood--Sobolev, Bernstein, Sobolev embedding, Young's inequality, \eqref{assume1.1}, \eqref{hc}, \eqref{lc}, \eqref{hc2}, and \eqref{pw1}. Writing 
		$G=F(u)-F(u_{\leq N})-F(u_{>N}),$
we find
	\begin{align*}
	&\eqref{b3}
	\\&\lesssim\xonorm{\vert x\vert^{-1}\ast\vert u_{>N}\vert^2}{4}{8}\xonorm{\vert\nabla\vert^{-1}u_{>N}\nabla P_{\leq N}(G)}{4/3}{8/7}
	\\ &\lesssim\xonorm{u_{>N}}{\infty}{2}\xonorm{u_{>N}}{4}{8/3}
		\xonorms{\vert\nabla\vert^{-1}u_{>N}}{4}\xonorm{\nabla P_{\leq N}(G)}{2}{8/5}
	\\ &\lesssim_u \eta N^{-2s_c}(1+N^{4s_c-1}K)^{1/4} N^{-3/4}\xonorms{\vert\nabla\vert^{-1/4}u_{>N}}{4}
	\\&\quad\quad\quad\times N^{3/2}\xonorm{\text{\O}(u_{\leq N}u_{>N}u^{p-1})}{2}{4/3}
	\\ &\lesssim_u \eta N^{3/4-2s_c}(1+N^{4s_c-1}K)^{1/4}\xonorms{\vert\nabla\vert^{-1/4}u_{>N}}{4}
	\\&\quad\quad\quad\times\xonorm{u_{\leq N}}{4}{\frac{8p}{4-p}}\xonorm{u_{>N}}{4}{8/3}\xonorm{u}{\infty}{2p}^{p-1}
	\\ &\lesssim_u \eta N^{3/4-3s_c}(1+N^{4s_c-1}K)^{1/2}\xonorms{\vert\nabla\vert^{-1/4}u_{>N}}{4}\xonorm{\nsc u_{\leq N}}{4}{8/3}
	\\ &\lesssim_u \eta^2 (N^{4s_c-1}+K)^{3/4}\xonorms{\vert\nabla\vert^{-1/4}u_{>N}}{4}
	\\ &\lesssim_u \eta(N^{4s_c-1}+K)+\eta\xonorms{\vert\nabla\vert^{-1/4}u_{>N}}{4}^4,
	\end{align*}
which is acceptable.

We next turn to \eqref{momb b two}. Again, we decompose the nonlinearity to write
	\begin{align}
	&\eqref{momb b two}\nonumber
	\\ &=\!\!\iiint \vert u_{>N}(t,y)\vert^2 \vert\nabla\vert^{-1}u_{>N}(t,x)
		\Delta P_{\leq N}(F(u_{\leq N}))(t,x)\nabla a(x-y)\,dx\,dy\,dt\label{b5}
	\\ &+\!\!\iiint\vert u_{>N}(t,y)\vert^2\vert\nabla\vert^{-1}u_{>N}(t,x)
		\Delta P_{\leq N}(F(u_{>N}))(t,x)\nabla a(x-y)\,dx\,dy\,dt \label{b4}
	\\ &+\!\!\iiint\vert u_{>N}(t,y)\vert^2\vert\nabla\vert^{-1}u_{>N}(t,x)
		\Delta P_{\leq N}(F(u)-F(u_{\leq N})-F(u_{>N}))(t,x)\nonumber
		\\ &\quad\quad\quad\quad\times\nabla a(x-y)\,dx\,dy\,dt\label{b6}
	\end{align}

To estimate \eqref{b5}, we will use H\"older, Bernstein, Sobolev embedding, Young's inequality, \eqref{assume1.1}, \eqref{properties of a}, \eqref{hc}, \eqref{lc}, \eqref{hc2} . We find	
	\begin{align*}
	&\eqref{b5}
	\\&\lesssim\sup_{x\in\R^4}\norm{\vert u_{>N}\vert^2\,\nabla a(x-y)}_{L_t^4L_y^1}\xonorm{\vert\nabla\vert^{-1} u_{>N}\Delta P_{\leq N}(F(u_{\leq N})}{4/3}{1}
	\\ &\lesssim\sup_{x\in\R^4}\norm{\nabla a(x-y)}_{L_y^8}\xonorm{u_{>N}}{\infty}{2}\xonorm{u_{>N}}{4}{8/3}\xonorm{\vert\nabla\vert^{-1}u_{>N}}{4}{4}
	\\&\quad\quad\quad\times\xonorm{\Delta P_{\leq N}(F(u_{\leq N}))}{2}{4/3}
	\\ &\lesssim_u \eta N^{-1/2-2s_c}(1+N^{4s_c-1}K)^{1/4}N^{-3/4}\xonorms{\vert\nabla\vert^{-1/4}u_{>N}}{4}
	\\&	\quad\quad\quad\times N^{2-s_c}\xonorm{\nsc F(u_{\leq N})}{2}{4/3}
	\\ &\lesssim_u \eta N^{3/4-3s_c}(1+N^{4s_c-1}K)^{1/4}\xonorms{\vert\nabla\vert^{-1/4}u_{>N}}{4}\xonorm{u}{\infty}{2p}^p\xonorm{\nsc u_{\leq N}}{2}{4}
	\\ &\lesssim_u \eta^2 (N^{4s_c-1}+K)^{3/4}\xonorms{\vert\nabla\vert^{-1/4}u_{>N}}{4}
	\\ &\lesssim_u \eta(N^{4s_c-1}+K)+\eta\xonorms{\vert\nabla\vert^{-1/4}u_{>N}}{4}^4,
	\end{align*}
which is acceptable.

To estimate \eqref{b4}, we will use H\"older, Young's convolution inequality, Bernstein, Lemma~\ref{pde}, Sobolev embedding, Young's inequality, \eqref{assume1.1}, \eqref{properties of a}, \eqref{hc}, and \eqref{hc2}. We find
	\begin{align*}
	&\eqref{b4}
	\\&\lesssim \xonorm{\nabla a\ast \vert u_{>N}\vert^2}{10/3}{4/3}
		\xonorm{\vert\nabla\vert^{-1}u_{>N}\Delta P_{\leq N}(F(u_{\leq N})}{10/7}{4}
	\\ &\lesssim\norm{\nabla a}_{L_x^{10/9}}\xonorm{u_{>N}}{\infty}{2}\xonorm{u_{>N}}{10/3}{20/7}\xonorms{\vert\nabla\vert^{-1}u_{>N}}{4}\xonorm{\Delta P_{\leq N}F(u_{>N})}{20/9}{\infty}
	\\ &\lesssim \eta N^{-18/5-2s_c}(1+N^{4s_c-1}K)^{3/10} N^{-3/4}\xonorms{\vert\nabla\vert^{-1/4}u_{>N}}{4}
	\\ &\quad\quad\quad\times N^{80/7-2/p}\xonorm{\vert\nabla\vert^{-1/4}F(u_{>N})}{20/9}{\frac{80p}{87p-40}}
	\\ &\lesssim_u \eta N^{1/4-s_c} (1+N^{4s_c-1}K)^{3/10}\xonorms{\vert\nabla\vert^{-1/4}u_{>N}}{4}^2\xonorm{\vert\nabla\vert^{1/4}\vert u_{>N}\vert^p}{5}{\frac{10p}{9p-5}}
	\\ &\lesssim_u \eta N^{1/4-s_c}(1+N^{4s_c-1}K)^{3/10}\xonorms{\vert\nabla\vert^{-1/4}u_{>N}}{4}^2\xonorm{u}{\infty}{2p}^{p-1}\xonorm{\vert\nabla\vert^{1/4}u_{>N}}{5}{5/2}
	\\ &\lesssim_u \eta (N^{1-4s_c}+K)^{1/2}\xonorms{\vert\nabla\vert^{-1/4}u_{>N}}{4}^2
	\\ &\lesssim_u \eta(N^{1-4s_c}+K)+\eta\xonorms{\vert\nabla\vert^{-1/4}u_{>N}}{4}^4,
	\end{align*}
which is acceptable. 

To estimate \eqref{b6}, we will use H\"older, Young's convolution inequality, Bernstein, Sobolev embedding, \eqref{properties of a}, \eqref{hc}, \eqref{lc}, \eqref{hc2}, \eqref{pw1}. Writing $G=F(u)-F(u_{\leq N})-F(u_{>N}),$ we find
	\begin{align*}
	&\eqref{b6}
	\\&\lesssim\xonorm{\nabla a\ast \vert u_{>N}\vert^2}{4}{\frac{2p}{p-1}}
		\xonorm{\vert\nabla\vert^{-1}u_{>N}\Delta P_{\leq N}(G)}{4/3}{\frac{2p}{p+1}}
	\\ &\lesssim \norm{\nabla a}_{L_x^{\frac{8p}{5p-4}}}\xonorm{u_{>N}}{\infty}{2}\xonorm{u_{>N}}{4}{8/3}\xonorms{\vert\nabla\vert^{-1}u_{>N}}{4}\xonorm{\Delta P_{\leq N}(G)}{2}{\frac{4p}{2+p}}
	\\ &\lesssim_u \eta N^{-1/2-3s_c}(1+N^{4s_c-1}K)^{1/4} N^{-3/4}\xonorms{\vert\nabla\vert^{-1/4}u_{>N}}{4}
	\\ &\quad\quad\quad\times N^{2+s_c}\xonorm{\text{\O}(u_{\leq N}u_{>N}u^p)}{2}{4/3}
	\\ &\lesssim_u \eta N^{3/4-2s_c}(1+N^{4s_c-1}K)^{1/4}\xonorms{\vert\nabla\vert^{-1/4}u_{>N}}{4}
	\\ &\quad\quad\quad\times\xonorm{u_{\leq N}}{4}{\frac{8p}{4-p}}\xonorm{u_{>N}}{4}{8/3}\xonorm{u}{\infty}{2p}^{p-1}
	\\ &\lesssim_u \eta N^{3/4-3s_c}(1+N^{4s_c-1}K)^{1/2}\xonorms{\vert\nabla\vert^{-1/4}u_{>N}}{4}\xonorm{\nsc u_{\leq N}}{4}{8/3}
	\\ &\lesssim_u \eta^2 (N^{4s_c-1}+K)^{3/4}\xonorms{\vert\nabla\vert^{-1/4}u_{>N}}{4}
	\\ &\lesssim_u \eta(N^{4s_c-1}+K)+\eta\xonorms{\vert\nabla\vert^{-1/4}u_{>N}}{4}^4,
	\end{align*}
which is acceptable. This completes our treatment of \eqref{momb b two}.

With \eqref{momb b one} and \eqref{momb b two} estimated, we have completed our estimation of \eqref{momb b}. Thus, we have finished our estimation of the contribution of term $III$ to \eqref{momb decomp}, which completes the proof of the Lemma~\ref{Pb terms}.
\end{proof}

%%%% potential energy term...%%%%
The final term we need to estimate is the potential energy term appearing in \eqref{Pb}. Similar to the term appearing in Lemma~\ref{nssc}, this term can be dangerous in the energy-critical (or energy-subcritical) regime; indeed, in that setting one cannot expect smallness from $\|u_{>N}\|_{L_t^\infty L_x^{p+2}}$ (as it scales like the $\dot{H}_x^1$ norm).  For this reason, this term also necessitates a carefully designed weight $a$ in \cite{KV3D}. In the energy-supercritical setting, however, we can use interpolation to control the $L_x^{p+2}$ norm by the $L_x^{2p}$ and $L_x^2$ norm, and hence we can exhibit smallness in this term (cf. \eqref{assume1.1}, Sobolev embedding, and \eqref{hc2}). In particular, we have the following

\begin{lemma}[Potential energy term]\label{potential energy} Let $\eta>0$ and take $N_0=N_0(\eta)$ as in Corollary~\ref{a priori bounds}. Let $N\leq N_0$. We estimate 
		$$\int_I\iint_{\R^4\times\R^4}\vert u_{>N}(t,y)\vert^2\Delta a(x-y)\vert u_{>N}(t,x)\vert^{p+2}\,dx\,dy\,dt$$
in two pieces. First,
	$$\int_I\iint_{\vert x-y\vert\leq N^{-1}}\vert u_{>N}(t,y)\vert^2\Delta a(x-y)\vert u_{>N}(t,x)\vert^{p+2}\,dx\,dy\,dt\geq 0,$$
as $\Delta a\geq 0$ in this region. Second,
	$$\int_I\iint_{\vert x-y\vert>N^{-1}}\vert u_{>N}(t,y)\vert^2\Delta a(x-y)\vert u_{>N}(t,x)\vert^{p+2}\,dx\,dy\,dt\lesssim_u \eta^{\frac{p-2}{p-1}}(N^{1-4s_c}+K).$$ 

\end{lemma}
\begin{proof}
We first notice that by \eqref{properties of a}, we have $a(x-y)=\vert x-y\vert$ for $\vert x-y\vert\leq N^{-1}$. Thus the first claim follows immediately.

Next, using \eqref{properties of a}, we can see that to prove the second claim it will suffice to show
	\begin{align}
	\int_I\!\iint_{N^{-1}<\vert x-y\vert\leq 2N^{-1}}\!\!\frac{\vert u_{>N}(t,x)\vert^2\vert u_{>N}(t,y)\vert^2}{\vert x-y\vert}\,dx\,dy\,dt
		\lesssim \eta^{\frac{p-2}{p-1}}(N^{1-4s_c}\!+K).
		\label{pe}
		\end{align}
Using interpolation, Sobolev embedding, \eqref{assume1.1}, \eqref{hc2}, and \eqref{mass control}, we find		
	\begin{align*}
	\LHS\eqref{pe}&\lesssim N\xonorm{u_{>N}}{\infty}{p+2}^{p+2}\int_I\sup_y\int_{\vert x-y\vert\leq 2N^{-1}}\vert u_{>N}(t,x)\vert^2\,dx\,dt
	\\&\lesssim_u N\xonorm{u_{>N}}{\infty}{2}^{\frac{p-2}{p-1}}\xonorm{u_{>N}}{\infty}{2p}^{\frac{p^2}{p-1}}
		N^{-2-2s_c}(1+N^{4s_c-1}K)
	\\&\lesssim_u (\eta N^{-s_c})^{\frac{p-2}{p-1}}N^{-1-2s_c}(1+N^{4s_c-1}K)
	\\&\lesssim_u \eta^{\frac{p-2}{p-1}}(N^{1-4s_c}+K),
	\end{align*}
as needed.
\end{proof}

Having addressed all of the terms in Proposition~\ref{morawetz identity}, we can finally turn to the

\begin{proof}[Proof of Theorem~\ref{uppbd}] Let $\eta>0$ and take $N\leq N_0(\eta)$, where $N_0$ is as in Corollary~\eqref{a priori bounds}. Recalling the definition of $M(t)$ in \eqref{morawetz action} and choosing $\varphi=u_{>N}$, we first note that
	$$\norm{M(t)}_{L_t^\infty(I)}\lesssim\xonorm{u_{>N}}{\infty}{2}^3\xonorm{\nabla u_{>N}}{\infty}{2}\lesssim_u \eta^4 N^{1-4s_c}.$$
We can thus use the identity for $\partial_t M$ in Proposition~\ref{morawetz identity}, the fundamental theorem of calculus, and all of the lemmas in this section to deduce
	\begin{align*}
	\int_I&\iint_{\R^4\times\R^4}\frac{\vert u_{>N}(t,y)\vert^2\vert u_{>N}(t,x)\vert^2}{\vert x-y\vert^3}\,dx\,dy\,dt
	\\ &\lesssim_u\eta(N^{1-4s_c}+K)+\eta\int_I\iint_{\R^4\times\R^4}\frac{\vert u_{>N}(t,y)\vert^2\vert u_{>N}(t,x)\vert^2}{\vert x-y\vert^3}\,dx\,dy\,dt.
	\end{align*}

We would now like to choose $\eta$ small enough that we can absorb the second term on the right-hand side of the inequality above into the left-hand side (thus completing the proof of Theorem~\ref{uppbd}). To justify this step, however, we need to show that this term is finite to begin with. This is indeed the case: by H\"older, Bernstein, Lemma~\ref{fspssn}, Lemma~\ref{pde}, Lemma~\ref{bootstrap}, and  \eqref{assume1.1}, we have
	\begin{align*}
	\int_I\int_{\R^4\times\R^4}&\frac{\vert u_{>N}(t,y)\vert^2\vert u_{>N}(t,x)\vert^2}{\vert x-y\vert^3}\,dx\,dy\,dt
	\\&\sim\xonorms{\vert\nabla\vert^{-1/2} \vert u_{>N}\vert^2}{2}^2
	\\&\lesssim\xonorm{\vert\nabla\vert^{-1/2}u_{>N}}{\infty}{8/3}^2\xonorm{\vert\nabla\vert^{1/2}u_{>N}}{2}{4}^2
	\\&\lesssim\xonorm{u_{>N}}{\infty}{2}^2N^{1-2s_c}\xonorm{\vert\nabla\vert^{s_c} u}{2}{4}^2
	\\&\lesssim_u N^{1-4s_c}\big(1+\smallint_I N(t)^2\,dt\big).
	\end{align*} 
This completes the proof of Theorem~\ref{uppbd}.
\end{proof}

\section{The quasi-soliton scenario}\label{qs section}
In this section, we preclude the existence of quasi-soliton soutions, that is, solutions as in Theorem~\ref{threenew} such that $\int_0^{T_{max}}N(t)^{3-4s_c}\,dt=\infty.$ The proof will rely primarily on the frequency-localized interaction Morawetz inequality established in the previous section.
\begin{theorem}[No quasi-solitons]\label{no quasi-solitons}
Let $1<s_c<3/2$. Then there are no almost periodic solutions $u:[0,T_{max})\times\R^4\to\C$ to \eqref{equ1.1} with $N(t)\equiv N_k\geq 1$ on each characteristic subinterval $J_k\subset[0,T_{max})$ that satisfy
			\begin{equation}\label{rfcmd12}
			\nonumber
			\xnorms{u}{3p}{[0,T_{max})}=\infty
			\end{equation}
and
			\begin{equation}\label{rfkjfs}
			\int_0^{T_{max}}N(t)^{3-4s_c}\,dt=\infty.
			\end{equation}
\end{theorem}
\begin{proof} We argue by contradiction. Suppose that $u:[0,T_{max})\times\R^4\to\C$ were such a solution. We first claim that we have the following

\begin{lemma}[Lower bound]\label{lowbound} There exists $N_1>0$ such that for any $N\leq N_1,$
\begin{equation}\label{lowbd}
\int_I\iint_{\R^4\times\R^4}\frac{\vert u_{>
N}(t,x)\vert^2\vert u_{> N}(t,y)\vert^2}{\vert x-y\vert^3}dx\,dy\,dt \gtrsim_u \int_I N(t)^{3-4s_c}\,dt
\end{equation}
for any compact interval $I\subset[0,T_{max})$ that is a contiguous union of characteristic subintervals $J_k.$ 
\end{lemma}

We will prove Lemma~\ref{lowbound} below; for now, let us take it for granted and use it to complete the proof of Theorem~\ref{no quasi-solitons}. 

Let $\eta>0$ and choose $N_0=N_0(\eta)$ as in Theorem~\ref{uppbd}. Let $I\subset[0,T_{max})$ be a compact time interval that is a contiguous union of characteristic subintervals $J_k$. Combining \eqref{limemd} and \eqref{lowbd}, we find that for $N\leq \min\{N_0,N_1\}$, we have
$$\int_I N(t)^{3-4s_c}\,dt\lesssim_u \eta\Big(N^{1-4s_c}+\int_I N(t)^{3-4s_c}\,dt\Big)$$
uniformly in $I$.

Choosing $\eta$ sufficiently small, we deduce
	$$\int_I N(t)^{3-4s_c}\,dt\lesssim_u N^{1-4s_c}$$
uniformly in $I.$ We can now contradict \eqref{rfkjfs} by taking $I$ to be sufficiently large inside $[0,T_{max}).$ This completes the proof of Theorem~\ref{no quasi-solitons}.
\end{proof}

It remains to establish Lemma~\ref{lowbound}. 
\begin{proof}[Proof of Lemma~\ref{lowbound}] The key will be to show that there exist $N_1>0$ and $C(u)>0$ such that
	\begin{equation}\label{qs lemma}
	\inf_{t\in[0,T_{max})}N(t)^{2s_c}\int_{\vert x-x(t)\vert\leq\frac{C(u)}{N(t)}}\vert u_{>N}(t,x)\vert^2\,dx\gtrsim_u 1
	\end{equation}
for all $N\leq N_1$. Indeed, with \eqref{qs lemma} in place, we can estimate
	\begin{align*}
	\int_I&\iint_{\R^4\times\R^4}\frac{\vert u_{> N}(t,x)\vert^2\vert u_{> N}(t,y)\vert^2}{\vert x-y\vert^3}dx\,dy\,dt
	\\ &\gtrsim \int_I\iint_{\vert x-y\vert\leq\frac{2C(u)}{N(t)}}\Big[\tfrac{N(t)}{2C(u)}\Big]^3\vert u_{>N}(t,x)\vert^2\vert u_{>N}(t,y)\vert^2\,dx\,dy\,dt
	\\ &\gtrsim_u\int_I N(t)^3\Bigg(\int_{\vert x-x(t)\vert\leq\frac{C(u)}{N(t)}} \vert u_{>N}(t,x)\vert^2\,dx\Bigg)^2\,dt
	\\ &\gtrsim_u \int_I N(t)^{3-4s_c}\,dt,
	\end{align*}
as needed. Thus it remains to establish \eqref{qs lemma}. 

We will first establish that for $C(u)$ sufficiently large, we have 
		\begin{equation}\label{stepping stone}
		\inf_{t\in[0,T_{max})}N(t)^{2s_c}\int_{\vert x-x(t)\vert\leq\frac{C(u)}{N(t)}}\vert u(t,x)\vert^{2}\,dx\gtrsim_u 1.
		\end{equation}
To this end, we let $\eta_0>0$ and (by almost periodicity) choose $C_0:=C(\eta_0)$ large enough that
		\begin{equation}\label{small0}
		\xonorm{\nsc u_{>C_0N(t)}}{\infty}{2}<\eta_0.
		\end{equation}	
Then using H\"older, Bernstein, and Sobolev embedding, we can estimate
		\begin{align*}
		\bigg\vert\int_{\vert x-x(t)\vert\leq\frac{C(u)}{N(t)}}& \vert u(t,x)\vert^2-\vert u_{\leq C_0N(t)}(t,x)\vert^2\,dx\bigg\vert
		\\ &\lesssim_u N(t)^{-s_c}\norm{u_{>C_0N(t)}(t)}_{L_x^2(\R^4)}\norm{u(t)}_{L_x^{2p}(\R^4)}
		\\ &\lesssim_u \eta_0 N(t)^{-2s_c}
		\end{align*}
for $t\in [0,T_{max})$. Thus, if we can show that for $C(u)$ sufficiently large, we have
		\begin{equation}
		\label{stepping stepping stone}
		\inf_{t\in[0,T_{max})}  N(t)^{2s_c}
		\int_{\vert x-x(t)\vert\leq\frac{C(u)}{N(t)}}\vert u_{\leq C_0N(t)}(t,x)\vert^2\,dx
		\gtrsim_u 1,
		\end{equation}
then we will have \eqref{stepping stone} by choosing $\eta_0=\eta_0(u)$ sufficiently small. 

To prove \eqref{stepping stepping stone}, we first choose $C(u)>0$ large enough that
		$$\inf_{t\in[0,T_{max})}\int_{\vert x-x(t)\vert\leq\frac{C(u)}{N(t)}}\vert u(t,x)\vert^{2p}\,dx\gtrsim_u 1$$
(cf. Remark~\ref{rem1.1}). We then use H\"older, Sobolev embedding, and \eqref{small0} to see
	\begin{align*}
		\bigg\vert\int_{\vert x-x(t)\vert\leq\frac{C(u)}{N(t)}} \vert u(t,x)\vert^{2p}-\vert u_{\leq C_0N(t)}\vert^{2p}\,dx\bigg\vert
		\lesssim\xonorm{u_{>C_0N(t)}}{\infty}{2p}\xonorm{u}{\infty}{2p}^{2p-1}
		\lesssim_u\eta_0
	\end{align*}
for $t\in[0,T_{max})$. Thus for $\eta_0=\eta_0(u)$ sufficiently small, we find 
		\begin{equation}\label{2p ulo}
		\inf_{t\in[0,T_{max})}\int_{\vert x-x(t)\vert\leq\frac{C(u)}{N(t)}}\vert u_{\leq C_0N(t)}(t,x)\vert^{2p}\,dx\gtrsim_u 1.
		\end{equation}
Finally, using H\"older and Bernstein, we can get
	\begin{align*}
	\int_{\vert x-x(t)\vert\leq\frac{C(u)}{N(t)}}&\vert u_{\leq C_0N(t)}(t,x)\vert^{2p}\,dx
	\\ &\lesssim \norm{u_{\leq C_0N(t)}(t)}_{L_x^\infty(\R^4)}^{2p-2}\int_{\vert x-x(t)\vert\leq\frac{C(u)}{N(t)}} \vert u_{\leq C_0N(t)}(t,x)\vert^2\,dx
	\\ &\lesssim_u N(t)^{2s_c}\norm{u(t)}_{L_x^{2p}(\R^4)}^{2p-2}\int_{\vert x-x(t)\vert\leq \frac{C(u)}{N(t)}}\vert u_{\leq C_0N(t)}(t,x)\vert^2\,dx
	\\ &\lesssim_u N(t)^{2s_c}\int_{\vert x-x(t)\vert\leq \frac{C(u)}{N(t)}}\vert u_{\leq C_0N(t)}(t,x)\vert^2\,dx.
	\end{align*}
Together with \eqref{2p ulo}, this implies \eqref{stepping stepping stone}, which in turn implies \eqref{stepping stone}.

With \eqref{stepping stone} in place, we are now in a position to establish \eqref{qs lemma} and complete the proof of Lemma~\ref{lowbound}. We let $\eta_1>0$ be a small parameter to be determined shortly. As $\inf_{t\in[0,T_{max})}N(t)\geq 1$, we may find $N_1=N_1(\eta_1)$ so that
		$$\xonorm{u_{\leq N}}{\infty}{2p}<\eta_1\quad\text{for }N\leq N_1.$$
We then use H\"older and Sobolev embedding to estimate 
	\begin{align*}
	\bigg\vert\int_{\vert x-x(t)\vert\leq\frac{C(u)}{N(t)}}\vert u(t,x)\vert^2-\vert u_{>N}(t,x)\vert^2\,dx\bigg\vert
		&\lesssim_u N(t)^{-2s_c}\xonorm{u_{\leq N}}{\infty}{2p}\xonorm{u}{\infty}{2p}
		\\&\lesssim_u\eta_1N(t)^{-2s_c}
	\end{align*}
for $t\in[0,T_{max})$ and $N\leq N_1.$ Thus, choosing $\eta_1=\eta_1(u)$ sufficiently small, we may use \eqref{stepping stone} to deduce \eqref{qs lemma}. This completes the proof of Lemma \ref{lowbound}. \end{proof}

\begin{center}

\end{center}

\end{document}